\documentclass[12pt,oneside,openany,article]{memoir}
\usepackage{mempatch}

\nouppercaseheads 
\usepackage[amsthm,thmmarks,hyperref]{ntheorem}
\usepackage[noTeX]{mmap}
\usepackage{etex}

\usepackage[english]{babel}
\usepackage[leqno]{mathtools}

\usepackage{mathrsfs}
\usepackage{upgreek}
\usepackage[compress,square,comma,numbers]{natbib}
\usepackage{hypernat} %required if hyperref is used as well

\usepackage{sfmath}

\usepackage{microtype}

\usepackage{subfig}
\usepackage{wrapfig}
\usepackage{array}

\usepackage{graphicx,color}
\definecolor{gray}{gray}{0}
\pagecolor{white}

\usepackage{hyperxmp}
\usepackage{xr-hyper}
\usepackage{nameref}
\usepackage[pdftex,bookmarks,pdfnewwindow,plainpages=false,unicode]{hyperref}

\usepackage{bookmark}

\usepackage{enumitem}

\usepackage{amssymb} % for \Subset, \lozenge, etc.

\hypersetup{
colorlinks=true,
linkcolor=black,
citecolor=black,
urlcolor=blue,
pdfauthor={Victor Ivrii},
pdftitle={Multidimensional magnetic Schr\"odinger operator. I.\\
Full-rank case},
pdfsubject={Sharp Spectral Asymptotics},
pdfkeywords={Microlocal Analysis,  Sharp Spectral Asymptotics, Schort loops},
bookmarksdepth={4}
}

\numberwithin{equation}{chapter}

\theoremstyle{plain}
\newtheorem{theorem}{Theorem}[section]

\newtheorem{proposition}[theorem]{Proposition}
\newtheorem{corollary}[theorem]{Corollary}

\theoremstyle{definition}

\newtheorem{definition}[theorem]{Definition}
\newtheorem{Problem}[theorem]{Problem}
\theoremstyle{remark}
\newtheorem{remark}[theorem]{Remark}
\newtheorem{example}[theorem]{Example}
\newtheorem{problem}[theorem]{Problem}

\DeclareMathAlphabet{\mathpzc}{OT1}{pzc}{m}{it}

\newcommand{\cA}{\mathcal{A}}
 \newcommand{\cB}{\mathcal{B}}
 \newcommand{\cC}{\mathcal{C}}
 \newcommand{\cD}{\mathcal{D}}

 \newcommand{\cG}{\mathcal{G}}

 \newcommand{\cK}{\mathcal{K}}
 \newcommand{\cL}{\mathcal{L}}
 
 \newcommand{\cN}{\mathcal{N}}

 \newcommand{\cR}{\mathcal{R}}
 \newcommand{\cQ}{\mathcal{Q}}
 
 \newcommand{\cT}{\mathcal{T}}

 \newcommand{\cW}{\mathcal{W}}

 \newcommand{\sC}{\mathscr{C}}

 \newcommand{\sF}{\mathscr{F}}
 
 \newcommand{\sH}{\mathscr{H}}

 \newcommand{\sL}{\mathscr{L}}

\newcommand{\corr}{{\mathsf{corr}}}

\newcommand{\D}{{\mathsf{D}}}

\newcommand{\MW}{{\mathsf{MW}}}
\newcommand{\N}{{\mathsf{N}}}

\newcommand{\R}{{\mathsf{R}}}

\newcommand{\T}{{\mathsf{T}}}

\newcommand{\W}{{{\mathsf{W}}}}
\newcommand{\w}{{\mathsf{w}}}

\newcommand{\const}{{\mathsf{const}}}
\newcommand{\dist}{{{\mathsf{dist}}}}

\newcommand{\new}{{{\mathsf{new}}}}

\newcommand{\bC}{{\mathbb{C}}}
\newcommand{\bH}{{\mathbb{H}}}
\newcommand{\bK}{{\mathbb{K}}}
\newcommand{\bL}{{\mathbb{L}}}

\newcommand{\bR}{{\mathbb{R}}}

\newcommand{\bZ}{{\mathbb{Z}}}

\newcommand{\fA}{{\mathfrak{A}}}

\newcommand{\fm}{{\mathfrak{m}}}
\newcommand{\fM}{{\mathfrak{M}}}
\newcommand{\fn}{{\mathfrak{n}}}
\newcommand{\fN}{{\mathfrak{N}}}

\newcommand{\fz}{{\mathfrak{z}}}

\def\1{\boldsymbol {|}}
\newcommand{\boldgamma}{{\boldsymbol{\gamma}}}

\newcommand{\boldlambda}{{\boldsymbol{\lambda}}}

\newcommand{\boldrho}{{\boldsymbol{\rho}}}

\newcommand{\boldtau}{{\boldsymbol{\tau}}}

\newcommand{\blangle}{{\boldsymbol{\langle}}}
\newcommand{\brangle}{{\boldsymbol{\rangle}}}

\newcommand{\Def}{\mathrel{\mathop:}=}

\newcommand{\Ad}{\operatorname{Ad}}

\newcommand{\Hess}{\operatorname{Hess}}
\renewcommand{\Im}{\operatorname{Im}}       % \Im is already frak I

\newcommand{\mes}{\operatorname{mes}}

\newcommand{\rank}{\operatorname{rank}}
\renewcommand{\Re}{\operatorname{Re}}       % \Re is already frak R

\newcommand{\Spec}{\operatorname{Spec}}

\newcommand{\supp}{\operatorname{supp}}

\newcommand{\Tr}{\operatorname{Tr}}

\newenvironment{claim}[1][{\textup{(\theequation)}}]{\refstepcounter{equation}\vglue10pt
\begin{trivlist}
\item[{\hskip\labelsep#1}]}{\vglue10pt\end{trivlist}}

\newenvironment{claim*}[1][{}]{\vglue10pt
\begin{trivlist}
\item[{\hskip\labelsep#1}]}{\vglue10pt\end{trivlist}}

\newenvironment{phantomequation}[1][]{\refstepcounter{equation}}{}
\newcounter{note}

\DeclareTextCommand{\textinfty}{PU}{\9042\036}

\DeclareTextCommand{\textge}{PU}{\9042\145}
\DeclareTextCommand{\textle}{PU}{\9042\144}
\DeclareTextCommand{\texthat}{PD1}{\136}

\begin{document}
\title{Multidimensional magnetic Schr\"odinger operator. I.\\
Full-rank case}
\author{Victor Ivrii}

\maketitle
{\abstract%
With derive sharp spectral asymptotics (with the remainder estimate
$O(\mu ^{-1}h^{1-d}+\mu ^{\frac{d} {2}-1}h^{1-\frac{d}{2}})$ for $d$-dimensional Schr\"odinger operator with a strong magnetic field; here $h$ and $\mu$ are Plank and binding constants respectively and magnetic intensity matrix has full rank at each point.

In comparison with version 1 of 4.5 year ago this version contains more results (we also study some degenerations), improvements and some minor corrections.
\endabstract}

\setcounter{chapter}{-1}

\chapter{Introduction}
\label{sect-19-1}

\section{Preface}
\label{sect-19-1-1}

In this Chapter we consider multidimensional Schr\"odinger operator
\begin{multline}
A=A_0+V(x),\qquad A_0=\sum_{j,k\le d}P_jg ^{jk}(x)P_k, \\
P_j=hD_j- \mu V_j(x),
\quad h\in (0,1],\ \mu \ge 1.
\label{19-1-1}
\end{multline}
It is characterized by the \emph{magnetic field intensity tensors\/}\index{magnetic field intensity tensors}  $(F_{jk})$ with
\begin{equation}
F_{jk}=\partial _kV_j-\partial _jV_k,
\label{19-1-2}%
\end{equation}
which is skew-symmetric $d\times d$-matrix, and $(F^j_p)=(g^{jk})(F_{kp})$ which is equivalent to the skew-symmetric matrix $(g^{jk})^{\frac{1}{2}}(F_{jk})(g^{jk})^{\frac{1}{2}}$.

Compare this with $2\D$ and $3\D$ cases when we could characterize intensity by a (pseudo)scalar $F$ and (pseudo)vector $\mathbf{F}$ respectively.

Then
\begin{claim}\label{19-1-3}%
All eigenvalues of $(F^j_k)$ are $\pm i f_p$ ($f_l>0$, $p=1,\ldots, r$) and $0$ of multiplicity $d-2r$ where $2r=\rank (F^j_k)$.
\end{claim}

In this Chapter we assume that the magnetic field intensity tensor has a full rank:
\begin{equation}
\rank (F^j_k )=2r= d \quad\text{and}\quad |(F^j_k )^{-1}|\le c_0
\label{19-1-4}
\end{equation}
and (under certain conditions) we derive sharp spectral asymptotics (with the remainder estimate $O\bigl(\mu ^{-1}h^{1-d}\bigr)$ as $\mu h\lesssim 1$; as $\mu h\gtrsim 1$ we consider corresponding Schr\"odinger-Pauli operator and derive remainder estimate $O\bigl(\mu ^{r-1}h^{1-r}\bigr)$. The typical (but different from the general case) example (already studied in Chapters~\ref{book_new-sect-13} and~\ref{book_new-sect-18} of \cite{futurebook}) is $2\D$ magnetic Schr\"odinger operator.

As usual, we consider operator in some domain or on some manifold $X$ with some boundary conditions, assuming that it is self-adjoint in $\sL^2(X)$ and denote by $e(x,y,\tau)$ Schwartz' kernel of it spectral projector.

As usual we assume that conditions (\ref{book_new-13-1-4}) and (\ref{book_new-13-1-5}) of \cite{futurebook} are fulfilled i.e. 
\begin{gather}
\epsilon _0\le \sum_{j,k} g^{jk}\eta _j\eta _k\cdot |\eta |^{-2}\le c
\quad \forall \eta \in \bR^d\setminus 0 \quad \forall x \in B(0,1)
\label{19-1-5}\\% \tag{6.1.3}
\shortintertext{and}
X \supset B(0,1).
\label{19-1-6}%\tag{6.1.4}
\end{gather}

\section{Canonical form}
\label{sect-19-1-2}

Recall (see e.g. Subsection~\ref{book_new-sect-13-1-2}  of  \cite{futurebook}) that if $X={\bR}^d$ and $g^{jk},F_{jk},V$ are constant then operator (\ref{19-1-1}) is unitary equivalent to
\begin{equation}
\sum_{1\le j\le r}  f_j(h^2D_j^2+\mu ^2 x_j^2) +
\sum_{r+1\le j\le d-2r} h^2D_j^2 +V
\label{19-1-7}
\end{equation}
and $e(x,x,\tau) =h^{-d}\cN_d^\MW (\tau)$ is defined by (\ref{book_new-13-1-9}) of  \cite{futurebook};  in particular,  under condition (\ref{19-1-4}) decomposition (\ref{19-1-7}) becomes
\begin{equation}
\sum_{1\le j\le r}  f_j(h^2D_j^2+\mu ^2 x_j^2)  +V
\label{19-1-8}
\end{equation}
and
\begin{multline}
h^{-d}\cN_d^\MW (\tau)\Def  \\
(2\pi )^{-r} \mu ^rh^{-r}
\sum _{\alpha \in \bZ^{+r}}
\uptheta \Bigl(\tau - \sum_j(2\alpha_j +1) f_j\mu h -V\Bigr)
f_1\cdots f_r\sqrt{g}
\label{19-1-9}
\end{multline}
with  $g=\det (g_{jk})$, $(g_{jk})=(g^{jk})^{-1}$ and Heaviside function $\uptheta$.

In particular,
\begin{claim}\label{19-1-10}
For the pilot-model operator under condition (\ref{19-1-4}) (and only in this case) the spectrum is pure-point (of infinite multiplicity) consisting of \emph{Landau levels\/}\index{Landau levels}
\begin{equation}
E_\alpha\Def \sum_j(2\alpha_j+1)f_j, \qquad\alpha \in \bZ^{+r}.
\label{19-1-11}
\end{equation}
\end{claim}
\vglue-10pt

Now in the general case (i.e. for $X\ne \bR^d$, and  variable $V$ and, may be, $g^{jk},F_{jk}$) we are  interested in asymptotics of the spatially mollified spectral i.e. 
\begin{equation}
\N^\MW _\psi \Def \int e(x,x,0) \psi (x)\,dx 
\label{19-1-12}
\end{equation}
as $h\to +0,\mu \to +\infty$ where $\psi $ is a fixed function, smooth and compactly supported in $X$. As $\mu h \ge C_0$ the result will be trivial and therefore instead of the original Schr\"odinger operator we will need to consider generalized Schr\"odinger-Pauli operator\index{Schr\"odinger-Pauli operator}
\begin{equation}
A-  \sum_j \fz_j f_j \mu h
\label{19-1-13}
\end{equation}
with constants $\fz_j\in \bR$.

The principal part of such asymptotics under some reasonable and obvious conditions is of magnitude $h^{-2r}$ as $\mu h\lesssim 1$ and  $\mu ^r h^{-r}$ as $\mu h\gtrsim 1$.

On the other hand, the remainder estimate is $O(\mu h^{1-2r}+\mu ^rh^{-r})$ unless one imposes some non-degeneracy assumption while under the strongest possible assumption one can expect  remainder estimate
$O\bigl( \mu ^{-1}h^{1-2r} +\mu ^{r-1} h^{1-r}\bigr)$.

\begin{remark}\label{rem-19-1-1}
However even if $g^{jk},F_{jk},V$ are constant the remainder estimate is as bad as $\mu h^{1-d}$ (for $\mu h\ll 1$) only if $f_1,\ldots, f_r$ are commensurable. Otherwise as $\mu h \lesssim 1$ the remainder estimate is between 
$O\bigl(\mu^r h^{-r}\bigr)$ and $O\bigl(\mu h^{1-2r}\bigr)$ depending on the non-commensurability of $f_1,\ldots, f_r$. 

Probably it is the same remainder estimate as an error in the ``equality'' 
\begin{equation}
h^{-d}\cN^\MW (0)\approx h^{-d}\cN^\W_\infty  (0)
\label{19-1-14}
\end{equation}
where $h^{-d}\cN^\W_\infty$ denotes  Weyl approximation with many terms. Considering the left-hand expression as an $r$-dimensional Riemannian sum, passing to  the integral and correcting it so if instead of $\uptheta$ there was a smooth function we would get an error $O(h^{-d}(\mu h)^\infty)$, we get a right-hand expression.

Also obviously all Hamiltonian trajectories are periodic if and only if $f_1,\ldots, f_r$ are commensurable.
\end{remark}

\begin{Problem}\label{Problem-19-1-2}
Explore error in (\ref{19-1-14})  when $f_1,\ldots,f_r$ are non-com\-men\-surable. It definitely would depend on how these numbers are non-commensurable (something related to their diophantine properties - as $r=2$ it is related to Liouville's exponent for $f_1/f_2$).
\end{Problem}

Recall from Section~\ref{book_new-sect-13-2} of \cite{futurebook} that in the smooth case for $d=2$ one can reduce operator in question to the canonical form
\begin{equation}
\sum _{m+k+j \ge 1} \mu^{2-2k-2m-j}h^j
a_{m,k,j}(x_2,\mu ^{-1}hD_2) \bigl(h^2D_1^2+\mu^2x_1^2\bigr)^m
\label{19-1-15}
\end{equation}
with $a_{1,0}=F\circ \Psi_0$, $a_{0,1}= V\circ \Psi_0$ and a certain diffeomorphism $\Psi_0:T\bR^1 \to \bR^2$. Then one can replace harmonic oscillator  $h^2D_1^2+\mu^2x_1^2$ by one of its eigenvalues
$(2\alpha +1)\mu h$ with $\alpha \in \bZ^+$; we ignore terms with $j\ge 1$.

So, for $\mu h\le \epsilon_0$ our operator looked like a family (with
$C_0(\mu h)^{-1}$ elements) of $1$-dimensional  $\mu ^{-1}h$-pseudo-differential operators with the principal symbols
$\bigl((2\alpha +1)\mu h F+V\bigr) \circ \Psi_0$.

Spectral asymptotics for such operators are very sensitive to a degeneracy of the symbol but  under the non-degeneracy condition $|b|+|\nabla b|\ge \epsilon_1$ with $b=a_{0,1}\lambda +a_{1,0}$ with $\lambda\in \bR^+$ this asymptotics has a principal part of the magnitude $\mu h^{-1}$ and  a remainder  $O(1)$. This non-degeneracy condition is equivalent to
$|\nabla V/F |\ge \epsilon_1$ and the final asymptotics has its principal part of magnitude $h^{-2}$ and a remainder $O(\mu^{-1}h^{-1})$\,\footnote{\label{foot-19-1} We could also consider a weaker non-degeneracy condition invoking second derivatives.}.

On the other hand, as $\mu h\ge \epsilon_0$ we had the family of no more than $C_0$ of $\mu^{-1}h$-pseudo-differential operators, then the non-degeneracy condition became  $|(2\alpha +1)\mu h + (V/F)| + |\nabla (V/F)|\ge \epsilon_1$ and the final asymptotics has its principal part of magnitude $\mu h ^{-1}$ and  a remainder $O(1)$\,\footref{foot-19-1}.

Non-smoothness prevented us from the complete canonical form but we had a ``poor man'' canonical form which was sufficient; see Section~\ref{book_new-sect-18-3} of \cite{futurebook}.

Multidimensional case is much more tricky. Under assumption (\ref{19-1-4}) one could expect a canonical form
\begin{equation}
\sum_{m\in \bZ^{+r}, k,j: |m|+k+j\ge 1}
\mu ^{2-2k-2|m|-j}h^j a_{m,k,j}(x'',\mu ^{-1}hD'') H^m
\label{19-1-16}
\end{equation}
with $H=(H_1,\ldots,H_r)$, harmonic oscillators $H_i=h^2D_j^2+\mu^2x_i^2$
$a_{(m,0,0)}=f_j\circ \Psi_0$ as $|m|=0$ with $m_i=\updelta_{ji}$ \underline{and}
$a_{(0,1,0)}=V\circ \Psi_0$ where now we consider diffeomorphism
$\Psi_0:T\bR^r=\bR^{2r}\to \bR^{2r}$ with $x''=(x_{r+1},\dots, x_{2r})$.

In this case for  $\mu h\le \epsilon_0$ one would get a family (with 
$C_0(\mu h)^{-r}$ elements) of  $r$-dimensional  
$\mu ^{-1}h$-pseudo-differential operators, and the principal part of asymptotics for them would be of magnitude $\mu ^r h^{-r}$,  and under  proper non-degeneracy assumption  remainder $O(\mu ^{r-1}h^{1-r})$ and the final asymptotics has its principal part of the magnitude 
$\mu ^r h^{-r} \times (\mu h)^{-r}h^{-2r}$ and a remainder 
$O(\mu ^{r-1}h^{1-r} \times (\mu h)^{-r}h^{-2r})= O(\mu^{-1}h^{1-2r})$.

As $\mu h\ge \epsilon_0$ one would have the family of no more than $C_0$ of $\mu^{-1}h$-pseudo-differential operators, under proper non-degeneracy condition
$|b|+\nabla b|\ge \epsilon_1$, $b\Def \sum_j (2\alpha_j+1)\mu h f_j + V$,
the final asymptotics has its principal part of magnitude $\mu ^r h^{-r}$ and  a remainder $O(\mu^{-1}h^{1-r})$.

\begin{remark}\label{rem-19-1-3}
As $f_j=\const$ non-degeneracy assumption becomes $|\nabla V|\ge \epsilon_0$ and we will be able to replace it by a weaker assumption $|\nabla V|\le \epsilon_0\implies |\det \Hess V|\ge \epsilon_0$.
\end{remark}

\section{Resonances}
\label{sect-19-1-3}

However, the resonances prevent us from reducing our operator to the desired canonical form (\ref{19-1-6}) even in the smooth case. In fact, \emph{$m$-th order resonances\/}\index{resonances!$m$-th order}
\begin{equation}
\sum _j f_j \gamma _j =0\qquad \text{as\ \ } \gamma \in \bZ^r\; \text{and\ \ }
|\gamma |\Def \sum_j|\gamma_j|=m\ge 2
\label{19-1-17}
\end{equation}
prevent us from reducing properly terms
\begin{equation}
\sum _{\alpha ,\beta:  |\alpha|+|\beta| =m;\  k,j }
\mu^{2-2k-|\alpha|-|\beta|-j} h^j a_{\alpha, \beta; k,j}
(x'',\mu^{-1}hD'') (hD')^\alpha (\mu x')^\beta.
\label{19-1-18}
\end{equation}
In particular, due to the \emph{$2$-nd order resonances\/}\index{resonances!$2$-nd order} $f_j\approx f_k$ ($k\ne j$) we can reduce the main part of operator only to
\begin{equation}
\sum_{\fm \in \fM } \sum _{j,k\in \fm } a_{jk}(x'',\mu ^{-1}hD'') Z_j^* Z_k +  a_0(x'',\mu ^{-1}hD'')
\label{19-1-19}
\end{equation}
with $Z_k=hD_k+i\mu x_k$ where $\fm \in \fM$ are disjoint subsets of
$\{1,\dots, r\}$ and eigenvalues of each of matrices
$(a_{jk})_{j,k\in \fm }$ are  close to one another (and to $f_j$). This leads to the necessity of the matrix rather than the scalar non-degeneracy (microhyperbolicity) condition.

Further, the \emph{$3$-rd order resonances\/}\index{resonances!$3$-rd order} $f_i\approx f_j+f_k$ (with possible $j=k$) prevent  us from getting rid off  the terms
\begin{gather}
\mu ^{-1} a_{jkl}(x'',\mu ^{-1}hD'') Z_j^* Z_k Z_l
\label{19-1-20}\\
\intertext{and their adjoints}
\mu ^{-1} a^*_{jkl}(x'',\mu ^{-1}hD'') Z_j Z^*_k Z^*_l
\tag*{$\textup{(\ref*{19-1-20})}^*$}\label{19-1-20-*}
\end{gather}
and these terms appear as the perturbations of the main part of the operator
unless $F_{jk}$ and $g^{jk}$ are constant.

Furthermore, the fourth order resonances $f_i=f_j+f_k+f_l$ and $f_i+f_j=f_k+f_l$ where $i,j,k,l$ may coincide but in the second case $(i,j)\ne (k,l)$ leave us with the whole bunch of terms instead of just $\mu ^{-2}a_{jk}H_j H_k$ and but these terms are smaller than those produced by the third order resonances and rather harmless. 

More precisely, the symbol of the reduced operator is
\begin{equation}
\mu ^2 \sum_{\fm \in \fM}\sum _{j,k \in \fm } a_{jk} \zeta_j^\dag \zeta_k +
a_0 +
\mu ^2 \sum_{\fn \in \fN } \sum _{j,k,l\in \fn} \Re a_{jkl} \zeta_j^\dag \zeta_k \zeta_l +\ldots
\label{19-1-21}
\end{equation}
where $\zeta_j = \xi_j+ ix_j$ are symbols of $\mu^{-1}h$ differential operators $\mu^{-1}Z_j$ and precise definition of $\fN$ will appear later. Recall that $^\dag$ means a complex conjugation for scalars and a Hermitian conjugation for matrices.

\begin{remark}\label{rem-19-1-4}
(i) Resonances of $m$-th order become important only as smoothness is large enough (at least $(m,0)$).

\medskip\noindent
(ii) If we assume that $g^{jk}$ and $F_{jk}$ are constant then $m$-th order resonances affect  only  terms (\ref{19-1-18}) without factor $\mu^2$ i.e.
\begin{equation}
\sum _{\alpha ,\beta:  |\alpha|+|\beta| =m;\  k,j }
\mu^{-2k-|\alpha|-|\beta|-j} h^j a_{\alpha, \beta; k,j}
(x'',\mu^{-1}hD'') (hD')^\alpha (\mu x')^\beta.
\tag*{$\textup{(\ref*{19-1-18})}'$}\label{19-1-18-'}
\end{equation}
\end{remark}

\section{Dynamics and microhyperbolicity}
\label{sect-19-1-4}

If operator in its canonical form  is 
\begin{equation}
\sum_{1\le j\le r} f_j(x'',\mu^{-1}hD'') (h^2D_j+\mu^2 x_j^2)+
V(x'',\mu^{-1}hD'')
\label{19-1-22}
\end{equation}
then in the classical dynamics $(\xi_j^2+\mu^2x_j)^2=\rho_j=\const$ and dynamics in $(x',\xi')$ is described by the Hamiltonian
\begin{equation}
b(x'',\xi'',\rho_1,\ldots,\rho_r)\Def \sum_{1\le j\le r} f_j(x'',\xi'') \rho_j +V (x'',\xi'');
\label{19-1-23}
\end{equation}
where with respect to $x''=(x_{r+1},\ldots , x_d)$ we use $\hslash$-quantization with $\hslash=\mu^{-1}h$ and therefore this dynamics is
\begin{multline}
\frac{d x''}{dt}= \mu^{-1}\nabla_{\xi''} b(x'',\xi''),\qquad 
\frac{d \xi''}{dt}= -\mu^{-1}\nabla_{x''} b(x'',\xi'')\\
b(x'',\xi'')=\tau
\label{19-1-24}
\end{multline}
and therefore microhyperbolicity assumption $|b-\tau|+|\nabla b|\ge \epsilon_0$ depends on $\boldrho$; for each $\boldrho$ with 
$|b(x'',\xi'',\boldrho)-\tau|\le \epsilon_0$ there exists direction $\ell=\ell(\boldrho)$ such that 
$\langle \ell,\nabla\rangle b(x'',\xi'',\boldrho)\ge \epsilon_0$ and we know that this is preserved until time $T^*=\epsilon \mu$. 

However if there are $3$-rd order resonances situation becomes more complicated: $\rho_j$ are evolving with the speed $O(1)$ and we can take only $T^*=\epsilon$ (which leads to the less sharp remainder estimate).

To remedy this situation we may assume that $\ell$ does not depend on $\boldrho$ which is the case when $f_j=\const$ (then microhyperbolicity condition is $|\nabla V|\ge \epsilon_0$). More generally, we need to assume that  $\{1,\ldots,r\}$ is partitioned into subsets $\fn\in \fN$ such that $2$-nd and $3$-rd order resonances involve indices from the same subset $\fn$ only, and that $\ell$ depends on $\boldlambda=(\lambda_\fn)_{\fn\in \fN}$ with  
$\lambda_\fn = \sum_{j\in \fn} f_j \rho_j$ rather than on $\boldrho$. 

More general definitions of microhyperbolicity and $\fN$-microhyperbolicity working also in the cases when even the main part is not in the form (\ref{19-1-22}) will be given later (see definitions~\ref{def-19-2-4} and~\ref{def-19-2-5}).

\section{Regularity assumptions and mollification}
\label{sect-19-1-5}

We assume that 
\begin{phantomequation}\label{19-1-25}\end{phantomequation}
\begin{gather}
g^{jk}\in \sC^{\bar{l},\bar{\sigma}}, \qquad V\in \sC^{l,\sigma}, \tag*{$\textup{(\ref*{19-1-25})}_{1,3}$}\label{19-1-25-1}\\
\shortintertext{and} 
V_j = \partial_j \phi _j, \qquad \phi_j \in \sC^{\bar{l}+2,\bar{\sigma}}
\tag*{$\textup{(\ref*{19-1-25})}_{2}$}\label{19-1-25-2}
\end{gather}
where the last assumption is a bit stronger than more natural 
$F_{jk}\in  \sC^{\bar{l},\bar{\sigma}}$ or 
$V_\in \sC^{\bar{l}+1,\bar{\sigma}}$.

Also due to problems appearing when we reduce to a canonical form, mollification parameter in the intermediate magnetic field case is $\varepsilon=C\mu^{-1}$ rather than $\varepsilon =C(\mu^{-1}h|\log h|)^{\frac{1}{2}}$ as in the previous Chapter~\ref{book_new-sect-18} of \cite{futurebook}, therefore the threshold between weak magnetic field ($\varepsilon=C\mu h|\log h|$) and intermediate one is $\mu = (h|\log h)^{-\frac{1}{2}}$ rather than $\mu = (h|\log h)^{-\frac{1}{3}}$  and we also need to assume a larger regularity than there. 

\section{Plan of the Chapter}
\label{sect-19-1-6}

This Chapter consists of six more Sections. 

As usual, we start from the Section~\ref{sect-19-2} devoted to the weak magnetic field when asymptotics is defined by evolution to time $T=\epsilon \mu^{-1}$. We also study the classical dynamics here.

Again, following the standard scheme, Section~\ref{sect-19-3} is devoted to the canonical form, in Section~\ref{sect-19-4} we consider intermediate magnetic field $\mu \le \epsilon (h|\log h|)^{-1}$ (overlapping with Section~\ref{sect-19-2}) and in Section~\ref{sect-19-5} we consider a strong magnetic field $\epsilon (h|\log h|)^{-1}\mu \le  \epsilon h^{-1}$. 

In the general case we assume some microhyperbolicity assumption (see definition~\ref{def-19-2-4} but there is also stronger  $\fN$-microhyperbolicity, see definition~\ref{def-19-2-5}). If there are no $2$-nd order resonances, then this condition means that $\lambda_1\nabla(V/f_1)+\ldots+\lambda_r\nabla(V/f_r)\ne 0$ as $\lambda_1\ge 0,\ldots ,\lambda_r\ge 0, \lambda_1+\ldots +\lambda_r>0$.

We also consider case of constant $g^{jk}$, $F_{jk}$ and in this case instead of microhyperbolicity assumption $\nabla V\ne 0$ we assume a weaker non-degeneracy assumption $\nabla V=0\implies \det\Hess V \ne 0$.

Finally, for the main course, in Section~\ref{sect-19-6} we consider a very strong and superstrong magnetic field $\epsilon h^{-1}\le \mu \le Ch^{-1}$ and $\mu \ge Ch^{-1}$ respectively.

Furthermore, in Section~\ref{sect-19-7} we consider the case when the rank of  $\{\nabla(V/f_1),\ldots, \nabla(V/f_r)\}$ is less than $r$ (usually $(r-1)$) at some points. There we assume that coefficients are regular.

\chapter{Weak magnetic field}
\label{sect-19-2}

As we mentioned, in the weak magnetic field we consider original $(x,\xi)$ coordinates. As we explained,  we can take $x$-scale exactly as in Chapter~\ref{book_new-sect-18} of \cite{futurebook}
\begin{equation}
\varepsilon = C\mu h |\log h|
\label{19-2-1}
\end{equation}
and therefore to satisfy logarithmic uncertainty principle  we can take $\xi$-scale $\mu ^{-1}$\,\footnote{\label{foot-19-2} Recall that we consider $\mu^{-1}h$-differential operators.}. By no means $\varepsilon$ is larger than
$C(\mu ^{-1}h|\log h|)^{\frac{1}{2}}$ needed to consider 
$\mu ^{-1}h$-pseudo-differential operator with the symbol smooth in $\varepsilon$-scale in both $x$ and $\xi$ and therefore we cannot make a reduction here; only some kind of quasi-reduction as in Chapter~\ref{book_new-sect-18} of \cite{futurebook} in case the of weak magnetic field it is possible.

\section{Heuristics. Classical dynamics}
\label{sect-19-2-1}

\subsection{Smooth theory. Classical dynamics}
\label{sect-19-2-1-1}

Here we are interested in the classical dynamics generated by symbol
$a(x,\xi)\Def \mu^2 a^0(x,\xi) +V(x)$ with
\begin{equation}
a^0(x,\xi) \Def \sum _{j,k}g^{jk}p_jp_k,\qquad p_j= \mu ^{-1}\xi_j - V_j(x)
\label{19-2-2}
\end{equation}
at the energy levels close to $0$ for time $t:|t|\le T^*$ with 
$T^*=\epsilon \mu $. On the energy levels below $c_0$ we have
\begin{equation}
|p_j|\le C\mu ^{-1}.
\label{19-2-3}
\end{equation}

At this heuristic stage we will not look at the smoothness (thus assuming that $l$ is large enough).

First of all, exactly as in the previous Chapter~\ref{book_new-sect-18} of \cite{futurebook}, we can find
$\upphi_{jk} \in \sF^{\bar{l},\bar{\sigma}}$ (see remark~\ref{book_new-rem-18-3-2} of \cite{futurebook}) such that for $q_j\Def x_j - \sum _k \upphi_{jk} (x) p_k$
\begin{equation}
\{p_m,q_j\}= O(\mu ^{-2}).
\label{19-2-4}
\end{equation}

Really,
\begin{equation}
\{p_m,p_k\}= \mu ^{-1}F_{km}
\label{19-2-5}
\end{equation}
and therefore
$\{p_m,q_j\}= \mu ^{-1}\updelta_{mj}-\mu ^{-1}\sum_{j,k} \upphi_{jk}F_{km}$
and one needs to take
\begin{equation}
(\upphi_{jk})=(F_{jk})^{-1}.
\label{19-2-6}
\end{equation}
Then one can see easily that
\begin{gather}
\{q_j,q_k\}\equiv \{x_j,q_k\}\equiv \mu ^{-1} \upphi_{jk} \qquad\mod O(\mu ^{-2})
\label{19-2-7}\\
\intertext{and therefore}
\mu ^2 \{a^0, q_j\} = O(\mu ^{-1}), \qquad \{V,q_j\}=O(\mu^{-1})
\label{19-2-8}
\end{gather}
which implies that
\begin{claim}\label{19-2-9}
In the classical evolution at energy levels below $c_0$ for time
$t:|t|\le \epsilon \mu$ increment of $q_j$ does not exceed
$C\mu ^{-1}|t|\le \epsilon _1$;
\end{claim}
\vglue-10pt
\begin{claim}\label{19-2-10}
In the classical evolution at energy levels below $c_0$ for time
$t:|t|\le \epsilon \mu$ increment of $x_j$ does not exceed
$C\mu ^{-1}(|t|+1)\le \epsilon_1$.
\end{claim}

One can see easily that in the case of constant $g^{jk},F_{jk}$
\begin{equation*}
\mu ^2 \{a^0, q_j\} = 0, \qquad \{V,q_j\}=\mu^{-1}\sum_k \upphi_{jk}\partial _kV
\tag*{$\textup{(\ref*{19-2-8})}'$}\label{19-2-8-'}
\end{equation*}
and therefore
\begin{claim}\label{19-2-11}
In the case of constant $g^{jk},F_{jk}$ the classical evolution for time  $t:|t|\le \epsilon \mu$ is approximately described by equation\footnote{\label{foot-19-3} This equation describes approximately coordinates $x_j$ but \emph{not} velocities.}
\begin{equation}
\frac {d\ }{dt}x_j=\mu^{-1}\sum_k \upphi_{jk}\partial _kV
\label{19-2-12}
\end{equation}
along which $V$ is preserved.
\end{claim}

To describe evolution more precisely in the general case let us consider a point $\bar{x}$. Consider $\fM (\bar{x})$  a partition  of $\{1,\ldots,r\}$ such that  $|(f_i-f_j)(\bar{x})|\le \epsilon $ iff $i,j$ belong to the same element $\fm$ of this partition (so $2$-nd  order resonances involve indices from the same subset $\fm$ only).

Then we can introduce
\begin{equation}
\zeta_j(x,\xi)=\sum_{1\le k\le 2r} \alpha_{jk}(x)p_k(x,\xi), \qquad j=1,\ldots,r
\label{19-2-13}
\end{equation}
with complex coefficients $\alpha_{jk}\in \sF^{\bar{l},\bar{\sigma}}$ such that
\begin{multline}
\{\zeta_j,\zeta_k\}=0,\quad
\{\zeta_j^\dag,\zeta _k\} =2i\mu ^{-1}\updelta_{jk}\\
\text{at\ \ } \Sigma^0\Def\bigl\{(x,\xi):\ p_1=\ldots =p_{2r}=0\bigr\}
\label{19-2-14}
\end{multline}
where as usual $^\dag$ means complex conjugation, and therefore
\begin{multline}
\{\zeta_j,\zeta_k\}\equiv 0,\quad
\{\zeta_j^\dag,\zeta_k\}\equiv 2i\mu ^{-1}\updelta_{jk}  \quad
\mod O(\mu^{-2})\\
\text{as\ \ } |p|\le c\mu^{-1}
\label{19-2-15}
\end{multline}
and also
\begin{equation}
a^0 =\frac{1}{2}\mu^2\sum _{\fm}
\sum _{j,k\in \fm} a_{jk}(x)\zeta_j^\dag(x,\xi)\zeta_k(x,\xi)
\label{19-2-16}
\end{equation}
with Hermitian matrices $(a_{jk})_{j,k\in \fm}$ close to scalar matrices $f_{\fm}I_{\#\fm}$ of the same dimensions.

Note that $\Sigma^0$ is a symplectic manifold and $\bR^d\ni x$ inherits its symplectic form, a volume form, and a Hamiltonian field (we call it ``Liouvillian field'') respectively:
\begin{gather}
\frac{1}{2}\sum_{j,k} F_{jk}\,dx_j\wedge \,dx_k, \label{19-2-17}\\
|\det F|^{\frac{1}{2}}\,dx_1\wedge\ldots\wedge dx_d=
f_1\cdots f_r \sqrt{g}\,dx,
\label{19-2-18}\\[2pt]
\cL_\psi = -\sum _{j,k} \upphi_{jk}(\partial_j\psi)\partial_k.\label{19-2-19}
\end{gather}
Note also that
\begin{gather}
a^0_{\fm}\Def \frac{1}{2}
\mu^2\sum _{j,k\in \fm } a_{jk}(x)\zeta_j^\dag(x,\xi)\zeta_k(x,\xi)
\label{19-2-20}\\
\shortintertext{satisfy}
\{a^0, a^0_\fm\} =O(1),\qquad\{V,a^0_\fm\}=O(1).
\label{19-2-21}
\end{gather}
Therefore an increment of $a^0_{\fm}$ for time $t$ does not exceed $C|t|$
and thus $a^0_{\fm}$ are \emph{moderate\/} but not necessarily \emph{slow variables\/} and their increment for $t=T^*\asymp \mu$ could be rather large. 

Now we can replace $\zeta_j$ by
\begin{equation}
\zeta_{j,\new}=\zeta_j- \sum _{kl}\beta_{jkm}(x) \eta_k\eta_m
\label{19-2-22}
\end{equation}
with $\beta_{jkm}\in \sF^{\bar{l}-1,\bar{\sigma}}$ to make
\begin{equation}
\{\zeta_j,\zeta_k\}\equiv 0,\quad
\{\zeta_j^\dag,\zeta_k\} \equiv 2i\mu ^{-1}\updelta_{jk}\qquad
\mod O(\mu ^{-3})
\label{19-2-23}
\end{equation}
where for the sake of simplicity we use notations
$\eta _{2k-1}\Def \Re \zeta_k$, $\eta _{2k}\Def \Im \zeta _k$, $k=1,\ldots,r$.

Then
\begin{align}
& a^0 =
\frac{1}{2}\mu^2\sum _{\fm}\sum _{j,k\in \fm } a_{jk}(q)\zeta_j^\dag \zeta_k +
\frac{1}{3} \mu^2 \sum_{j,k,m} a_{jkl}\, \eta _j\eta_k\eta_m +\ldots, \label{19-2-24}\\[2pt]
&V= V(q)+ \sum _j a_j(q)p_j+\ldots
\label{19-2-25}
\end{align}
where dots denote terms with $O(\mu^{-1})$ gradients.

These cubic terms in (\ref{19-2-24}) and linear terms in (\ref{19-2-25}) are the only sources of the trouble because Poisson brackets
$ \{\tilde{a}^0_{\fm}, \tilde{a}^0_{\fm'}\}=O(\mu ^{-1})$
and $ \{ \tilde{a}^0_{\fm}, V(q)\}=O(\mu ^{-1})$
where $\tilde{a}^0_{\fm}$ are defined by (\ref{19-2-20}) with  $a_{jk}(x)$ replaced  by $a_{jk}(q)$.

To get rid of these terms we can redefine $\zeta_j$ again replacing them
by $\zeta_j + \{S,\zeta_j\}$ with
\begin{equation}
S= \mu \sum _{k,m,n } \gamma_{kmn} \eta_k\eta_m\eta_n +
\mu ^{-1} \sum _m \gamma_m\eta_m.
\label{19-2-26}
\end{equation}
This replacement preserves (\ref{19-2-14}), (\ref{19-2-15}) and (\ref{19-2-19}) and, modulo terms with $O(\mu^{-1})$ gradients, it is equivalent to the replacement of $ a^0 +V $ by
\begin{equation}
 a^0+V +\{S,  a^0 +V\}.
\label{19-2-27}
\end{equation}

We can always choose $\gamma_j$ killing $a_j$ and if there is no $3$-rd order resonances we can choose $\gamma_{jkm}$ killing $a_{jkm}$; in the general case we can reduce $a^0+V$ to
\begin{equation}
\sum _{\fm}   \tilde{a}^0_{\fm} +
\mu ^2 \Re  \sum _{j,k,m: |f_j-f_k-f_m|\le \epsilon} a_{jkm}(q)\,\zeta _j^\dag \zeta_k \zeta_m +V (q).
\label{19-2-28}
\end{equation}

\begin{example}\label{example-19-2-1}
Consider $p_j=\xi_j$ as $j=1,\ldots,r$ and $p_{j+r}= \xi_{r+j}+x_j$; then $\zeta_j = \xi_j + i (\xi_{r+j}+x_j)$. Let consider a symbol which is the real part of the quadratic form of $(\zeta_1,\zeta^\dag_1,\ldots \zeta_r,\zeta^\dag_r)$ with a linear coefficients coinciding at $x=0$ with $\sum_j f_j\zeta^\dag_j\zeta_j$ with constant $f_1,\ldots,f_r$.

Let  the only term with non-constant coefficients be
\begin{multline*}
2\Re \bigl(\alpha x_1 \zeta_2\zeta_3\bigr) =
\Re \bigl(\alpha i (\zeta_1^\dag -\zeta_1 +2i \xi_{r+1})\zeta_2\zeta_3\bigr) =\\
\Re \bigl(\alpha i \zeta_1^\dag \zeta_2\zeta_3\bigr)-
\Re \bigl(\alpha i\zeta_1\zeta_2\zeta_3\bigr)- 2  \Re \bigl(\alpha\xi_{r+j}\zeta_2\zeta_3\bigr)
\end{multline*}
with $\alpha \in \bC$. Here $\xi_{1+r}$ is nothing but a constant parameter in the dynamics; let pick it equal to $0$. 

Let 
\begin{enumerate}[label=(\roman*), leftmargin=*]
\item $r=2$, $\zeta_3\Def \zeta_2$; one can see easily that we arrive exactly to the case with non-zero cubic terms (the last term in the right-hand expression could be removed but not the first one as $f_1=2f_2$); 

\item $r=3$; one can see easily that we arrive exactly to the case with non-zero cubic terms (the last term in the right-hand expression could be removed but not the first one as $f_1=f_2+f_3$).
\end{enumerate}
\end{example}

\begin{example}\label{example-19-2-2}
(i) Consider
\begin{equation*}
\mu^{-2}a = f_1 \zeta^\dag_1 \zeta^\dag_1 + f_2 \zeta^\dag_2 \zeta^\dag_2 +
\Re \bigl(\alpha \zeta_1^\dag \zeta_2^2\bigr)
\end{equation*}
with $f_1=2f_2$; let $\varrho_j= \mu^2 \zeta_j^\dag\zeta_j$; then
\begin{gather*}
\frac{d\varrho_2}{dt}= \{a, \varrho_2\}= \mu^4 \Re \{\alpha \zeta_1^\dag \zeta_2^2, \zeta_2^\dag \zeta_2\}= -4\mu^3 \Re \bigl(\alpha i \zeta_1^\dag \zeta_2^2\bigr);\\[2pt]
\frac{d\ }{dt}\mu^3 \Re (\alpha \zeta_1^\dag \zeta_2^2)=
\{a,\mu^3 \Re (\alpha \zeta_1^\dag \zeta_2^2)\}=0\implies
\mu^3\Re (\alpha \zeta_1^\dag \zeta_2^2) = |\alpha|\sigma (=\const)
\end{gather*}
where $\sigma$ can accept any value in the interval $[-\varrho_1^{\frac{1}{2}}\varrho_2, \varrho_1^{\frac{1}{2}}\varrho_2]$ and then
\begin{equation*}
\frac{d^2\varrho_2}{dt^2}=\{a,\{a, \varrho_2\}\}=
-4\mu^5  \{ \Re (\alpha \zeta_1^\dag \zeta_2^2),
\Re (\alpha i\zeta_1^\dag \zeta_2^2)\}=
 -2|\alpha|^2 \varrho_2 (\varrho_2 -4 \varrho_1);
\end{equation*}
therefore as
$2\varrho_1 + \varrho_2= \lambda=(\const)$ we arrive to
\begin{equation}
\frac{d^2\varrho_2}{dt^2} = -2|\alpha|^2  (3\varrho_2^2 -2\lambda).
\label{19-2-29}
\end{equation}
Equation (\ref{19-2-29}) describes a kind of oscillations in the real time in the interval
\begin{equation*}
J_\lambda= 
\{\varrho_2 \ge 0, \ \varrho_2^3-\lambda \varrho_2^2 +2 \sigma^2 \le 0\}.
\end{equation*}
Thus $a^0_{j}=\varrho_j$ are really moderate but not slow variables.

\medskip\noindent
(ii) Consider
\begin{equation*}
\mu^{-2}a = f_1 \zeta^\dag_1 \zeta^\dag_1 + f_2 \zeta^\dag_2 \zeta^\dag_2 +
f_3 \zeta^\dag_3 \zeta^\dag_3+
\Re \bigl(\alpha \zeta_1^\dag \zeta_2\zeta_3\bigr)
\end{equation*}
with $f_1=f_2+f_3$; let $\varrho_j= \mu^2 \zeta_j^\dag\zeta_j$; then for $j=2,3$
\begin{equation*}
\frac{d\varrho_j}{dt}=\{a, \varrho_j\}=
\mu^4 \Re \{\alpha \zeta_1^\dag \zeta_2\zeta_3, \zeta_j^\dag \zeta_j\}=
-2\mu^3 \Re \bigl(\alpha i \zeta_1^\dag \zeta_2\zeta_3\bigr)
\end{equation*}
and again $\mu^3\Re (\alpha \zeta_1^\dag \zeta_2\zeta_3) = |\alpha|\sigma (=\const)$ where $\sigma$ can accept any value in the interval
$[-(\varrho_1\varrho_2\varrho_3)^{\frac{1}{2}}, (\varrho_1\varrho_2\varrho_3)^{\frac{1}{2}}]$
\begin{multline*}
\frac{d^2\varrho_2}{dt^2}=\{a,\{a, \varrho_j\}\}=
-2\mu^5  \{ \Re (\alpha \zeta_1^\dag \zeta_2\zeta_3),
\Re (\alpha i\zeta_1^\dag \zeta_2\zeta_3)\}=\\
 -|\alpha|^2 (\varrho_2\varrho_3 - \varrho_1\varrho_2-\varrho_1\varrho_3).
\end{multline*}
Due to the first equation with $j=2,3$ \
$\varrho_2-\varrho_3= \const$ and $\varrho_{2,3}=\varrho \pm \eta$
with constant $\eta$ and $\varrho=\frac{1}{2}(\varrho_2+\varrho_3)$;
therefore as $f_1\varrho_1 + f_2\varrho_2+f_3 \varrho_3=\const$ we conclude that $\varrho_1= -\varrho + \lambda$ and arrive to equation
\begin{equation}
\frac{d^2\varrho}{dt^2} = -|\alpha|^2  (3\varrho^2 -2\lambda \varrho +\eta^2).
\label{19-2-30}
\end{equation}

Again, equation (\ref{19-2-30}) describes a kind of oscillations in the real time in the interval
\begin{equation*}
J_{\lambda,\eta}=\{\varrho \ge 0, \
\varrho ^3-\lambda \varrho ^2 + \eta^2 \varrho  -\lambda\eta^2+\sigma^2\le 0\}.
\end{equation*}
Thus $a^0_{j}=\varrho_j$ again are really moderate but not slow variables.
\end{example}

To have $a^0_j$ as slow variables we need \underline{either} to assume that there are no cubic terms (i.e. either $F_{jk},g^{jk}$ are constant or there are no third order resonances) \underline{or} to group them together.

The important role will be played by partition $\fN=\fN_\epsilon (\bar{x})$ of $\{1,\ldots,r\}$ such that

\begin{claim}\label{19-2-31}
$j\ne k$ belong to the same element $\fn\in \fN$ if there exists $l$ such that either $|f_j- f_k-f_l|\le \epsilon$ or $|f_k- f_j-f_l|\le \epsilon$ or
$|f_l- f_j-f_k|\le \epsilon$.
\end{claim}

\begin{remark}\label{rem-19-2-3}
(i) Then, according to this definition $l$ also belongs to the same element $\fn$. We do not exclude $l= j$ or $l=k$;

\medskip\noindent
(ii) Obviously $\fM$ is a subpartition of $\fN$; recall that $\fM$ groups together only indices $j,k$ such that $f_j\approx f_k$;

\medskip\noindent
(iii) It may happen that indices $j,k$ must belong to the same element $\fn\in \fN$ even if they are not part of the same resonance equation (f.e. $1,2$ are swapped together by $f_1\approx 2f_2$ and $2,3$ are swapped together by $f_2\approx f_3$). This was not the case with $\fM$; 

\medskip\noindent
(iv) We do not assume that $\fN$ is the finest partition satisfying (\ref{19-2-31}) as we do not want to exclude possibility to take $\#\fN=1$ (i.e. no partition at all).
\end{remark}

Then, after we got rid of all non-resonant cubic term let us define
\begin{gather}
\tilde{a}_{\fn}^0\Def \sum _{\fm\subset\fn}
\tilde{a}_{\fm} +
\mu^2 \Re  \sum _{j,k,m\in \fn } a_{jkm}(q)\zeta _j^\dag \zeta_k \zeta_m.
\label{19-2-32}
\\
\shortintertext{Then}
\{\tilde{a}_{\fn}^0, \tilde{a}^0_{{\fn}'}\}=O(\mu ^{-1}),\qquad
\{\tilde{a}_{\fn}^0, V(q)\}=O(\mu ^{-1})
\label{19-2-33}
\end{gather}
and therefore
\begin{claim}\label{19-2-34}
In the classical evolution at energy levels below $c_0$ for time
$t:|t|\le \epsilon \mu$ increment of $\tilde{a}_{\fn}^0$ does not exceed
$C\mu ^{-1}|t|\le \epsilon _1$ (so $\tilde{a}_{\fn}^0$ are slow variables);
\end{claim}
and
\begin{claim}\label{19-2-35}
In the classical evolution at energy levels below $c_0$ for time
$t:|t|\le \epsilon \mu$ increment of
$a^0_{\fn}\Def \sum_{{\fm}\subset {\fn}} a^0_{\fm} $ does not exceed
$C\mu ^{-1}(|t|+1)\le \epsilon_1$.
\end{claim}
It follows  from (\ref{19-2-7}) that
\begin{equation}
\{q_j,q_k\}\equiv  \mu ^{-1} f_{jk} (q)+ O(\mu ^{-2})
\tag*{$\textup{(\ref*{19-2-7})}'$}\label{19-2-7-'}
\end{equation}
and then one can reintroduce
$q_j\Def x_j - \sum _k f_{jk} (x) \eta_k-\sum _{k,m} \rho_{jkm} (x) \eta_k\eta_m$ to make
\begin{gather}
\{\eta_k,q_j\}= O( \mu ^{-3}).
\label{19-2-36}\\
\shortintertext{Then}
\{ a^0 +V , q_j\} \equiv \mu ^2 \sum _{k,m} \{a_{km} ,q_j\}
\zeta_k\zeta_m ^\dag
+ \{V(q), q_j\} \mod O(\mu ^{-3}).
\label{19-2-37}
\end{gather}

\subsection{Microhyperbolicity}
\label{sect-19-2-1-2}

Let us introduce microhyperbolicity condition provided $a_0$ does not vanish i.e.
\begin{equation}
V\le - \epsilon_0;
\label{19-2-38}
\end{equation}
we are not interested in the classically forbidden case of $V\ge \epsilon_0$.

\begin{definition}\label{def-19-2-4}
We call magnetic Schr\"odinger operator \emph{microhyperbolic\/}\index{Schr\"odinger operator!microhyperbolic} at point $\bar{x}$ if for each $\boldtau=(\tau_{\fm})_{{\fm}\in \fM}$ such that
$|\sum _{{\fm}\in \fM} \tau _{\fm} +V|\le \epsilon$  there exists vector
$\ell= \ell (\bar{z}, {\boldtau})\in \bR^{2r}$ such that
\begin{gather}
-\mu^2 \sum_{j,k} \bigl(\ell (a_{jk}a_0^{-1}) \bigr)\zeta _j^\dag \zeta _k\ge \epsilon_0
\label{19-2-39}\\
\shortintertext{as long as}
\mu^2\sum _{j,k\in {\fm} }a_{jk}\zeta_j^\dag  \zeta_k =  \tau_{\fm}
\qquad \forall \zeta\in \bC^r\quad \forall {\fm}\in \fM.
\label{19-2-40}
\end{gather}
\end{definition}

\begin{definition}\label{def-19-2-5}
Let $\fN=\fN(\bar{x})$ satisfy (\ref{19-2-31}). We call magnetic Schr\"odinger operator \emph{$\fN$-microhyperbolic\/}\index{Schr\"odinger operator!microhyperbolic!$\fN$-microhyperbolic}  at point $\bar{x}$ if for each $\boldtau=(\tau_{\fn})_{{\fn}\in \fN}$ such that
$|\sum _{{\fn}\in \fN} \tau _{\fn} +V|\le \epsilon$  there exists vector
$\ell= \ell (\bar{z}, {\boldtau})\in \bR^{2r}$ such that (\ref{19-2-39}) holds as long as
\begin{equation}
\mu^2\sum _{j,k\in {\fn} }a_{jk}\zeta_j^\dag  \zeta_k =  \tau_{\fn}
\qquad \forall \zeta\in \bC^r\quad \forall {\fn}\in \fN.
\label{19-2-41}
\end{equation}
\end{definition}

Then we conclude that
\begin{claim}\label{19-2-42}
If operator is microhyperbolic at point $\bar{x}$ then  $\sum_j \ell_jq_j$ will  increase with the exact rate $\asymp \mu ^{-1}$ for time $|t|\le T^*=\epsilon $
\end{claim}
and
\begin{claim}\label{19-2-43}
If  operator is \underline{either} $\fN$-microhyperbolic at point $\bar{x}$ \underline{or} contains no cubic terms\footnote{\label{foot-19-4} Which happens for example if either there are no $3$-rd order resonances or if $F_{jk},g^{jk}$ are constant.} then in the frames of the smooth theory one can take  $T^*=\epsilon \mu$.
\end{claim}

Here both microhyperbolicity and $\fN$-microhyperbolicity obviously mean ``on the energy level $0$''. 

\begin{example}\label{example-19-2-6}
(i) If $F_{jk},g^{jk}$ are constant then the microhyperbolicity condition means exactly that $|\nabla V|\ge \epsilon$. In this case also $\fN$-microhyperbolicity condition is fulfilled with $\#\fN=1$;

\medskip\noindent
(ii) If there are no $2$-nd order resonances\footnote{\label{foot-19-5} Or if eigenvalues $f_j$ of $(F^j_k)$ have constant multiplicities.} then the microhyperbolicity condition means exactly that
\begin{multline}
|\sum _k \lambda_k \nabla \log (-V/f_k)|\ge \epsilon_0 \\
\forall
\lambda_1\ge 0,\ldots ,\lambda_r\ge 0:\ \lambda_1+\ldots+\lambda_r=1.
\label{19-2-44}
\end{multline}
Then there exists vector $\ell$ such that
\begin{multline*}
\sum _k \lambda_k \langle \ell,\nabla\rangle \log (-V/f_k) \ge \epsilon_0 \qquad
\forall
\lambda_1\ge 0,\ldots ,\lambda_r\ge 0:\ \lambda_1+\ldots+\lambda_r=1
%\label{19-2-44}
\end{multline*}
and then we have $\fN$-microhyperbolicity with $\#\fN=1$. So in this case the notions of microhyperbolicity and $\fN$-microhyperbolicity also coincide.

This assumption (\ref{19-2-44}) is fulfilled provided $\nabla (V/f_1),\ldots, \nabla (V/f_r)$ are linearly independent. 
\end{example}

\begin{Problem}\label{Problem-19-2-7}
Either prove that the the notions of microhyperbolicity and $\fN$-microhyperbolicity coincide in the general case or construct a counter-example.
\end{Problem}

\subsection{Non smooth theory}
\label{sect-19-2-1-3}

Consider what happens for not very large smoothness. First of all, all our analysis above obviously holds as $(l,\sigma)\succeq (2,0)$.

Furthermore, as long as we do not need to analyze $a^0_{\fn}$ (i.e. if we have $\#\fN=1$ in the $\fN$-microhyperbolicity condition) $\sC^1$-smoothness is sufficient: to estimate by $O(\mu ^{-1})$ the rate of change of $q$ we need to have only $\{\zeta_j,q_k\}=O(\mu ^{-2})$ and $\{x_j,q_k\}=O(\mu ^{-1})$.

Microhyperbolicity arguments require $\{\zeta_j,q_k\}=o(\mu ^{-2})$ and calculation of $\{x_j,q_k\}$ modulo $o(\mu ^{-1})$. To do it we just need to mollify $\sC$-coefficients $\rho_{jkm}$ with the mollification parameter $o(1)$.

If $(l,\sigma)\prec  (2,0)$ and $\#\fN>1$ in the $\fN$-microhyperbolicity condition, we need to analyze evolution of $a^0_{\fn}$ more carefully. One can see easily that in this case
\begin{equation}
\{a^0_{\fn},  a \}=  O\bigl(\mu ^{-l+1} |\log \mu |^{-\sigma}\bigr)
\label{19-2-45}
\end{equation}
because one can consider $3$-rd order and linear terms and replace them by their
$\mu ^{-1}$-mollifications with respect to $\eta$ and $Q$; then for mollified
terms this estimate would hold and for mollification error it would hold as well since the gradient of the error would be $O(\mu ^{1-l}|\log \mu |^{-\sigma})$.

Therefore under $\fN$-microhyperbolicity assumption with $\#\fN>1$ one can take
\begin{equation}
T^*= \epsilon\left\{\begin{aligned}
& \mu\qquad && (l,\sigma)\succeq (2,0),\\
&\mu ^{l -1}|\log \mu |^{ \sigma }\qquad
&& (l, \sigma )\prec (2,0)
\end{aligned}\right.
\label{19-2-46}
\end{equation}
and in the end of the day semiclassical error will be $O\bigl(h^{1-d}T^{*\,-1}\bigr)$ which is exactly  the first term in estimate (\ref{19-2-87}) below.

\section{Semiclassical propagation}
\label{sect-19-2-2}

Now moving from the classical dynamics to the rigorous analysis we prove few statements, assuming that
\begin{gather}
\mu_0 \le \mu \le \mu_2^*=\epsilon_1 (h|\log h|)^{-1},
\label{19-2-47}\\[2pt]
\varepsilon\ge C h|\log h|
\label{19-2-48}
\end{gather}
where here and below $\mu_0,C$ are large enough constants.

\begin{claim}\label{19-2-49}
In the statements of this Subsubsection  $\phi _1$ is supported in $B(0,1)$, $\phi _2=1$ in $B(0,2)$ and $\chi $ is supported in $[-1,1]$\,\footnote {\label{foot-19-6} Recall that all such auxiliary functions are appropriate in the sense of Section~\ref{book_new-sect-2-3} of \cite{futurebook}.}.
\end{claim}

\subsection{General theory}
\label{sect-19-2-2-1}

\begin{proposition}\label{prop-19-2-8}\footnote{\label{foot-19-7} Finite speed of propagation with respect to $x$; cf. proposition~\ref{book_new-prop-18-2-4} of \cite{futurebook}.}
Let  $\mu \le \epsilon h^{-1}|\log h|^{-1}$\,\footnote{\label{foot-19-8} Lower bound condition is not needed in this statement.} and let
\begin{equation}
M \ge
\sup _{\{\sum_{j,k} g^{jk}\xi_j\xi_k +V = 0\}} 2\sum _k g^{jk}\xi _k +\epsilon
\label{19-2-50}
\end{equation}
with arbitrarily small constant $\epsilon >0$.

Let
\begin{equation}
Ch |\log h| \le T \le T^{*\prime}\Def\epsilon _0.
\label{19-2-51}
\end{equation}
Then
\begin{multline}
|F_{t\to h^{-1}\tau } \chi _T(t) \bigl(1- \phi _{2, MT} (x- \bar{x})\bigr)
\phi_{1, MT} (y-\bar{x}) U(x,y,t) | \le Ch^s\\
\forall \tau \le \epsilon_1
\label{19-2-52}
\end{multline}
where here and below  $\epsilon_1>0 $ is a small enough constant.
\end{proposition}

\begin{proposition}\label{prop-19-2-9}\footnote {\label{foot-19-9} Magnetic propagation; cf. proposition~\ref{book_new-prop-18-2-6} of \cite{futurebook}.}
Let condition \textup{(\ref{19-1-4})} be fulfilled and let
\begin{equation}
T_* \Def C\varepsilon ^{-1}h|\log h|  +C (\mu h |\log h|)^{\frac{1}{2}}
\le T \le T^* = \epsilon _0\mu.
\label{19-2-53}
\end{equation}
Then
\begin{multline}
|F_{t\to h^{-1}\tau } \chi _T(t)
\bigl(1-\phi _{2, \mu^{-1}MT}(q_1-\bar{x}_1,\ldots, q_d-\bar{x}_d) \bigr)^\w_x \times \\[3pt]
U(x,y,t)\,^t\!\bigl(\phi_{1,\mu^{-1}MT}(q_1-\bar{x}_1,\ldots,q_d-\bar{x}_d) \bigr)^{\w}_y |
\le Ch^s\qquad \forall \tau : |\tau |\le \epsilon_1
\label{19-2-54}
\end{multline}
where here and below  $\mu _0>0$, $C_1$, $M$ are large enough constants. Here and below $b^\w$ means Weyl $\mu^{-1}h$-quantization of symbol $b$ and due to condition \textup{(\ref{19-2-53})} logarithmic uncertainty principle holds and this quantization of the symbols involved is justified;

\medskip\noindent
(ii) In particular, for $T\ge C_1$ this inequality holds for $q_1,\ldots,q_d$ replaced by $x_1,\ldots,x_d$.
\end{proposition}

\begin{proposition}\label{prop-19-2-10}
Let $\nu=C\mu h\varepsilon^{-1}|\log h|+ C(\mu h|\log h|)^{\frac{1}{2}}$.

\medskip\noindent
(i) Let $\#\fM \ge 2$. Then
\begin{multline}
|F_{t\to h^{-1}\tau } \chi _T(t)
\bigl(1-\phi _{2, MT+\nu }(a_\fm-\boldtau_\fm) \bigr)^\w_x \times \\[3pt]
U(x,y,t) \,^t\!\bigl(\phi_{1,MT+\nu }(a_\fm-\boldtau_\fm) \bigr)^{\w}_y |
\le Ch^s \qquad \forall \tau : |\tau |\le \epsilon_1.
\label{19-2-55}
\end{multline}
(ii) Let $\#\fN \ge 2$ and  $v=T^{*\,-1}$ with $T^*$ defined by \textup{(\ref{19-2-46})}. Then
\begin{multline}
|F_{t\to h^{-1}\tau } \chi _T(t)
\bigl(1-\phi _{2, v T+\nu }(a_\fn-\boldtau_\fn) \bigr)^\w_x \times \\[3pt]
U(x,y,t) \,^t\!\bigl(\phi_{1,vT+\nu }(a_\fn-\boldtau_\fm) \bigr)^{\w}_y |
\le Ch^s \qquad \forall \tau : |\tau |\le \epsilon_1.
\label{19-2-56}
\end{multline}
\end{proposition}

\begin{proof}[Proofs of propositions \ref{prop-19-2-8}--\ref{prop-19-2-10}]
All propositions \ref{prop-19-2-8}--\ref{prop-19-2-10} are proven by the same scheme as in Section~\ref{book_new-sect-2-3} of \cite{futurebook}:

To prove that the speed with respect to $x$ does not exceed $M$ one can use function
\begin{equation}
\upchi\Bigl(
\bigl( \frac {|x-\bar{x}|^2} {T^2}  + \epsilon^2\bigr)^{\frac{1}{2}} - C\varsigma \frac{t}{T} \Bigr)
\label{19-2-57}
\end{equation}
with $\upchi$ function of the same type as used in theorem~\ref{book_new-thm-2-3-1} of \cite{futurebook},  $\varsigma =\pm 1$ depending on the time direction and arbitrarily small
constant $\epsilon >0$.

To prove that the speed with respect to $Q$ does not exceed $C\mu ^{-1}$
one can use function
\begin{equation}
\upchi \Bigl(
\bigl( \frac {\mu ^2 |q-\bar{q}|^2}{T^2} +\epsilon^2\bigr)^{\frac{1}{2}} - C\varsigma \frac{t}{T} \Bigr).
\label{19-2-58}
\end{equation}

To prove that the speed with respect to $a_{\fn}$ does not exceed $v=T^{*\,-1}$ with $T^*$ defined by  (\ref{19-2-46}) one can use function
\begin{equation}
\upchi \Bigl(\bigl(
\frac {|a_\fn -\tau_\fn |^2 } {v^2T^2} + \epsilon^2\bigr)^{\frac{1}{2}} -
C\varsigma {\frac{t}{T}}\Bigr)
\label{19-2-59}
\end{equation}
and to prove that the speed with respect to $a_{\fm}$ does not exceed $v=C$ one can use the same function (\ref{19-2-59}) with $\fn$ replaced by $\fm$.
\end{proof}

\begin{corollary}\label{cor-19-2-11}
(i) Let $\#\fM\ge 2$. Then  (in the microlocal sense) $\{a_\fm\}$ stays  in the $\epsilon$-vicinity of value $\boldtau$ for time $T^*=\epsilon$;

\medskip\noindent
(ii) Let $\#\fN\ge 2$. Then (in the microlocal sense) $\{a_\fn\}$ stays in the $\epsilon$-vicinity of value $\boldtau$ for time $T^*$  defined by \textup{(\ref{19-2-46})}.
\end{corollary}

\subsection{Microhyperbolic theory}
\label{sect-19-2-2-2}

In this Subsubsection $\chi$ is supported in $[-1,-\frac{1}{2}]\cup[\frac{1}{2},1]$ and  $\bar{\chi}$ is supported in $[-1,1]$ and $\bar{\chi}=1$ on $[-\frac{1}{2},\frac{1}{2}]$\,\footref{foot-19-6}.

\begin{proposition}\label{prop-19-2-12}\footnote{\label{foot-19-10} Weak magnetic field case.}
Let $ (l,\sigma)\succeq (1,2)$ and
\begin{gather}
\mu_0\le \mu \le \mu^*_1 \Def \epsilon_1 (h|\log h|)^{-\frac{1}{2}},
\label{19-2-60}\\[3pt]
\varepsilon =  C'\mu h|\log h|,\label{19-2-61}\\[3pt]
\bar{T}= \epsilon \mu ^{-1}\label{19-2-62}
\end{gather}
where large enough constant $C'$ is chosen in the very last moment. Then

\medskip\noindent
(i)  Under microhyperbolicity condition (see definition \ref{def-19-2-4}) let us pick up $T^*=\epsilon$;

\medskip\noindent
(ii) Under $\fN$-microhyperbolicity condition (see definition \ref{def-19-2-5}) let us define $T^*$ by \textup{(\ref{19-2-46})} as $\#\fN\ge 2$ and
$T^*=\epsilon \mu$ as $\#\fN=1$\,\footnote{\label{foot-19-11} So under $\fN$-microhyperbolicity condition $T^*=\epsilon \mu$ unless $\#\fN\ge 2$ and $(l,\sigma)\prec (2,0)$.}.

\medskip
Then for $T\in [\bar{T},T^*]$, $\tau \in [-\epsilon_1,\epsilon_1]$
\begin{equation}
|F_{t \to h^{-1}\tau} \chi_T(t) \Gamma \psi U | \le Ch^s.
\label{19-2-63}
\end{equation}
\end{proposition}

\begin{proof}  (i) Let us first try function
\begin{equation}
\upchi \Bigl( \frac{\mu}{T}\langle  \ell, q-\bar{q}\rangle +
\epsilon\varsigma \frac{t}{T}\Bigr)
\label{19-2-64}
\end{equation}
with $\ell= \ell (y;\tau _1,\ldots, \tau_{\#\fN} )$.

One can check easily that it is an admissible symbol as long in scale $\rho=\gamma= \mu T$ with respect to $x$ and $\xi$ and the logarithmic uncertainty principle $\rho \gamma \ge C\mu^{-1}h|\log h|$ is fulfilled\footnote{\label{foot-19-12} To fulfill logarithmic uncertainty principle we also need
$\rho \varepsilon \ge C\mu^{-1}h|\log h|$, but it will be fulfilled automatically as $T\ge \bar{T}$ and $\varepsilon \ge C'\mu h|\log h|$.} as
\begin{gather}
T\ge T_*\Def C(\mu h |\log h|)^{\frac{1}{2}};
\label{19-2-65}\\
\intertext{then $T_* \le \bar{T}=\epsilon_1 \mu^{-1}$ as}
\mu \le \epsilon _3 (h|\log h|)^{-\frac{1}{3}}.
\label{19-2-66}
\end{gather}

Then our standard analysis (like in theorem~\ref{book_new-thm-2-3-1} of \cite{futurebook}, together with  propositions \ref{prop-19-2-8}--\ref{prop-19-2-9} and corollary \ref{cor-19-2-11}\footnote{\label{foot-19-13} If we are in framework of assertion (ii) with $\#\fN\ge 2$.} imply estimate (\ref{19-2-63}).

\medskip\noindent
(ii) Now let us consider the case when
\begin{gather}
\epsilon_3 (h|\log h|)^{-\frac{1}{3}}\le \mu
\le \epsilon_1 (h|\log h|)^{-\frac{1}{2}},
\label{19-2-67}\\
\epsilon\mu ^{-1} \le T\le C (\mu h |\log h|)^{\frac{1}{2}};
\label{19-2-68}\\
\shortintertext{then}
\varepsilon \ge C(\mu ^{-1}h|\log h|)^{\frac 1 2}\ge \mu ^{-1}T,
\label{19-2-69}
\end{gather}
function (\ref{19-2-64}) is not necessarily admissible symbol and one needs a bit more subtle arguments.

Let us consider a small vicinity of $y$. Without any loss of the generality one can assume that $\upphi_{jk}(y)=0$ as $|j-k|\ne r$ (otherwise one can achieve it by a rotation). After this, without any loss of the generality one can assume that $\ell = (\ell',0)$ where $x'=(x_1,\ldots, x_r)$, 
$x''=(x_{r+1},\ldots, x_{2r})$ etc. as we can achieve it by rotations in the planes $(x_j,x_{j+r})$, $j=1,\ldots,r$.

Recall that by definition
\begin{equation*}
q_j=x_j-\sum_k \upphi_{jk}p_k + \frac{1}{2}\sum_{i,k}\beta_{jik}p_ip_k,
\end{equation*}
where  $\beta_{jik}$ are (mollified) coefficients chosen to eliminate linear with respect to $(p_1,\ldots,p_d)$ terms in  $\{p_k,q_j\}$. Consider function
\begin{multline}
\upchi\Bigl( \frac{\mu}{T} \langle \ell, q-\bar{q}\rangle +
\epsilon\varsigma \frac{t}{T} \Bigr)=\\
\upchi\Bigl( \frac{\mu}{T} \bigl(\langle \ell', x-\bar{x}'\rangle - \sum_{j,k}\ell_j \upphi_{jk}p_k + 
\frac{1}{2}\sum_{i,j, k}\ell_j \beta_{jik}p_ip_k \bigr)\Bigr)
\label{19-2-70}
\end{multline}
and use a vector scale $(\boldgamma,\boldrho)$ with respect to $(x,\xi)$ with
\begin{equation}
\rho_j\gamma_j \ge C\mu ^{-1}h|\log h|.
\label{19-2-71}
\end{equation}
Consider the first derivatives of the argument of $\upchi$; we want them to be bounded \emph{after rescaling\/}.

Let us take $\gamma_j= \gamma $, $\rho_j=\rho $ as $j=1,\ldots, r$ and $\gamma_j= \rho $, $\rho_j= \gamma $ as $j=r+1,\ldots,2r$ with $\rho \ge \gamma$; then $|\upphi_{jk}-\upphi_{jk}(y)|\le c\rho$.

Note that to have the first derivatives with respect to $x'$ bounded by $C_0\gamma^{-1}$ one needs to take $\gamma = \mu^{-1}T$\,\footnote{\label{foot-19-14} Or less; we always take the maximal value.} as $j=1,\ldots, r$.

Further, to have the first derivatives with respect to $\xi'$ bounded by $C_0\rho^{-1}$ one needs to take
$\rho = (\mu^{-1}T)^{\frac{1}{2}}$ as $j=1,\ldots, r$. Then the first derivatives with respect to $\xi''$ are also bounded by $C_0\gamma^{-1}$ and to have the first derivatives with respect to $x''$ bounded by $C_0\rho^{-1}$ one needs to have $|\partial_j p_k|\le c\rho $ as $j=1,\ldots, r$, $k=j+r$ or, equivalently, $\partial _jV_k(y)=0$ for indicated $j,k$; this could be achieved by a gauge transformation.

To satisfy (\ref{19-2-71}) one needs then
\begin{equation}
T\ge   C\mu (\mu ^{-1}h|\log h|)^{\frac{2}{3}};
\label{19-2-72}
\end{equation}
one can see easily that $T=\bar{T}=\epsilon \mu^{-1}$ fits if and only if
$\mu $ satisfies (\ref{19-2-60}). Then picking up
$\gamma= C\min \bigl(\mu h|\log h|,\mu^{-1}T\bigr)$ and
$\rho =\gamma^{\frac{1}{2}}$ we see that (\ref{19-2-71}) holds and (\ref{19-2-70}) is an admissible symbol.

Then our standard analysis (like in theorem~\ref{book_new-thm-2-3-1} of \cite{futurebook}, together with  propositions \ref{prop-19-2-8}--\ref{prop-19-2-9} and corollary \ref{cor-19-2-11}\footref{foot-19-13} imply estimate (\ref{19-2-63}).
\end{proof}

Applying the same arguments as in (i), (ii) in the cases
$T\ge C(\mu h|\log h|)^{\frac{1}{2}}$ and $T\le C(\mu h|\log h|)^{\frac{1}{2}}$
respectively one can prove easily

\begin{proposition}\label{prop-19-2-13}\footnote{\label{foot-19-15} Weak magnetic field approach.}
Let $(l,\sigma)\succeq (1,2)$, microhyperbolicity condition be fulfilled and
\begin{gather}
 \epsilon_1 (h|\log h|)^{-\frac{1}{2}}\le \mu \le
\epsilon (h|\log h|)^{-1},\label{19-2-73}\\[3pt]
C'(\mu ^{-1}h|\log h|)^{\frac 1 2}\le \varepsilon \le  C'\mu h|\log h|, \label{19-2-74}\\[3pt]
T_*= C  \max\Bigl(\varepsilon^{-1} h|\log h|,\ 
\mu (\mu^{-1}h|\log h|)^{\frac{2}{3}}\Bigr)\qquad\label{19-2-75}
\end{gather}
where large enough constant $C'$ is chosen in the very last moment. Then

\medskip\noindent
(i) For $T\in [T_*,T^*]$, $\tau \in [-\epsilon_1,\epsilon_1]$ estimate \textup{(\ref{19-2-63})} holds where $T^*=\epsilon$; 

\medskip\noindent
(ii) Furthermore, under $\fN$-microhyperbolicity this estimate holds with $T^*$  defined  as in proposition~\ref{prop-19-2-12}(ii).
\end{proposition}

Then we arrive immediately to

\begin{corollary}\label{cor-19-2-14} 
In the framework of propositions \ref{prop-19-2-12}, \ref{prop-19-2-13}
\begin{equation}
|F_{t \to h^{-1}\tau} \bigl(\bar{\chi}_T(t) - \bar{\chi}_{T'}(t)\bigr)\Gamma \psi U | \le Ch^s
\label{19-2-76}
\end{equation}
with $T_*\le T'\le T\le T^*$ with $T_*$ defined by \textup{(\ref{19-2-75})} and $T^*$ defined in proposition~\ref{prop-19-2-12}(i), (ii). 
\end{corollary}
Recall that $\bar{\chi}$ is admissible function supported in $[-1,1]$ and equal $1$ at $[-\frac{1}{2}, \frac{1}{2}]$.

\section{Tauberian theory}
\label{sect-19-2-3}

\begin{corollary}\label{cor-19-2-15}
In the framework of proposition \ref{prop-19-2-12}
\begin{equation}
|F_{t \to h^{-1}\tau} \bar{\chi}_T(t) \Gamma (\psi U) | \le Ch^{1-d}\qquad \forall \tau:|\tau|\le \epsilon
\label{19-2-77}
\end{equation}
and
\begin{multline}
\R^\T\Def |\Gamma (\psi \tilde{e}) (\tau)- h^{-1} \int_{-\infty}^\tau
\Bigl(F_{t \to h^{-1}\tau'} \bar{\chi}_T(t) \Gamma (\psi U)\Bigr)\,d\tau'| \le \\[2pt]
CT^{*\,-1}h^{1-d}\qquad \forall \tau:|\tau|\le \epsilon
\label{19-2-78}
\end{multline}
for any $T\in [\bar{T},T^*]$.
\end{corollary}

\begin{proof}
Rescaling $x\mapsto \mu x$, $\varepsilon \mapsto \mu\varepsilon$,
$t\mapsto \mu t$, $h\mapsto \mu h$, $T_0\mapsto \epsilon$  we arrive  to a standard Schr\"odinger operator (i.e. with $\mu =1$),  we can apply standard methods. Then estimate (\ref{19-2-77}) holds for $T=\bar{T}$. Combining with (\ref{19-2-76}) we arrive to (\ref{19-2-77}) with arbitrary
$T\in [\bar{T}, T^*]$.

Applying Tauberian arguments we arrive then to the Tauberian estimate (\ref{19-2-78}).
\end{proof}

\section{Main theorem}
\label{sect-19-2-4}

Now under assumption (\ref{19-2-60}) we need to calculate the Tauberian expression
\begin{equation}
h^{-1} \int_{-\infty}^\tau
\Bigl(F_{t \to h^{-1}\tau'} \bar{\chi}_{T}(t) \Gamma (\psi U)\Bigr)\,d\tau'.
\label{19-2-79}
\end{equation}
with $T= Ch|\log h|$.

Using rescaling to the standard Schr\"odinger operator we get that it is equal to
\begin{equation}
\sum_{m,n\ge 0} \kappa_{nm}h^{-d+2m+2n}\mu^{2n}
\label{19-2-80}
\end{equation}
as all other coefficients vanish under assumption (\ref{19-2-38}).

Taking only term $\kappa_{00}h^{-d}$ we get the standard Weyl expression
\begin{equation}
h^{-d}\int \cN^\W (\tau,x)\psi (x)\,dx
\label{19-2-81}
\end{equation}
with an error $O(\mu^2h^{2-d})$ which is $O(\mu^{-1}h^{1-d})$ as $\mu \le h^{-\frac{1}{3}}$.

Further, taking only terms $\kappa_{00}h^{-d}+\kappa_{10}h^{2-d}\mu ^2$ we commit an error $O(h^{2-d}+\mu^4 h^{4-d})$ which is $O(\mu^{-1}h^{1-d})$ as
$\mu \le h^{-\frac{3}{5}}$. One can calculate this correction term $\kappa_{10}h^{2-d}\mu ^2$ easily. More generally, taking $n$-terms we make an error $O(\mu^{2n}h^{2n-d})$ arriving to 
\begin{multline}
|\Gamma (\psi \tilde{e}) \tau)- 
h^{-d}\int \tilde{\cN}^\W_{(n)}(x,\tau)\psi (x)\,dx | \le CT^{*\,-1}h^{1-d}+ C\mu^{2n}h^{2n-d}\\
\forall \tau:|\tau|\le \epsilon
\label{19-2-82}
\end{multline}
where $h^{-d}\cN^\W_{(n)}$ denotes $n$-term Weyl expression i.e. cut the asymptotic sum
\begin{equation}
h^{-d}\cN^\W_{(\infty)} \Def \sum_{n\ge 0} \varkappa_{n0}(x) h^{-d}(\mu h)^{2n}.
\label{19-2-83}
\end{equation}

Instead let us notice that
$\kappa_{n0}=\int \varkappa_{n0}(x,\tau)\psi (x)\, dx$ where $\varkappa_{n0}(x,\tau)$ depends only on $g^{jk},V$ and $\partial_kV_j$ calculated at point $x$, and thus with indicated error to calculate
(\ref{19-2-79}) with $\Gamma_y$ instead of $\Gamma$ and with $\psi=1$ and one can consider Schr\"odinger operator
\begin{multline}
\bar{A}= \sum_{j,k} \bar{g}^{jk}  \bar{P}_j\bar{P}_k+\bar{V},\qquad
\bar{g}^{jk}=g^{jk}(y), \bar{V}=V(y), \\
\bar{P}_j=hD_j-V_j(y)-\sum_k(\partial_kV_j)(y)(x_k-y_k)
\label{19-2-84}
\end{multline}
and for this operator expression in question is exactly $h^{-d}\tilde{\cN}^\MW(y,\tau)$.

So, under assumption (\ref{19-2-60}) we have proven
\begin{multline}
|\Gamma (\psi \tilde{e}) \tau)- h^{-d}\int \tilde{\cN}^\MW(x,\tau)\psi (x)\,dx | \le CT^{*\,-1}h^{1-d}\\
\forall \tau:|\tau|\le \epsilon
\label{19-2-85}
\end{multline}
where for $\mu \le Ch^{-\frac{1}{3}}$ one can replace $\tilde{\cN}^\MW$ by
$\tilde{\cN}^\W$.

Thus tasking into account that the difference between Weyl or magnetic Weyl expressions for the original and mollified operator does not exceed $C\varepsilon^l |\log \varepsilon|^{-\sigma} h^{-d}$ we arrive to

\begin{theorem}\label{thm-19-2-16}
Let assumptions \textup{(\ref{19-1-4})}--\textup{(\ref{19-1-6})}, $\textup{(\ref{19-1-25})}_{1-3}$ with \linebreak $(\bar{l},\bar{\sigma})\succeq(l,\sigma)\succeq (1,2)$ and \textup{(\ref{19-2-38})}  be fulfilled. Let 
$\mu  \le \mu^*_2 = C(h|\log h|)^{-\frac{1}{2}}$. Then there are two framing approximations\footref{book_new-foot-18-16} (see Chapter~\ref{book_new-sect-18} of \cite{futurebook}) such that the following statements are true:

\medskip \noindent
(i) Let  $\fN$-microhyperbolicity condition (see definition~\ref{def-19-2-5}) be fulfilled. Then if \underline{either}  $(l,\sigma)\succeq (2,0)$ \underline{or} $\#\fN=1$, then 
\begin{multline}
\R^\MW \Def
|\int \Bigl( \tilde{e}(x,x,0) - h^{-d}\cN^\MW (x,0)\Bigr)\psi (x)\, dx |\le\\[3pt]
C\mu ^{-1}h^{1-d} + C(\mu h)^l |\log h|^{l-\sigma} h^{-d}
\label{19-2-86}
\end{multline}
and if  $(l,\sigma)\prec (2,0)$ \underline{and} $\#\fN\ge 2$ then
\begin{equation}
\R^\MW \le 
C\mu ^{l-l} h^{1-d}|\log \mu |^{- \sigma } + 
C(\mu h)^l |\log h|^{l-\sigma} h^{-d};
\label{19-2-87}
\end{equation}

\medskip\noindent
(ii) Let  microhyperbolicity condition (see definition~\ref{def-19-2-4}) be fulfilled. Then
\begin{equation}
\R^\MW \le C h^{1-d} + C(\mu h)^l |\log h|^{l-\sigma} h^{-d};
\label{19-2-88}
\end{equation}
\end{theorem}

\begin{corollary}\label{cor-19-2-17}
If $(l,\sigma)\succeq (3,1)$ then in the framework of (i) asymptotics
\begin{equation}
\R^\MW \le C\mu ^{-1}h^{1-d}
\label{19-2-89}
\end{equation}
holds as $\mu \le C(h|\log h|)^{-\frac{1}{2}}$.
\end{corollary}

\begin{remark}\label{cor-19-2-18}
(i) The above estimates hold for $\R^\W_{(2)}$ as well;

\medskip\noindent
(ii) Theorem \ref{thm-19-2-16} can be extended to the case
$C (h|\log h|)^{-\frac{1}{2}} \le \mu \le h^{\delta-1}$ but it will  be our theorem of choice only as $\mu \le h^{-\frac{1}{2}}$;

\medskip\noindent
(iii) In the case $C (h|\log h|)^{-\frac{1}{2}} \le \mu \le h^{\delta-1}$ we need to take $n$ such that $(2n,0)\succeq (l,\sigma)$.

\medskip\noindent
(iv) This remark also applies to theorem~\ref{thm-19-2-19} below.
\end{remark}

The general case of weak magnetic field is completed. 

\section{Special case of constant $g^{jk}$, $F_{jk}$}
\label{sect-19-2-5}

Assume now that $g^{jk}$ and $F_{jk}$ are constant. Then without any loss of the generality one can assume that 
\begin{gather}
g^{jk}=\updelta_{jk},\qquad F_{jk}= f_j \updelta_{k,j+r}-f_k\updelta_{j,k-r}
\label{19-2-90}\\
\shortintertext{and}
V_j=0 \quad {as\ \ } j=1,\ldots,r\qquad V_j= f_{j-r}x_{j-r}\quad {as\ \ } j=r+1,\ldots,2r.\label{19-2-91}
\end{gather}
Then all the arguments of the previous Subsections could be easily simplified as we can assume without any loss of the generality that $\ell=(0,\ell'')$ and consider the shift with respect to $\xi''$. Then the logarithmic uncertainty principle would mean $ T \varepsilon \ge Ch|\log h|$ and thus we can take $\varepsilon = C\mu h|\log h|$ to accommodate $T=\bar{T}\asymp  \mu^{-1}$.

However we can make better than this and  allow  critical points. Let us first consider a point $\bar{x}$ such that at this point
\begin{equation}
|\nabla V| \asymp \nu \ge \bar{\nu} \Def 
C_0\max \bigl(\mu^{-1},\mu h|\log h|\bigr)
\label{19-2-92}
\end{equation}
where we pick $\bar{\nu}$ exactly as in the $2\D$-case.

Then $|\nabla V|\asymp \nu$ in $B(\bar{x},2\epsilon_1 \nu)$ and it follows from our analysis that the evolution (including microlocal one) which starts at $B(\bar{x},\epsilon_1 \nu)$  does not leave $B(\bar{x},2\epsilon_1 \nu)$ for time $T^*=\epsilon \mu$.

On the other hand, the shift in $\xi''$ is $\asymp T\nu$ and the logarithmic uncertainty principle requires $T\nu \cdot \varepsilon \ge Ch|\log h|$ and we can take
\begin{equation}
\varepsilon = C\left\{\begin{aligned}
&\mu h\nu^{-1}|\log h|\qquad && \nu \ge \bar{\nu},\\
&\mu h \bar{\nu}^{-1} \qquad && \nu \le \bar{\nu}
\end{aligned}\right.
\label{19-2-93}
\end{equation}
where $\nu \le \bar{\nu}$ means that $|\nabla V|\le C\bar{\nu}$ in $\bar{x}$ and its $2\bar{\nu}$-vicinity.

Then contribution of $B(x,\nu)$ to the Tauberian remainder does not exceed 
$C \mu^{-1}h^{1-d} \nu^d $ as $\nu \ge \bar{\nu}$ and $C\mu h^{1-d}\bar{\nu}^d$ as $\nu \le \bar{\nu}$ and in the latter case we take $T^*=\epsilon \mu^{-1}$. Summation over partition returns
\begin{equation}
\R^\T \le C\mu^{-1}h^{1-d} +
C\mu h^{1-d}\mes \bigl(\{x:\ |\nabla V(x)|\le \bar{\nu}\}\bigr)
\label{19-2-94}
\end{equation}
which in turn implies the same  estimate for $\R^\W_{(\infty)}$ calculated for mollified operator and under assumption (\ref{19-2-38}) an error when we pass from mollified to the original operator error becomes
\begin{equation}
Ch^{-d} \int \varepsilon ^l |\log \varepsilon | ^{-\sigma}\, dx \asymp
Ch^{-d} \int (\mu h) ^l |\log (\mu h/\nu)| ^{-\sigma} \nu^{-l} \,dx .
\label{19-2-95}
\end{equation}
Assuming that $(l,\sigma)\succeq (2,0)$ and imposing generic non-degeneracy condition
\begin{equation}
|\nabla V|\le 0\epsilon \implies |\det \Hess V |\ge\epsilon
\label{19-2-96}
\end{equation}
we conclude that expression (\ref{19-2-94}) does not exceed 
\begin{equation}
C\mu^{-1}h^{1-d}+ C\mu (\mu h|\log h|)^{d/2} h^{1-d}
\label{19-2-97}
\end{equation}
which even in the worst case $d=4$ does not exceed $C\mu^{-1}h^{1-d}$ as 
$\mu \le (h|\log h|)^{-\frac{1}{2}}$. 

Meanwhile as $l<d$ expression (\ref{19-2-95}) does not exceed 
$C(\mu h)^l |\log h|^{-\sigma} h^{-d}$ and as $l\ge d$ this expression does not exceed (\ref{19-2-97}) and therefore both expressions (\ref{19-2-94})  and (\ref{19-2-95}) do not exceed
\begin{equation}
C\mu^{-1}h^{1-d}+ C\mu (\mu h|\log h|)^{d/2} h^{1-d} + 
C(\mu h)^l |\log h|^{-\sigma} h^{-d}.
\label{19-2-98}
\end{equation}

Then we arrive to

\begin{theorem}\label{thm-19-2-19}
Let assumptions \textup{(\ref{19-1-4})}--\textup{(\ref{19-1-6})}, $\textup{(\ref{19-1-6})}_{3}$ with $(l,\sigma)\succeq (1,2)$ and  \textup{(\ref{19-2-38})} be fulfilled. Let $g^{jk}$, $F_{jk}$ be constant. 

Then there are two framing approximations\footref{book_new-foot-18-16} (see Chapter~\ref{book_new-sect-18} of \cite{futurebook}) such that the following statements are true:

\medskip \noindent
(i) Under microhyperbolicity assumption  $|\nabla V|\ge \epsilon$ estimate 
\textup{(\ref{19-2-86})} holds;

\medskip\noindent
(ii) Under assumptions $(l,\sigma)\succeq (2,0)$ and \textup{(\ref{19-2-96})} $\R^\MW$ does not exceed \textup{(\ref{19-2-98})}.
\end{theorem}

\chapter{Canonical form}
\label{sect-19-3}

\section{Pilot-model}
\label{sect-19-3-1}

Assume now that $g^{jk}$ and $F_{jk}$ are constant. Then without any loss of the generality we can assume that $V_j(x)$ are linear functions. 

Then $A^0$ is transformed into exactly
\begin{equation}
\sum_{1\le j\le r} f_j\Bigl(h^2 D_j^2+\mu ^2 x_j^2\Bigr)
\label{19-3-1}
\end{equation}
by $\mu ^{-1}h$-metaplectic transformation which consists of the following steps:

\medskip\noindent
(i) Change of variables $(x,\mu^{-1}hD) \mapsto (Qx, \,^t\!Q^{-1}\mu^{-1}hD)$ transforming $g^{jk}$ into $\updelta_{jk}$ and $F$ into matrix with $F_{j,j+r}=f_j$, $F_{j+r,j}=-f_j$ and other elements $0$. It transforms $V(x)$ into $V(Qx)$. 

\medskip\noindent
(ii) Gauge transformation (multiplication by $e^{i\mu h^{-1}S(x)}$ with quadratic form $S(x)$, transforming $A_j(x)$ into $0$  and $A_{j+r}(x)$ into $f_j x_j$ for $j=1,\ldots,r$\,\footnote{\label{foot-19-16} Thus we achieve (\ref{19-2-90})--(\ref{19-2-91}).}. Then $A^0$ is transformed into
\begin{equation}
h^2 \sum_{1\le j\le r} \Bigl( D_j^2 + \bigl(D_{j+r}- f_j\mu h^{-1}x_j\bigr)^2\Bigr)
\label{19-3-2}
\end{equation}
and $V(x)$ is preserved.

\medskip\noindent
(iii) Partial $\mu^{-1}h$-Fourier transform:
\begin{equation*}
(x', x''; \mu ^{-1}h D', \mu ^{-1}h D')\mapsto
(x', -\mu ^{-1}h D''; \mu ^{-1}h D', x'')
\end{equation*}
transforming $A^0$ into
\begin{equation}
h^2 \sum_{1\le j\le r}
\Bigl( D_j^2 + \mu ^2 h^{-2} \bigl(x_{j+r}-f_jx_j\bigr)^2\Bigr)
\label{19-3-3}
\end{equation}
and $V(x)$ into $V(x', -\mu^{-1}hD'')$.

\medskip\noindent
(iv) Change of variables 
\begin{equation*}
(x', x''; \mu ^{-1}h D', \mu ^{-1}h D')\mapsto
(x'-Kx'',x''; \mu^{-1}hD', \mu ^{-1}hD'' + \,^t\!K \mu^{-1}hD')
\end{equation*}
with $K_{j,j+r}=f_j\updelta_{j,k-r}$, transforming $A^0$ into
\begin{equation*}
\sum_{1\le j\le r} \Bigl(h^2 D_j^2 +  f_j^2\mu^2x_j^2\Bigr)
\end{equation*}
and $V$ into 
$\tilde{V}\Def V(x'-Kx'', \mu ^{-1}h D''- \,^t\!K \mu ^{-1}hD')$ with Weyl quantization. Finally, $x_j\mapsto f_j^{\frac{1}{2}}x_j$, 
$D_j\mapsto f_j^{-\frac{1}{2}}D_j$ reduces operator $A^0$ into (\ref{19-3-1}) 
and transforms $\tilde{V}$ accordingly.

This example already demonstrates why we need $\varepsilon \gtrsim \mu^{-1}$: we  need to fulfill equation similar to one considered in Section~\ref{book_new-sect-18-3} of \cite{futurebook}
\begin{equation}
\sum_j f_j \partial _{\phi_j} L \approx  \tilde{V}-W 
\label{19-3-4}
\end{equation}
where $L$ is a symbol of the operator generating transformation, $\tilde{V}$ is defined above, $W$ is its replacement in the canonical form and $(\rho_j,\phi_j)$ are polar coordinates in $(x_j,\xi_j)$-plane.

But then if $f_1,\ldots,f_r$ are commensurable, for generic $\tilde{V}$ one cannot satisfy this equation with an error better than $O(\mu^{-l}|\log \mu|^{-\sigma})$ and with $L$ which is $2\pi$-periodic with respect to all arguments $\phi_j$ and $W=W(\rho_1,\ldots, \rho_r)$ as $W$ would need to depend on some linear combinations of $\phi_j$. 

\begin{Problem}\label{Problem-19-3-1}
Explore what is possible when $f_1,\ldots,f_r$ are non-com\-men\-surable, and we try to satisfy (\ref{19-3-4}) with a better error. It definitely would depend on how these numbers are non-commensurable (something related to their diophantine properties - as $r=2$ it is related to Liouville's exponent for $f_1/f_2$).
\end{Problem}

\section{General case: framework}
\label{sect-19-3-2}

Now we need to reduce our operator to a canonical form assuming that
$(\bar{l},\bar{\sigma})\succeq (2,1)$ and
\begin{equation}
\mu_1^*\Def \epsilon_1(h|\log h|)^{- \frac{1}{2}}\le \mu \le 
\mu_2^*\Def \epsilon (h|\log h|)^{-1}
\label{19-3-5}
\end{equation}
(case $\mu  \ge  \mu_2^*$ and especially case $\mu \ge \mu_3^*\Def \epsilon h^{-1}$ we consider later). As we mentioned we will need to assume that
\begin{gather}
\varepsilon \ge C_0\mu ^{-1}
\label{19-3-6}\\
\shortintertext{rather than}
\varepsilon \ge C(\mu ^{-1}h|\log h|)^{\frac{1}{2}}
\label{19-3-7}
\end{gather}
as in Chapter~\ref{book_new-sect-18} of \cite{futurebook}. This larger $\varepsilon$ makes certain parts of our construction much simpler but leads to a larger error. In this and following Sections  we consider $\mu ^{-1}h$-quantization.

\section{Reducing the main part}
\label{sect-19-3-3}

First we need to reduce the main part. Consider point $x$ in the vicinity of
$\bar{x}$ and consider Hamiltonian fields $H_{\eta_j}$ of $\eta_j$ reduced to $\Sigma^0_x$ where we recall that $\eta_j$ are linear combinations of $p_j=\xi_j-V_j$: $\eta _{2k-1}\Def \Re \zeta_k$, $\eta _{2k}\Def \Im \zeta _k$, $k=1,\ldots,r$ and $\zeta_k$ are defined by (\ref{19-2-13})--(\ref{19-2-14}) (rather than redefined  by (\ref{19-2-22})--(\ref{19-2-23})). Also recall that $\Sigma^0=\{\eta_1=\ldots=\eta_{2r}=0\}$. 

Consider symplectic map $e^{t H_q}$ where $q= q(x; \eta_1,\dots, \eta_d)$ is a quadratic form
\begin{equation}
q (x;  \eta_1,\ldots, \eta_d)= \frac{1}{2}\sum _{i,j} q_{ij}(x)\eta_i\eta_j
\label{19-3-8}
\end{equation}
with $q_{ij}\in \sF^{\bar{l},\bar{\sigma}}$ where $\sF^{l,\sigma}$ are defined in remark~\ref{book_new-rem-18-3-2} of \cite{futurebook}. Note that
\begin{equation}
\frac {d\ }{dt} (\eta_k\circ e^{t H_q}) =  \{q,\eta_k\}\circ e^{tH_q}
\equiv \sum _{i,j} q_{ij}\{\eta_i,\eta _k\} \eta _i \circ e^{tH_q}
\label{19-3-9}
\end{equation}
modulo quadratic form. Then  using arguments of Subsection \ref{book_new-sect-18-3-1} of \cite{futurebook} one can prove easily that
\begin{equation}
\eta_k\circ e^{H_q} =\sum_j \beta _{kj}(x) \eta_j +
\sum _{i,j} \beta'_{kij} (x,\xi)\eta_i\eta_j
\label{19-3-10}
\end{equation}
with $\beta_{jk}\in \sF^{\bar{l},\bar{\sigma}}$,
$\beta'_{kij}\in \sF^{\bar{l}-1,\bar{\sigma}}$. Here
$\cB \Def (\beta _{jk})=e^{\cQ\Lambda}$  with $\cQ=(q_{kj})$ and 
$\Lambda = (\Lambda _{jm})$ where $\Lambda_{jm}=1$ as $j=2i-1,m=2i$, $\Lambda_{jm}=-1$ as $j=2i,m=2i-1$ and $\Lambda_{jm}=0$ in all other cases. Recall that $\{\eta_j,\eta_m\}\bigr|_{\Sigma^0}=\Lambda_{jm}$.

Obviously, 
\begin{claim}\label{19-3-11}
For a given matrix $\cB\in \sF^{\bar{l},\bar{\sigma}}$ one can find a symmetric matrix $\cQ=(q_{ij})$ of the same regularity such that $\cB =e^{\cQ\Lambda}$  if and only if  $\cB$ is a \emph{symplectic matrix\/} i.e.
\begin{equation}
\cB ^\dag \Lambda\cB=\Lambda.
\label{19-3-12}
\end{equation}
\end{claim}
In particular, one can transform (modulo quadratic forms): $\eta_1$ into  $p_1$,  $\eta_2$ into $\beta _{22} p_2$
with  $\beta_2$ disjoint from $0$\,\footnote{\label{foot-19-17} Since $\{p_1,p_k\}$ is  disjoint from 0 for some $k$ and we just rename $x_k$ into $x_2$ and v.v.}, and $\eta_3$ into
$p_{33}+\beta_{31}p_1+\beta_{32}p_2$, and $\eta_4$ into
$\beta_{44} p_3+\beta_{41}p_1+\beta_{42}p_2$ etc:
\begin{align}
&\cB\eta_m= \eta_m'\Def \hphantom{\beta_{mm}}p_m+\sum_{j\le 2k-2} \beta_{mj}p_j \quad 
&&\text{as\ \ }m=2k-1,\label{19-3-13}\\
&\cB\eta_m=\eta_m'\Def  \beta_{mm} p_m+\sum_{j\le 2k-2}\beta_{mj} p_j\quad &&\text{as\ \ }m=2k\label{19-3-14}
\end{align}
where coefficients $\beta_{jk}\in \sF^{\bar{l},\bar{\sigma}}$ are chosen to satisfy 
\begin{equation}
\{\eta'_j,\eta'_m\}=\Lambda_{jm}.
\label{19-3-15}
\end{equation}

Now let us  consider $\mu^{-1}h$-differential operator $Q=q^\w$ and transformation $T(t)=e^{it\mu h^{-1} Q}$ (``poor man's Fourier integral operator'').

Using arguments  of Subsection~\ref{book_new-sect-18-3-1} of \cite{futurebook} one can prove easily that
\begin{equation}
T(1) \eta_k ^\w T(-1) \equiv (\eta'_k)^\w+
\Bigl(\sum _{i,j} \beta '_{kij}\eta '_i\eta '_j\Bigr)^\w
\label{19-3-16}
\end{equation}
modulo operator with the norm not exceeding $C\mu ^{-2}h$. Then for
\begin{gather}
A^0=\mu^2\Bigl(\sum _{jk} a_{jk}(x)\, \eta_j\eta_k\Bigr)^\w
\label{19-3-17}\\
\shortintertext{we have}
T(1) A^0 T(-1) \equiv
\mu ^2 \Bigl(\sum _{j,k}  a_{jk}(x)\, \eta'_j\eta'_k \Bigr)^\w +
\mu ^2 \Bigl(\sum _{i,j,k}  b'_{ijk}\, \eta'_i \eta'_j\eta'_k \Bigr)^\w
\label{19-3-18}
\end{gather}
modulo operator with the norm not exceeding $C\mu^{-1}h$ where
$b'_{ijk}\in \sF^{\bar{l}-1,\bar{\sigma}}$; for simplicity of notations we  skip $'$ in what follows.

Also, the same arguments show that
\begin{equation}
T(-1)V(x)T(1) \equiv V(x) +\sum _j (b''_j\, \eta_j)^\w
\label{19-3-19}
\end{equation}
modulo operator with the norm not exceeding
$C\mu ^{-2} + C\mu ^{-l}|\log h|^{-\sigma}$ where
$b''_j\in \sF^{l-1,\sigma}$.

Recall that while $a_{jk}$, $V$ are functions depending on $x$ only, $b'_{ijk},b''_j$ are complete symbols. In the end of the day we decompose them in the powers of $\eta_1,\dots,\eta_{2r}$ as far as their smoothness allows.

Due to construction $\eta_1=p_1=\xi_1 - V_1(x)$\,\footnote{\label{foot-19-B-16} Recall that $\eta_k$ means $\eta'_k$ defined by (\ref{19-3-13}).}. Note that we can  get rid of $V_1$  using gauge transform, after which condition \ref{19-1-25-2} still holds. 

\begin{remark}\label{rem-19-3-2}
Exactly here we need this condition  \ref{19-1-25-2} rather than more natural 
condition $V_j  \in\sC^{\bar{l}+1,\bar{\sigma}}$. Otherwise we would be forced to confine ourselves to $C\varepsilon$-vicinity of $\bar{x}$ which would cause problems in Section~\ref{sect-19-4}.
\end{remark}

Then according to commutation property $\partial_{x_1}\eta_2 =1$ on $\Sigma^0$. Further, by its construction $\eta_2$ does not depend on $\xi_1,\xi_3,\ldots, \xi_d$ and therefore
\begin{equation}
\eta_2 = \alpha (x,\xi_2)
\bigl(x_1 - \lambda (x '  ,\xi_2)\bigr), \qquad
\alpha \bigr|_{\Sigma^0}=1
\label{19-3-20}
\end{equation}
with $\alpha \in \sF^{\bar{l},\bar{\sigma}}$, 
$\lambda \in \sF^{\bar{l}+1,\bar{\sigma}}$ and $x'=(x_2,\dots,x_d)$. 

Let us redefine
\begin{equation}
\eta_2\Def \bigl(x_1 - \lambda (x '  ,\xi_2)\bigr).
\label{19-3-21}
\end{equation}

\begin{claim}\label{19-3-22}
Without any loss of the generality we assume that $\bar{x}=0$; otherwise
in what follows one should replace $x$ by $x-\bar{x}$.
\end{claim}

Now we want to transform $\eta_2$ into $x_1$. To achieve this goal let us consider  $T'(t)= e^{i\mu h^{-1}(\xi_1\lambda )^w}$. Then we have the series of exact equalities
\begin{multline}
T'(-t)D_1T(t)= D_1, \quad T'(-t)\lambda^\w T'(t)=\lambda^\w,\\
T'(-t)x_1T'(t) = x_1 + t \lambda ^\w.
\label{19-3-23}
\end{multline}
Therefore $T'(1)$ will transform ``new'' $\eta_2^\w$ into $x_1$. Now we need to check how this transformation will affect $\eta_j^\w$ with $j\ge 3$ and also $\alpha^\w$, $b_{jk}^\w $, $b^{\prime \w}_{ijk}$, and finally $V$ and $(b''_j\eta_j)^\w$.

Note that Hamiltonian function $p_1\lambda $ belongs to
$\sF^{\bar{l}+2,\bar{\sigma}}$\,\footnote{\label{foot-19-19} Really, $\lambda \in \sF^{\bar{l}+1,\sigma}$ and multiplication by $p_1$ increases regularity by $1$ due to inequalities $|p_j|\le c\mu^{-1}\le  \varepsilon$ due to (\ref{19-2-2}).}. Therefore
Hamiltonian map $\Phi_t= e^{t H_{p_1\lambda}}$ belongs to $\sF^{\bar{l}+1,\bar{\sigma}}$ as $t\in [0,1]$ and then
\begin{multline}
\|T'(-t)\eta_j^\w T(t) - (\eta_j\circ \Phi_t )^\w \|  \le \\[2pt]
C(\mu ^{-1}h)^2
\Bigl(1+\varepsilon ^{\bar{l}+2-3}|\log \varepsilon |^{-\bar{\sigma}} \Bigr)\times
\Bigl(1+\varepsilon ^{\bar{l}+1-3}|\log \varepsilon |^{-\bar{\sigma}}\Bigr)
\le C\mu ^{-2}h,
\label{19-3-24}
\end{multline}
\begin{multline}
\| T'(-t)b_{jk}^\w T(t) - (b_{jk}\circ \Phi_t )^\w \| \le\\[2pt]
C(\mu ^{-1}h)^2
\Bigl(1+\varepsilon ^{\bar{l} -3}|\log \varepsilon |^{-\bar{\sigma}}\Bigr) \times
\Bigl(1+\varepsilon ^{\bar{l}+1-3}|\log \varepsilon |^{-\bar{\sigma}}\Bigr)
\le C\mu ^{-1}h,
\label{19-3-25}
\end{multline}
\begin{multline}
\| T'(-t)b_{ijk}^\w T(t) - (b_{ijk}\circ \Phi_t)^\w \| \le \\[2pt]
C(\mu ^{-1}h)^2
\Bigl(1+\varepsilon ^{\bar{l} -4}|\log \varepsilon |^{-\bar{\sigma}}\Bigr) \times
\Bigl(1+\varepsilon ^{\bar{l}+1-3}|\log \varepsilon |^{-\bar{\sigma}}\Bigr)
\le C h,
\label{19-3-26}
\end{multline}
Furthermore,  (\ref{19-3-25}) also holds for $\alpha$ (coefficient in (\ref{19-3-20})) instead of $b_{jk}$.

Therefore 
\begin{claim}\label{19-3-27}
Transformation of $A^0$ given by (\ref{19-3-17}) by $T'(1)$ leads to the same expression (\ref{19-3-17}) but all the symbols $\eta'_j$, $b'_{jk}$ and  $b'_{ijk}$ are replaced by  $\eta'_j\circ \Phi_1$,  $b'_{jk}\circ \Phi_1$ and  $b'_{ijk}\circ \Phi_1$ respectively and the total error does not exceed $C\mu^{-1}h$.
\end{claim}

\begin{remark}\label{rem-19-3-3}
Now $\eta_1$, $\eta_2$  are redefined as $p_1\circ \Phi_1=\xi_1$, $p_2\circ \Phi_1=x_1$ respectively.
\end{remark}

Let us consider $\eta_j$ with  $j\ge 3$. Note that $\{p_1,\eta_j\}=0$ on $\Sigma^0$  and therefore
\begin{equation}
\eta_j \circ \Phi _t =\eta_j + k_j\eta_1+ \sum _{k,i}\gamma_{jik}p_ip_k.
\label{19-3-28}
\end{equation}
The same arguments hold for $x_j$ with $j\ge 2$. Therefore, moving  errors arising in the quadratic part of $A^0$ into its cubic part, we can replace
$\eta_j \circ \Phi _t $ by $\eta_j + k_jp_1$ for $j\ge 3$, and $x_j\circ \Phi_t$ by $x_j$ for $j\ge 2$; moreover, in $b_{jk}$ we can replace
$x_1\circ \Phi_1 = x_1 -\lambda(x',\xi_2) $ just by $-\lambda (x',\xi_2)$.

Further, from the beginning $\eta_j$ with $j\ge 3$ could depend on $\xi_1$ (see (\ref{19-3-13}), (\ref{19-3-14})). After above transformation $\eta_j\circ \Phi_1$ could acquire ``more'' $\xi_1$. Note however that originally 
$\{\eta_j ,\eta_2\}=0$ on $\Sigma^0$ and $\eta_2\circ \Phi_1=x_1$; therefore
$\{\eta_j\circ \Phi_1, x_1\}=0$ on $\Sigma^0$ and therefore in $\eta_j\circ \Phi_1$, which are linear combination of  $p_k\circ \Phi_1$, we can replace  $p_k\circ \Phi_1$ by $\xi_k - V_k(0,x'')$ as $k\ge 3$  (modifying again cubic terms) where temporarily $x''=(x_2,\ldots,x_d)$. 

Thus we arrive to the same quadratic expression but with $p_1=\xi_1$, $p_2=x_1$, $p_k=\xi_k - V_k(x_2,\ldots, x_d;\xi_2)$, 
$b_{jk}=b_{jk}(x_2, x_3,\ldots, x_d;\xi_2)$.

\medskip
Repeating this process in the end of the day we arrive to the same quadratic expression with $\eta_{2k-1}=\xi_k$, $\eta_{2k}=x_k$ $(k=1,\ldots, r)$,
$b_{jk}=b_{jk}(x'',\xi'')$ where now $x''\Def (x_{r+1},\ldots, x_{2r})$,  $\xi''\Def(\xi_{r+1},\ldots, \xi_{2r})$:
\begin{equation}
\mu^2  \sum_{j,k} b_{jk}(x'',\mu^{-1}hD'') \, \eta_j^\w \eta_k^\w +
\mu^2\Bigl(\sum_{j,k,m} b_{jkm}\, \eta_j\eta_k\eta_m\Bigr)^\w
\label{19-3-29}
\end{equation}
with  $b_{jkm}=b_{jkm}(x,\xi) \in \sF^{\bar{l}-1,\bar{\sigma}}$.

Also, considering transformation of $V(x)$ we see that at each step the semiclassical error does not exceed
$C(\mu ^{-1}h)^2 \varepsilon ^{l-3}|\log h|^{-\sigma}\le C\mu ^{-1}h$
and in the end of the day we arrive to
\begin{equation}
W(x'',\mu^{-1}hD'') +\Bigl(\sum_j b_j \, \eta_j\Bigr)^\w
\label{19-3-30}
\end{equation}
with $W\in \sF^{l,\sigma}$ and $b_j=b_j (x,\xi) \in \sF^{l-1,\sigma}$.

Therefore we have proven

\begin{proposition}\label{prop-19-3-4} 
Let conditions \textup{(\ref{19-1-4})},  \textup{(\ref{19-3-5})} and \textup{(\ref{19-3-6})} be fulfilled, $\bar{x}$ be a fixed point and
$\varepsilon \le R \le \epsilon_2$ with a small enough constant $\epsilon_2$. Then there exists  a bounded operator $\cT$ such that

\medskip\noindent
(i) For operators $Q=q^\w$  with symbol $q$ supported in 
$B(0,R)\subset \bR^{2d}$ and  $Q'=q_1^{\prime\,\w}$  with $q'_1$ supported in $B((\bar{x},0),2C_1R)\subset \bR^{2d}$ and equal $1$ in
$B((\bar{x},0),C_1R)\subset \bR^{2d}$ the following equalities hold modulo negligible operators
\begin{gather}
(I-Q'_1 )\cT Q \equiv 0, \label{19-3-31}\\[2pt]
\cT ^* \cT  Q \equiv Q\label{19-3-32}\\
\intertext{and modulo operators with norm not exceeding $C\mu ^{-1}h$}
\cT ^* \psi \cT \equiv {\bar \psi}^\w,\label{19-3-33}\\[2pt]
\cT ^* A\cT  Q \equiv \cA  Q\label{19-3-34}
\end{gather}
with
\begin{multline}
\cA  =  \sum _{i,j}  b_{ij}(x'',\mu ^{-1}hD'')
\bar{P}_i\bar{P}_j+
\mu^2 \Bigl(\sum _{i,j,k}  b_{ijk} \bar{p}_i\bar{p}_j \bar{p}_k\Bigr)^\w +\\
b_0(x'',\mu ^{-1}hD'')  + \Bigl(\sum _j  b_i \bar{p}_i\Bigr)^\w
\label{19-3-35}
\end{multline}
where $b_{ij}\in \sF^{\bar{l},\bar{\sigma}}$,
$b_{ijk}=b_{ijk}(x,\xi) \in \sF^{\bar{l}-1,\bar{\sigma}}$,
$b_0\in \sF^{l,\sigma}$, $b_i=b_i(x,\xi)\in \sF^{l-1,\sigma}$ are real-valued,
\begin{equation}
\bar{p}_{2k-1}=\xi_k,\quad \bar{p}_{2k}=x_k, \qquad k=1,\ldots, r,
\label{19-3-36}
\end{equation}
$\bar{P}_i=\mu {\bar p}_i^\w$,  $x'=(x_1,\dots, x_r)$,
$x''=(x_{r+1},\dots, x_{2r})$ etc; recall that $\psi=\psi (x)$ is a smooth function;

\medskip\noindent
(ii) For operators $Q'=q^{\prime\,\w}$  with symbol $q'$ supported in
$B((\bar{x},0),R)\subset \bR^{2d}$ and  $Q_1=q_1^\w$  with symbol $q_1$ supported in $B(0,2C_1R)\subset \bR^{2d}$ and equal $1$ in
$B(0,C_1R)\subset \bR^{2d}$
\begin{equation}
\|Q' \cT (1-Q_1  )\| \le C h^s.
\label{19-3-37}
\end{equation}
\end{proposition}

\begin{remark}\label{rem-19-3-5}  
(i) Calculations show that
\begin{equation}
b_{ij}= b'_{ij}\circ \Psi_0, \qquad b_0 =  V\circ \Psi_0, \quad
\bar{\psi} = \psi \circ \Psi
\label{19-3-38}
\end{equation}
where $\Psi $ is a symplectic map
$\bR^{2d}\supset B(0, R) \to B((\bar{x},0), 2C_1R)$ and
$\Psi_0= \Pi _x \Psi \bigr|_{\bar{\Sigma}^0}$ with
$\bar{\Sigma}^0\Def \{x'=\xi'=0\}$ and
\begin{equation}
|\det D\Psi _0 | = f_1\cdots f_r;
\label{19-3-39}
\end{equation}
(ii) Recall that one can rewrite quadratic part as
\begin{equation}
\mu^2 \sum_{ \fm \in \fM }\ \sum_{j,k\in  \fm } a_{jk}(x'',\mu^{-1}hD'')\, Z_j^*Z_k- \mu h\sum_k  a_{kk}(x'',\mu^{-1}hD'')
\label{19-3-40}
\end{equation}
with $Z_k= \mu^{-1}hD_k - i x_k$, $Z^*_k= \mu^{-1}hD_k + i x_k$;

\medskip\noindent
(iii) In particular, if $\#\fm=1$ for all $\fm\in \fM$ (no $2$-nd order resonances) one can rewrite quadratic form as
\begin{equation}
\sum_{1\le k\le r}  
f_j(x'',\mu^{-1}hD'')\,  \bigl( h^2D_k^2 + \mu^2 x_k^2\bigr).
\label{19-3-41}
\end{equation}
\end{remark}

\section{Reducing next terms}
\label{sect-19-3-4}

There is very little what we can do else because of all the resonances. We already finished the transformation of the ``main part''
\begin{equation}
\cA _0\Def  \sum _{i,j}  b_{ij}(x'',\mu ^{-1}hD'')
\bar{P}_i\bar{P}_j+ b_0(x'',\mu ^{-1}hD'')
\label{19-3-42}
\end{equation}
where one can rewrite quadratic part as (\ref{19-3-40}); in what follows the quadratic part and the potential of operator generate corrections  which could be written in the form
\begin{gather}
\mu ^2 \sum_{\alpha :\,3\le |\alpha|\le \lceil (\bar{l},\bar{\sigma})\rceil +1}
b''_{\alpha} (x'',\mu ^{-1}hD'')\,  (\bar{p}^\alpha)^\w \label{19-3-43}\\
\shortintertext{and}
\sum _{\alpha : 1\le |\alpha|\le  \lceil ( l,\sigma)\rceil -1}
b'''_{\alpha} (x'',\mu ^{-1}hD'')\,  (\bar{p}^\alpha)^\w \label{19-3-44}
\end{gather}
modulo terms with norms not exceeding 
$C\mu ^{-\bar{l}}|\log \mu |^{-\bar{\sigma}}$ and 
$C\mu ^{-l}|\log \mu |^{-\sigma}$ respectively  with
$b''_\alpha \in \sF^{\bar{l} +2-|\alpha|,\bar{\sigma}}$ and
$b'''_\alpha \in \sF^{l-|\alpha|,\sigma}$. 

Here and below $\lceil (l,\sigma)\rceil = \lceil l\rceil$ unless 
$l\in \bZ, \sigma >0$ in which case $\lceil (l,\sigma)\rceil = l+1$ etc.

The major problem are resonances. The $2$-nd order resonances ($f_i=f_j$ as $i\ne j$) prevent us from the diagonalization of the quadratic form as in (\ref{19-3-41}). We can only rewrite quadratic form as (\ref{19-3-40}).

Now it is easy to get rid of non-resonant terms of the $3$-rd order in the main part and of $1$-st order terms in $V$. Let us consider transformation
\begin{equation}
U= \Bigl(e^{i\mu ^{-1}h^{-1} \bigl(\mu ^2 S_3 (x'',\xi'';x', \xi ' )+
S_1(x'',\xi''; x',\xi ' )\bigr)  }\Bigr)^\w
\label{19-3-45}
\end{equation}
which is a pseudo-differential operator as
$\mu \ge C h^{-\frac{1}{3}}|\log h|^{\frac{1}{3}}$. As a result the principal part $\cA _0$ does not change, and the next terms 
\begin{equation}
\cA _1 \Def \mu^2a_1^\w  \qquad 
a_1\Def \sum _{\alpha :|\alpha|= 3} b''_{\alpha} (x'',\xi '')\, \bar{p}^\alpha +
\sum _{\alpha:|\alpha|=1} b'''_{\alpha} (x'',\xi '')\,  \bar{p}^\alpha
\label{19-3-46}
\end{equation}
(which have norms $O(\mu ^{-1})$) are replaced by
\begin{equation}
\cA _1 + \bigl\{\tilde{a}^0, \mu ^2 S_3 + S_1\bigr\} ^{\prime\,\w} + \ldots
\label{19-3-47}
\end{equation}
with $\tilde{a}^0 = \sum _{i,j} b_{ij}(x'',\xi'')\, p_i p_j$ and the ``short'' Poisson brackets $\{.,.\}'$ (i. e. with respect to $(x',\xi')$ only). 

Here dots denote terms of the same type (\ref{19-3-43}) with $|\alpha |\ge 4$ (modulo $O(\mu ^{-\bar{l}}|\log \mu |^{-\bar{\sigma}})$) and of the type (\ref{19-3-44}) with $|\alpha|\ge 2$ (modulo  
$O(\mu ^{-l}|\log \mu |^{-\sigma})$).

One can see easily that by an appropriate choice of $S_3$ and $S_1$ we can eliminate all terms in the ``main part'' 
$\mu^2 a _1 + \bigl\{\tilde{a}^0, \mu ^2 S_3 + S_1\bigr\} '$ save corresponding to the $3$-rd order resonances, namely
\begin{equation}
\Re \Bigl( \mu ^2 \sum _{m,j,k} b''_{mjk}(x'',\xi '') \,
\zeta _m^\dag \zeta_j\zeta_k\Bigr)
\label{19-3-48}
\end{equation}
where $\zeta_m = \xi _m + ix_m$, $m=1,\ldots,r$ as we can assume without any loss of the generality that
\begin{equation}
\{\tilde{a}^0,\zeta_m\}' = \sum _n if_{mn} \zeta_n
\label{19-3-49}
\end{equation}
where matrix $(f_{mn})$ has only positive eigenvalues $f_k$.

Note that  sum (\ref{19-3-48}) is restricted to indices $m,j,k$ such that $f_m$ is not disjoint from $f_j+f_k$; in particular $b''_{mjk}=0$ unless
$m,j,k\in \fn$ with $\fn\in \fN$.

Further, one can get rid of the non-resonant $4$-th order terms in (\ref{19-3-43}) and $2$-nd order terms in (\ref{19-3-44})  in the same way but while  the $2$-nd order resonant terms are only
\begin{enumerate}[label=(\roman*),leftmargin=*]
\item $\zeta_j^\dag \zeta_k$ with $j$, $k$ belonging to the same
$\fm$-group,

\hskip-30pt
the $4$-th order  resonant terms  include
%\begin{description}
\item $\zeta_j^\dag \zeta_k \zeta_m \zeta_n$ and 
$\zeta_j\zeta_k^\dag \zeta_m^\dag \zeta_n ^\dag  $ with
$f_j$ not disjoint from  $f_k+f_m+f_n$  

\hskip-3pt
and  also
\item
$\zeta_j^\dag \zeta_k^\dag \zeta_m \zeta_n$ with
$f_j+f+k$ not disjoint from $f_m+f_n$
\end{enumerate}
which we never intended to cover. The same is true for $5$-th order resonant terms as well etc.

Still we need to have some concern about only the $3$-rd, the $4$-th and the $5$-th order terms in (\ref{19-3-43}) and the $2$-nd and the $3$-rd  order terms in (\ref{19-3-44}) because for  $\mu \ge \epsilon (h|\log h|)^{-\frac{1}{2}}$ only these terms could be  larger than $C\mu ^{-1}h$.

Thus we arrive to

\begin{proposition}\label{prop-19-3-6} 
Let conditions \textup{(\ref{19-1-4})},  \textup{(\ref{19-3-5})} and \textup{(\ref{19-3-6})} be fulfilled, $\bar{x}$ be a fixed point. Then there exists  a bounded operator $\cT$ such that

\medskip\noindent
(i) For operators $Q=q^\w$  with the symbol $q$ supported in 
$B(0,R)\subset \bR^{2d}$ and  $Q'=q_1^{\prime\,\w}$  with the symbol $q'_1$ supported in $B((\bar{x},0),1)\subset \bR^{2d}$ and equal $1$ in
$B((\bar{x},0),1)\subset \bR^{2d}$ equalities \textup{(\ref{19-3-31})} and \textup{(\ref{19-3-32})} hold modulo negligible operators and equalities \textup{(\ref{19-3-33})} and \textup{(\ref{19-3-34})}  hold modulo operators with the norms not exceeding
$C\bigl(\mu^{-1}h + \mu^{-l}|\log h|^{-\sigma}\bigr)$ with
\begin{multline}
\cA = \cA_0 +
\mu ^2 \Re 
\sum _{\alpha,\beta :\,3\le |\alpha|+|\beta| \le \lceil (l,\sigma)\rceil +1}
b''_{\alpha\beta} (x'',\mu ^{-1}hD'')\,
\bigl(\zeta^\alpha \zeta ^{\dag\,\beta}\bigr)^\w +\\[2pt]
\Re \sum _{\alpha,\beta :\, 1\le |\alpha|+|\beta|\le  \lceil ( l,\sigma)\rceil -1}
b'''_{\alpha\beta} (x'',\mu ^{-1}hD'')\,
\bigl(\zeta^\alpha \zeta ^{\dag\,\beta}\bigr)^\w
\label{19-3-50}
\end{multline}
with $\cA_0$ defined by \textup{(\ref{19-3-42})} (and equal to \textup{(\ref{19-3-40})}) and with the symbols
$b''_{\alpha\beta},b'''_{\alpha\beta} \in \sF^{l -|\alpha|-|\beta|,\sigma}$ and
$b''_{\alpha\beta}=b'''_{\alpha\beta}=0$ unless
\begin{equation}
|\sum _j (\alpha_j -\beta_j) f_j|\le \epsilon
\label{19-3-51}
\end{equation}
with arbitrarily small constant $\epsilon >0$;

\medskip\noindent
(ii) Further,
\begin{equation}
\bar{\psi} =
\sum _{\alpha,\beta :\, 1\le |\alpha|+|\beta|\le  \lceil ( l,\sigma)\rceil -1}
\bar{\psi}_{\alpha\beta} (x'',\xi'')\, \zeta^\alpha \zeta ^{\dag\,\beta},
\label{19-3-52}
\end{equation}
where $ \bar{\psi}_{\alpha\beta} \in \sF^{l-|\alpha|-|\beta|,\sigma}$ and
\begin{align}
&\bar{\psi}_{00}= \psi \circ\Psi_0, \label{19-3-53}\\
&\bar{\psi}_{\alpha\beta} = \sum_{\gamma:\, 1\le |\gamma|\le |\alpha|+|\beta|}
c_{\alpha\beta\gamma}(\partial_x^\gamma \psi)\circ \Psi_0 \quad 
\text{as\ \ }|\alpha|+|\beta|\ge 1;\label{19-3-54}
\end{align}
(iii) Statement (ii) of proposition \ref{prop-19-3-4} remains true.
\end{proposition}

Recall that only cubic terms are obstacles in the proof that
$\tilde{a}^0_{\fm}$ evolve slowly i.e. with the speed $O(\mu^{-1})$; we  prove this statement  only for $\tilde{a}^0_{\fn}$.

\section{Strong and superstrong magnetic field case}
\label{sect-19-3-5}

\subsection{Strong  magnetic field case}
\label{sect-19-3-5-1}
The same construction works in the strong magnetic field case 
\begin{equation}
\mu_2^*\Def\epsilon (|\log h|)^{-1}\le \mu \le \mu_3^*\Def \epsilon h^{-1}
\label{19-3-55}
\end{equation}
but in this case one needs to take $\varepsilon$ according to (\ref{19-3-7}) i.e. $\varepsilon= C(\mu^{-1}h|\log h|)^{\frac{1}{2}}$ or larger as it is larger now than $C\mu^{-1}$. Then restriction to $p_j$ is now
\begin{equation}
|p_j|\le C\varepsilon .
\label{19-3-56}
\end{equation}
However, in the operator rather than microlocal sense we have still
\begin{equation}
\|p^\w_j\|\le C\mu ^{-1}
\label{19-3-57}
\end{equation}
on the energy levels below $c_0$ and we should estimate perturbations based on this estimate rather than (\ref{19-3-56}).

Note, that now we need to take in account only the $3$-rd and the $4$-th order terms in the first sum in the right-hand expression of (\ref{19-3-50}) and only the $1$-st and the $2$-nd order terms in the second sum in the right-hand expression of (\ref{19-3-50}). Similarly   in (\ref{19-3-52}) we need to sum as
$|\alpha|+|\beta|\le 2$.

\subsection{Very strong and superstrong magnetic field case}
\label{sect-19-3-5-2}

The same construction also works in the very strong and superstrong magnetic field cases
\begin{equation}
\mu \ge \mu_3^{*} \Def \epsilon h^{-1}
\label{19-3-58}
\end{equation}
we should take again $\varepsilon$ by (\ref{19-3-7}) and inequality (\ref{19-3-56}) again holds in the microlocal sense but (\ref{19-3-57}) is replaced by
\begin{equation}
\|p^\w_j\|\le \varsigma \Def C(\mu ^{-1}h)^{\frac{1}{2}}
\label{19-3-59}
\end{equation}

Really, then
\begin{equation}
\cA \equiv \mu^2 \sum_{ \fm \in \fM }\ \sum_{j,k\in {\fm}} a_{jk}(x'',\mu^{-1}hD'')\, Z_j^*Z_k + b'_0 (x'',\mu^{-1}hD'')
\label{19-3-60}
\end{equation}
with symbol $b'_0 \in \sF^{l,\sigma}$.

Then the principal part of operator $\bar{A}_0$ is of magnitude $C\mu h$
while norms of perturbations
$\bar{A}_1$, $\bar{A}_2$, $\bar{A}_3$, $\bar{A}_4$
do not exceed $C\mu h \varsigma$, $C\mu h \varsigma ^2 = Ch^2$,
$C\mu h \varsigma ^3$, $C\mu h \varsigma ^4 = C\mu ^{-1}h^3$ respectively and only the last one could be ignored for really large $\mu$. This is the reason beyond correction $h^{-d}\cN^\MW_{2\corr}$ in theorem \ref{thm-19-6-26} below. It is an also one of the reasons beyond the last term in the right-hand expression of estimate (\ref{19-6-65}) and the corresponding extra-smoothness requirement there to get the best possible estimate; another reason is the mollification error.

Anyway, in both strong and superstrong cases we arrive to

\begin{proposition}\label{prop-19-3-7} 
All statements of proposition \ref{prop-19-3-6} remains true in both strong and superstrong magnetic field cases with \textup{(\ref{19-3-34})} fulfilled modulo operator $\cR$ such that
\begin{multline}
\|\cR v\| \le \\ C\Bigl(\varepsilon ^l |\log \mu |^{-\sigma}+
\mu h \varepsilon ^{\bar{l}} |\log \mu |^{-\bar{\sigma}} \Bigr)\,  \|v\| +
\varepsilon^{-\bar{l}}|\log \mu|^{-\bar{\sigma}}\, \|\cA  v\| \qquad \forall v
\label{19-3-61}
\end{multline}
with $\varepsilon$ given by \textup{(\ref{19-3-7})}.
\end{proposition}

\chapter{Intermediate magnetic field}
\label{sect-19-4}

Let us consider intermediate magnetic field case (\ref{19-3-5}). While generally we assume that microhyperbolicity condition is fulfilled, in Subsection~\ref{sect-19-4-7} we consider the case of constant $g^{jk}$, $F_{jk}$ when $V$ has critical points.

\section{Mid-range propagation}
\label{sect-19-4-1}

By mid-range propagation we assume propagation for operator $\cA$ and  thus for
operator $\cT^* \cA \cT$ as  $T\in [T'_*,T^{*\,\prime}]$ with
\begin{gather}
T^{*\,\prime} = C(\mu h |\log h|)^{\frac{1}{2}},\label{19-4-1}\\[2pt]
T'_*= C\varepsilon^{-1} h |\log h| \le C\mu h|\log h|\label{19-4-2};
\end{gather}
we are interested almost exclusively in the extending proposition \ref{prop-19-2-12}. Here $T^{*\,\prime}$ is what used to be $T_*$\,\footnote{\label{foot-19-20} See (\ref{19-2-53}) and (\ref{19-2-65}) as, for a intermediate magnetic field we have already proven that trace is negligible as $T'_*\le T\le T^*$.} and $T'_*$ is the same as $T_*$ as if $\cA$ was $\mu^{-1}h$-pseudo-differential operator $\cA(x'',\mu^{-1}hD'')$ satisfying microhyperbolicity condition. 

\begin{remark}\label{rem-19-4-1}
Recall that then  $\varepsilon \gtrsim \mu^{-1}$ but we also assume that 
$\varepsilon \le h^\delta$ to avoid some complications as otherwise logarithmic uncertainty principle should be replaced by a microlocal uncertainty principle and below $T'_*= C h^{1-\delta}$.
\end{remark}

However $T'_*$ is still pretty large unless we increase $\varepsilon$ and thus smoothness assumptions. Instead  in the next Section we consider (smoothness-dependent) rate of decay of the left-hand expression of (\ref{19-2-63}) on interval $[T''_*,T^{*\,\prime\prime}]$ with $T^{*\,\prime\prime}=T'_*$ and $T''_*= Ch $.

Consider some point $(\bar{x},\bar{\xi})$. During time $T^{*\prime}$ (as now we assume that $(\bar{\ell},\bar{\sigma})\succeq (2,0)$)  $(x'',\xi'')$ will stay in $C(\mu ^{-1}h|\log h|)^{\frac{1}{2}}$-vicinity of $(\bar{x}'',\bar{\xi}'')$ in both classical and microlocal senses (we can follow proposition \ref{prop-19-2-9} for the proof in the microlocal sense).

According to Subsection~\ref{sect-19-3-4}
\begin{gather}
A_\cT \Def \cT^* A \cT= \cA +\cR,
\label{19-4-3}\\
\shortintertext{with}
\cA = \sum _{0\le m \le \lceil (l,\sigma)\rceil -1}\cA _m,\label{19-4-4}\\[3pt]
\cR =\rho (x,\mu^{-1}hD),
\qquad \rho \in \mu^{-l}|\log \mu|^{-\sigma} \sF^{0,0},
\label{19-4-5}
\end{gather}
$\cA_0$ given by (\ref{19-3-42}) (where one can rewrite quadratic part as (\ref{19-3-40}) in the general case and as (\ref{19-3-41}) as $\# \fm =1$ for all ${\fm}$) and
\begin{multline}
\cA_m =
\mu ^2 \Re \sum _{\alpha,\beta :\,  3\le |\alpha|+|\beta| =m+2 }
b''_{\alpha\beta} (x'',\mu ^{-1}hD'')\,
\bigl(\zeta^\alpha \zeta ^{\dag\,\beta}\bigr)^\w +\\
\Re \sum _{\alpha,\beta : 1\le |\alpha|+|\beta| =m }
b'''_{\alpha\beta} (x'',\mu ^{-1}hD'')\,
\bigl(\zeta^\alpha \zeta ^{\dag\,\beta}\bigr)^\w
\label{19-4-6}
\end{multline}
(see (\ref{19-3-50})) and $\zeta_j =\xi_j+ ix_j$, $\zeta^\dag _j =\xi_j- ix_j$,
$j=1,\ldots,r$.

\begin{remark}\label{rem-19-4-2}
Surely  we need to remember that $\cA$ and $A_\cT\Def \cT^* A \cT$ are close but not equal. In Chapter~\ref{book_new-sect-18} of \cite{futurebook} we considered $\cA$ and took $\cT \cA \cT^*$ as approximation\footnote{\label{foot-19-21} Recall that actually we consider $A_\varepsilon$
instead of $A$, so we already have an approximation.} but now we have the global
reduction only in the special case of condition (\ref{19-1-25}) while until the next Section we need only a weaker assumption 
$V_j\in \sC^{\bar{l}+1,\bar{\sigma}}$.

The following statements hold both for propagators
$\mathbf{U} '=\mathbf{U}_\cT \Def \cT^*e^{-h^{-1}At}\cT$ and its Schwartz kernel $U'\Def U_\cT$ and for $\mathbf{U} '= e^{ih^{-1}\cA t}$ and its Schwartz kernel $U'$.
\end{remark}

\begin{proposition}\label{prop-19-4-3} 
Let $U'$ be defined in remark \ref{rem-19-4-2}. Let $\phi _1$ be supported in $B(0,1)$, $\phi _2=1$ in $B(0,2)$, $\chi $ be supported in $[-1,1]$. Let $Q_k=\phi _{k, M\mu^{-1} T} (x''-\bar{x}'', \xi ''-\bar{\xi}'')^\w$.  

\medskip\noindent
(i) Then for $T^{*\, \prime}\le T\le C $ and large enough constant $M$
\begin{equation}
|F_{t\to h^{-1}\tau }
\bar{\chi} _T(t) \bigl(1- Q_{2x}\bigr)U_{\alpha\beta}(x,y,t)\,^t\!Q_{1,y} |
\le C\mu^{-s}\qquad \forall \tau \le c.
\label{19-4-7}
\end{equation}

\medskip\noindent
(ii) Further, the same estimate holds for $T'_*\le T \le C$ if we replace $\phi_*(x'',\xi'')$ by $\phi_*(x''_I,\mu^{-1}hD''_{II})$ with an arbitrary partition $x''=(x''_I;x''_{II})$.
\end{proposition}

\begin{proposition}\label{prop-19-4-4} 
Let $U'$ be defined in remark \ref{rem-19-4-2}.
Let  $m\Def \#\fM \ge 2$ and let   $\phi_j$ be functions as before. Then for $C(\mu h|\log h|)^{\frac{1}{2}}\le T \le c_0$
\begin{multline}
|F_{t\to h^{-1}\tau } \chi _T(t)
\bigl(1-\phi_{2,T} (\mu ^2 a^0_{{\fm}_1}-\tau_1,\ldots,
\mu^2a^0_{{\fm}_m}-\tau_m) \bigr)^\w \times \\[4pt]
\times U'(x,y,t)
\,^t\!\bigl(\phi _{1,T} (\mu ^2 a^0_{{\fm}_1}-\tau_1,\ldots,
\mu^2a^0_{{\fm}_m}-\tau_m)  \bigr)^\w_y|
\le C\mu^{-s}\qquad \forall \tau :\ |\tau |\le \epsilon_1.
\label{19-4-8}
\end{multline}
\end{proposition}

\begin{proof}[Proof propositions~\ref{prop-19-4-3}--\ref{prop-19-4-4}] Both propositions \ref{prop-19-4-3}--\ref{prop-19-4-4} are proven by the same scheme as in Chapter~\ref{book_new-sect-18} of \cite{futurebook}:

\medskip
To prove proposition \ref{prop-19-4-3}(i), (ii) one can use functions
\begin{gather}
\upchi \Bigl(\bigl( \frac {\mu^2 (|x''-\bar{x}''|^2+|\xi''-\bar{\xi}''|^2)} {T^2} +
\epsilon^2\bigr)^{\frac{1}{2}} - C\varsigma  \frac{t}{T}\Bigr)\label{19-4-9}\\
\shortintertext{and}
\upchi \Bigl(\bigl( \frac {\mu^2 (|x''_I-\bar{x}_I''|^2+|\xi_{II}''-\bar{\xi}_{II} ''|^2)} {T^2} +
\epsilon^2\bigr)^{\frac{1}{2}} - C\varsigma  \frac{t}{T}\Bigr)
\label{19-4-10}
\end{gather}
respectively with $\upchi$ function of the same type as used in theorem~\ref{book_new-thm-2-3-1} of \cite{futurebook} and  $\varsigma =\pm 1$ depending on time direction and arbitrarily small constant $\epsilon >0$.

\medskip
To prove proposition \ref{prop-19-4-4} one can use function
\begin{equation}
\upchi \Bigl(\bigl(
\frac{\mu^4 \sum_j |a_{\fm _j}-\tau_{\fm _j}|^2 }{T^2} + \epsilon^2\bigr)^{\frac{1}{2}} -
C\varsigma \frac{t}{T}\Bigr).
\label{19-4-11}
\end{equation}
Because $T \le \epsilon_0$ and we picked up velocity with respect to 
$a_{\fm _j}$, the $3$-rd order resonances do not pose any problem here.
\end{proof}

Even if at this stage the decomposition into Hermitian functions is not very useful, we can study propagation using only functions of $(x'',\xi'',t)$.

\begin{proposition}\label{prop-19-4-5} 
Let $U'$ be defined in remark \ref{rem-19-4-2}. Let microhyperbolicity condition (see definition~\ref{def-19-2-4}) be fulfilled.
Then for $T\in [T'_*,c_0]$, $\tau \in [-\epsilon',\epsilon']$
\begin{equation}
|F_{t \to h^{-1}\tau} \chi_T(t) \Gamma  Q U | \le C\mu^{-s}
\label{19-4-12}
\end{equation}
where $Q=Q(x,\mu^{-1}hD)$.
\end{proposition}

\begin{proof}
Without any loss of the generality one can assume that 
$\ell (\bar{x}'', \bar{\xi}'' ;\bar{\tau}_1,\ldots,\bar{\tau}_m)$ is equal to $(\ell ',0)$. Otherwise one can achieve it by a linear symplectic transformation in $(x'',\xi'')$ and corresponding FIO (which will be a metaplectic operator in this case). Then  using
\begin{equation}
\upchi
\Bigl(\mu T^{-1} \langle \ell ', x''-y''\rangle \pm \epsilon tT^{-1}\Bigr)
\label{19-4-13}
\end{equation}
which is an admissible pseudo-differential operator-symbol one can prove in the standard way that
\begin{equation}
|F_{t \to h^{-1}\tau} \chi_T(t) \Gamma'' Q u | \le C\mu^{-s}
\label{19-4-14}
\end{equation}
where  here and below $\Gamma''$ is a partial trace (with respect to $x''$ only).
\end{proof}

Propositions \ref{prop-19-2-13} and \ref{prop-19-4-5} imply immediately

\begin{corollary}\label{cor-19-4-6} 
(i) Under microhyperbolicity condition  estimate \textup{(\ref{19-2-63})} holds for $T\in [T'_*,T^*]$ with $T^*=\epsilon$; 

\medskip\noindent
(ii) Under $\fN$-microhyperbolicity condition one can take $T^*$ defined by \textup{(\ref{19-2-46})};

\medskip\noindent
(iii) Also estimate \textup{(\ref{19-2-76})} holds with $T_*$ replaced by the lesser value $T'_*$.
\end{corollary}

\begin{corollary}\label{cor-19-4-7}
Estimates \textup{(\ref{19-2-86})}--\textup{(\ref{19-2-88})} hold under corresponding assumptions of theorem~\ref{thm-19-2-16}  in the case of $\mu \le h^{\delta-1}$ (albeit they are not optimal).
\end{corollary}

\begin{proof} 
Picking $\varepsilon= C\mu h|\log h|$ we get approximation error
$Ch^{-d}(\mu h)^l |\log h|^{l-\sigma}$ and $T'_*=\epsilon \mu^{-1}$. Then we again have estimates (\ref{19-2-77}), (\ref{19-2-78}) and need just calculate expression (\ref{19-2-79}) which is equivalent (modulo negligible) to (\ref{19-2-80}) where without spoiling remainder estimate we can take only terms with $m=0$ and $n\le n(\delta)$ as $\mu \le h^{-1+\delta}$ with arbitrarily small exponent $\delta>0$. The rest repeats arguments of Subsection~\ref{sect-19-2-4}.
\end{proof}

\section{Short-range theory: framework}
\label{sect-19-4-2}
Still, $T_*=C\varepsilon^{-1} h|\log h|$ is too large. Really, the Tauberian method implies that the remainder does not exceed
\begin{equation}
\frac{1}{T} \sup_{|\tau|\le \epsilon }|F_{t\to h^{-1}\tau} \bigl(\bar{\chi}_T(t) \Gamma U\psi \bigr)|
\label{19-4-15}
\end{equation}
while the principal part is
\begin{equation}
\frac{1}{h} \int_{-\infty}^0 \Bigl(F_{t\to h^{-1}\tau}
\bigl(\bar{\chi}_T(t) \Gamma U\psi \bigr)\Bigr)\,d\tau.
\label{19-4-16}
\end{equation}
We will reduce (\ref{19-4-16}) to a more explicit form in Subsection~\ref{sect-19-4-5} (where this Subsection approach will be crucial as well) but now let us consider expression (\ref{19-4-15}). To decrease it one should increase $T$ but  generally supremum also grows proportionally $T$ and there is no improvement; the following estimate
\begin{equation}
|F_{t\to h^{-1}\tau} \bigl(\bar{\chi}_T(t) \Gamma U\psi \bigr)|\le
CT \Bigl(h^{-2r}+\mu^rh^{-r}\Bigr)\qquad \forall \tau\le c
\label{19-4-17}
\end{equation}
rather easily follows from
\begin{equation}
|\phi(hD_t) \bigl(\bar{\chi}_T(t) \Gamma U\psi \bigr)|\le
C \Bigl(h^{-2r}+\mu^rh^{-r}\Bigr).
\label{19-4-18}
\end{equation}
However results of the previous Subsection shows that while taking $T=T^*$
one can take $\bar{\chi}_{T_*}(t)$ and then  (\ref{19-4-15}) becomes
$CT_*T^{*\,-1} \bigl(h^{-2r}+\mu^rh^{-r}\bigr)$. It equals to our dream remainder estimate as $T_*\Def C_0h$. Recall that 
$T_*\Def C\varepsilon^{-1}h|\log h|$.

Note that even in the smooth case $T_*$ is larger\footnote{\label{foot-19-22} Since ``no mollifications needed'' has its own pitfalls, we would need actually request $T_*=h^{1-\delta}$ with arbitrarily small exponent $\delta>0$ rather than $T_*=Ch|\log h|$.}. 

We cannot just take $T_*=C_0h $ because  
$F_{t\to h^{-1}\tau} \bigl(\chi_T(t) \Gamma U\psi \bigr)$ is not negligible for $T\in [C_0h,T_*]$; instead our goal is to prove under microhyperbolicity condition that for $T\in (C_0h,T_*)$ 
\begin{multline}
|F_{t\to h^{-1}\tau} \bigl(\chi_T(t) \Gamma U\psi \bigr)|\le
CT \Bigl(h^{-2r}+\mu^rh^{-r}\Bigr) \times
({\frac h T})^l|\log ({\frac h T})|^{-\sigma}\\[3pt]
\forall \tau:|\tau |\le \epsilon'.
\label{19-4-19}
\end{multline}
As $(l,\sigma)\succ (1,1)$ it would imply that left-hand expression of (\ref{19-4-17}) with $T=T_*$ does not exceeds its right-hand expression with $T=C_0h$ i.e. that in the estimates effectively $T_*=C_0h$. 

To achieve this goal we apply $\eta$-mollification to $A_\cT$; this will lead
to approximation error in  
$F_{t\to h^{-1}\tau} \bigl(\bar{\chi}_T(t) \Gamma U\psi \bigr)$ which we are going to estimate and to (kind of) negligibility estimate of it after mollification. Minimizing the sum by $\eta$ (which will depend on $T$) we will arrive to (\ref{19-4-19}) for non-mollified operator (surely, the original mollification will be still here).

Alternatively we can try the partition approach like in the previous Chapter.

\section{Pilot-Model: some classes of pseudo-differential operators}
\label{sect-19-4-3}

We start from related simple results (we will need them anyway):

\begin{proposition}\label{prop-19-4-8} 
Let $A=A(x,hD)$ be a (matrix) self-adjoint $h$-pseudo-differential operator in $\bR^r$ with symbol  $a\in \sF^{l,\sigma}$ and satisfying microhyperbolicity condition on level $\tau=0$:
\begin{equation}
\langle (\ell a)v,v\rangle \ge \epsilon \|v\|^2 -C\|(a-\tau)v\|^2\qquad
\forall v\quad \forall (x,\xi).
\label{19-4-20}
\end{equation}
Let $U$ be a Schwartz kernel of $e^{-ih^{-1}A}$. Finally, let 
$\psi\in \sC_0^\infty$ be supported in the small vicinity of $0$.

\medskip\noindent
(i) Let $l>1$. Then for $T\in [Ch^{1-\delta}, T^*]$ with arbitrarily small exponent $\delta >0$ and small enough constant $T^*$
\begin{align}
& |\phi (hD_t) \chi_T(t) ( \Gamma U\psi_y )|\le
C h^{-r}\bigl( \frac{h} {T}\bigr)^{l-1}
\bigl(1+  \frac{T\varepsilon}{h} \bigr)^{-s}|\log \frac{h} {T}|^{-\sigma}, \label{19-4-21}\\[3pt]
&|F_{t\to h^{-1}\tau} \chi_T(t) ( \Gamma U\psi_y)|\le
C h^{1-r}\bigl(\frac{h}{T}\bigr)^{l-1}
\bigl(1+\frac{T\varepsilon}{h} \bigr)^{-s}
|\log \frac{h}{T}|^{-\sigma}\label{19-4-22}
\end{align}
as $\tau:|\tau|\le\epsilon$ with arbitrarily large exponent $s$.

\medskip\noindent
(ii) Let $l=1,\sigma \ge 2$. Then for $T\in [Ch^{1-\delta}, T^*]$ with arbitrarily small exponent $\delta >0$ and small enough constant $T^*$
\begin{multline}
|\phi (hD_t) \chi_T(t) ( \Gamma U\psi_y )|\le\\
C h^{-r}|\log \frac{h} {T}|^{-\sigma}
\bigl(1+\frac {T\varepsilon}{h}\bigr)^{-s}+
C h^{-r}|\log \frac{h} {T}|^{-s},
\label{19-4-23}
\end{multline}\vglue-15pt
\begin{multline}
|F_{t\to h^{-1}\tau} \chi_T(t) ( \Gamma U\psi_y)|\le\\
C h^{1-r}|\log \frac{h} {T}|^{-\sigma}
\bigl(1+ \frac {T\varepsilon}{h}\bigr)^{-s}+C h^{1-r}|\log \frac{h} {T}|^{-s}
\label{19-4-24}
\end{multline}
as $\tau:|\tau|\le\epsilon$ with arbitrarily large exponent $s$.
\end{proposition}

\begin{proof}[Pilot-model] Before an actual proof, let us start with a simple example (actually arising as $d=2$\,\footnote {\label{foot-19-23} In this case one can apply successive approximations on this stage and to arrive to oscillatory integral expression.}):

In the case of scalar operator  $\Gamma (U\psi )$ is something like
$h^{-d}I_1$ with 
\begin{gather}
I_1=\int e^{ih^{-1}\lambda (z)t}q(z)\,dz\label{19-4-25}\\
\intertext{and $F_{t\to h^{-1}\tau} \chi_T(t) \Gamma (U\psi )$ is $h^{1-d}I_2$ with}
I_2= h^{-1}\iint \chi (t/T)e^{ih^{-1}\lambda (z)t}q(z)\,dz\, dt.
\label{19-4-26}
\end{gather}
Let us consider the case $T\le \varepsilon^{-1}h$, $l>1$.

Replacing $\lambda $ by its  $\eta$-mollification we make an error
$O\bigl(\vartheta (\eta)\bigr)$ in $\lambda$ which leads to error
$O\bigl(Th^{-1}\vartheta (\eta)\bigr)$ in $I_1$. In the same time a (multiple) integration by parts with respect to $t$ shows that the error in $I_2$
is does not exceed
\begin{equation}
Ch^{-2}T^2 \int (Th^{-1}|\lambda -\tau |+1)^{-2}  \,dz \times
\vartheta (\eta) \le C h^{-1}T  \vartheta (\eta).
\label{19-4-27}
\end{equation}
Consider now $I_1$ and $I_2$ with mollified $\lambda$ assuming that $|\partial_{z_1} \lambda |\ge \epsilon$. A multiple integration by parts with respect to $z_1$   will transform modified $I_1$ into
$(h/T)^k\int e^{ih^{-1}\lambda (z)t}q_k(z)\,dz$ where
$q_k=O\bigl(\eta ^{-k-1}\vartheta (\eta)\bigr)$  as $k$ is large enough and therefore modified expression $I_1$ is 
$O\bigl((h/T)^k \eta ^{-k-1} \vartheta (\eta)\bigr)$.

We can treat expression $I_2$ in the similar way, but we add an integration by parts with respect to $t$ as well and use of nondegeneracy condition as above. Then modified expression $I_2$ has the same upper estimate.

Therefore both original expressions $I_1$ and $I_2$ do not not exceed
\begin{gather*}
CTh^{-1}\vartheta (\eta) +
\bigl(\frac{h}{T}\bigr)^k \eta ^{-k-1} \vartheta (\eta).\\
\intertext{Picking (near) optimal $\eta =h/T$ we get estimates}
|I_k|\le C\bigl(\frac{h}{T}\bigr)^{-1}\vartheta\bigl(\frac{h}{T}\bigr).
\end{gather*}
It proves (\ref{19-4-21}), (\ref{19-4-22})  in this special case.
\end{proof}

\begin{proof}[Proof proposition \ref{prop-19-4-8}]
(a) Note first that due to the standard propagation results for
$T\ge Ch\varepsilon ^{-1}|\log h|$ left-hand expressions in (\ref{19-4-21}), (\ref{19-4-22}) are negligible. Thus we need to consider only
$T\le T'_* \Def Ch\varepsilon ^{-1}|\log h|$.  We can consider $\chi $ supported in $[\frac{1}{2},1]$ and until the end of the proof  $t,T_k$ are positive.

Consider $\eta $-mollification $A_\eta$ of operator $A$ with respect to $(x,\xi)$\,\footnote{\label{foot-19-24} It means that we mollify its Weyl symbol; recall that original symbol is regular in $\varepsilon$-scale.}:
\begin{equation}
1\gg \eta \ge \varepsilon.
\label{19-4-28}
\end{equation}
Then the  error in the operator (the difference between $A_\eta$ and $A$)  will not exceed $C\vartheta(\eta)$ in the sense that
\begin{gather}
\| (A_\eta -A) \|\le C\vartheta(\eta)
\label{19-4-29}\\
\shortintertext{and identity}
e^{ih^{-1}tA}=e^{ih^{-1}t A _\eta}+ih^{-1}\int_0^t
e^{ih^{-1}(t-t')A }(A  - A _\eta)e^{ih^{-1}t'A _\eta}
\, dt'
\label{19-4-30}
\end{gather}
implies that
\begin{align}
& \|\bigl(e^{ih^{-1}tA }-e^{ih^{-1}t A _\eta}\bigr) \|\le CTh^{-1}\vartheta (\eta)\qquad  \text{as\ \ }|t|\le T,\label{19-4-31}\\[3pt]
&| \chi_T (t)\Gamma  \bigl((U _\eta -U)\psi_y \bigr)| \le Ch^{-d-1}T\vartheta (\eta)
\label{19-4-32}
\end{align}
for the standard cut-off $\psi $ in $B(\bigl(\bar{x},\bar{\xi}), 1\bigr)\subset \bR^{2r}_{(x,\xi)}$; this is exactly what we got in the pilot-model.

On the other hand,  microhyperbolicity condition (\ref{19-4-20}) and theorem~\ref{book_new-thm-2-3-1} of \cite{futurebook} imply that
$\phi (hD_t)\bigl(\chi_T (t)\Gamma   U _\eta \psi_y \bigr)$ is negligible provided
\begin{equation}
h^\delta \ge \eta \ge ChT^{-1}|\log h|
\label{19-4-33}
\end{equation}
(and $\supp \phi\subset [-\epsilon, \epsilon])$  and picking minimal
$\eta =ChT^{-1}|\log h|$ we arrive to inequality
\begin{equation}
|\phi (hD_t)\Gamma  \chi_T (t) U  \psi| \le
C h^{-d}\bigl( \frac{h}{T}\bigr)^{l-1}\cdot |\log h|^{l-\sigma}
\label{19-4-34}
\end{equation}
due to $T\ge h^{1-\delta}$.
This inequality  is almost as good as (\ref{19-4-21}): we have an extra factor $|\log h|^l$ in the right-hand expression.

\medskip\noindent
(b) To get gid  of this $|\log h|^l$ factor we need a more delicate analysis. First, let us iterate (\ref{19-4-30}):
\begin{multline}
e^{ih^{-1}tA }= \sum_{0\le p\le m-1}\frac{1} {p!} i^ph^{-p}\int_{\Delta_p} e^{ih^{-1}(t-t_1-\ldots-t_p)A_\eta}(A-A_\eta)\times \\[3pt]
e^{ih^{-1}t_1A_\eta}(A-A_\eta)e^{ih^{-1}t_2A_\eta}
\cdots e^{ih^{-1}t_{p-1}A_\eta}(A -A_\eta) e^{ih^{-1}t_p A_\eta } \,dt_1\cdots\,dt_p +R_m
\label{19-4-35}
\end{multline}
with the remainder
\begin{multline}
R_m=\frac{1}{m!} i^m h^{-m}\int_{\Delta_m} e^{ih^{-1}(t-t_1-\ldots-t_m)A_\eta}(A-A_\eta)\times \\[3pt]
e^{ih^{-1}t_1A_\eta}(A-A_\eta)e^{ih^{-1}t_2 A_\eta}
\cdots e^{ih^{-1}t_{m-1}A_\eta} (A -A_\eta) e^{ih^{-1}t_m A } \,dt_1\cdots\,dt_m.
\label{19-4-36}
\end{multline}
Then the  trace norm of the remainder does not exceed
$C\bigl(Th^{-1}\vartheta (\eta )\bigr)^mh^{-d}$ and it is less than the right-hand expression of (\ref{19-4-21}) for large enough $m$ (because we gain factor $(h/T)^{l-1}|\log (h/T)|^{l-\sigma}$ on each step).

Now as we take $p=1$ we need to consider
\begin{equation}
e^{ith^{-1}A_\eta} \bigl( ih^{-1}t (A -A_\eta )\bigr)\psi .
\label{19-4-37}
\end{equation}
Without any loss of the generality one can assume that the direction of microhyperbolicity is $\ell =\partial_{\xi_1}$. Let us apply
$\epsilon  T$-admissible partition  with respect to
$x_1$: $\psi '= \psi '\sum _\nu  \varphi_k(x_1)$. Here and below
$\varphi $ (with different indices) are $\epsilon T$-admissible functions.

Due to (\ref{19-4-33}) the logarithmic uncertainty principle
$T \times \eta \ge C'h|\log h|$ holds and therefore due to the microhyperbolicity
$\varphi _\nu  e^{ith^{-1}A}\equiv
\varphi _\nu  e^{ith^{-1}A } \tilde{\varphi}_\nu$ (modulo operator, negligible after cut-off $\phi (hD_t)$) as $\frac{1}{2}T\le t\le T$ where the distance between $\supp \tilde{\varphi}_\nu $  and  $\supp \varphi_\nu $ is exactly of magnitude $\mu^{-1}T$. Then to prove (\ref{19-4-21}) it is  sufficient to estimate properly a trace norm of
$ \sum  _\nu {\tilde\varphi}_\nu  B \psi  \varphi'_\nu \times
(h/T)^{-1}$. To do this it is sufficient to prove that

\begin{claim}\label{19-4-38}
The operator norm of each ``sandwich''
$\tilde{\varphi}_\nu  B \psi ' \varphi_\nu $ does not exceed
\begin{equation*}
M  \Def \bigl(\frac{h}{T}\bigr)^l\bigl(1+\frac{h}{T\varepsilon}\bigr)^{-s}
|\log \bigl(\frac{h}{T}\bigr) |^{-\sigma}.
\end{equation*}
\end{claim}\vglue-10pt

Really, then the trace norm of each ``sandwich'' will be
$CT  h^{-r}  \times M_k$  where $T$ is $x_1$-size of $\supp \varphi_\nu$. Then
after summation with respect to $\nu$ factor $T$ disappears and we estimate the trace norm of (\ref{19-4-37}) (after $\phi(hD_t)$ cut-off) by 
$h^{-r} M  \bigl(h/T\bigr)^{-1}$, which implies (\ref{19-4-21}) and (\ref{19-4-23}).

For $p\ge 2$ we replace the integral over $\Delta_p$ by the sum of the integrals over $\Delta_p\cap \{|t_j-\bar{t}_j|\le \epsilon' T,\ j=1,\ldots,p\}$ with very small constants $\epsilon'$. 

We also place partitions elements around each copy of the operator exponent, so $e^{ih^{-1}t_jA_\eta}$ is now replaced by
$\varphi_{\nu_{2j}}e^{ih^{-1}t_jA_\eta}\varphi_{\nu_{2j+1}}$, $t_0=t-t_1-\ldots-t_p$.

Then \underline{either} supports of  $\varphi_{\nu_{2j+1}}$ and $\phi_{\nu_{2j+2}}$ are disjoint by a distance at least $\epsilon\mu^{-1}T$ for some $j$ and we have a sandwich  $\varphi_{\nu_{2j+1}}B\varphi_{\nu_{2j+2}}$ of the type described in (\ref{19-4-38}), \underline{or} supports of  $\varphi_{\nu_0}$ and $\varphi_{\nu_{2p+1}}$ are disjoint by a distance at least $\epsilon\mu^{-1}T$ (and then after taking trace we get negligible term) \underline{or} $\varphi_{\nu_{2j}}e^{ih^{-1}t_jA_\eta}\varphi_{\nu_{2j+1}}$ is a negligible operator.

Since any sandwich has a norm not exceeding $C\eta^l |\log \eta|^{-\sigma}$ and
is accompanied by a factor $T/h$, it brings as a result a factor
$(T/h)^{l-1}|\log (T/h)|^{l-\sigma}$.

To estimate the norm of a sandwich coming from (\ref{19-4-36})-term let us consider the sequence of mollifications $A_\rho$ with 
$\rho=\rho_n = 2^n \rho_0$, $\rho_0=\max\bigl( h/T,\varepsilon\bigr)$ and
$n=1,2, \ldots , \bar{n}$ such that $\rho_{\bar{n}+1} =\eta$. 

Then operator $B = \sum_{0\le n \le \bar{n}} B_n$ with 
$B_n= A_{\rho_n}-A_{\rho_{n+1}}$ has a $\vartheta (\rho) \sF^{0,0}_\rho$ symbol
with $\rho\Def \rho_n$. Then  the norm of the sandwich
$\tilde{\varphi}_\nu B_n  \varphi_\nu $ does not exceed
$C\vartheta (\rho) \times (h/T\rho)^s$ with arbitrarily large $s$ and summation with respect to $n$ does not exceed the same expression with
$\rho =\rho_0= hT^{-1}$; so we get
$\vartheta (h/T) \bigl(1+T\varepsilon/h \bigr)^{-s}$ which is exactly $M$ as defined in (\ref{19-4-38}).

Thus (\ref{19-4-21}) and (\ref{19-4-23}) are proven.

\medskip\noindent
(b) To prove (\ref{19-4-22})  consider first the case $l>1$. Then plugging (\ref{19-4-35})--(\ref{19-4-36}) we see that for large enough $m$ we can skip $R_m$. Also we can skip term with $p=0$ because it involves only $\eta$-admissible operators. Applying the above arguments we need to estimate
\begin{equation}
|F_{t\to h^{-1}\tau} \chi_T(t) ( \Gamma v_p\psi _y)|
\label{19-4-39}
\end{equation}
with $v_p$ a Schwartz kernel of the term in (\ref{19-4-35}) with the
same notations as above. Let us consider $p=1$ first. One can see easily that
\begin{multline}
F_{t\to h^{-1}\tau}\Bigl(\chi_T(t)  e^{ith^{-1}A_\eta } \Bigr)=\\ 
\begin{aligned}
&(2\pi )^{-1} \int \chi_T(t) e^{ith^{-1}(A_\eta -\tau)} \,dt=\\
&(2\pi )^{-1} \int \chi_T(t) \bigl(-h^2\partial_t^2 +\gamma^2 \bigr)^q
e^{ith^{-1}(A_\eta -\tau)}\bigl((A_\eta  -\tau )^2+ \gamma ^2)^{-q}\,dt=\\
&(2\pi )^{-1} \bigl(\frac{h} {T}\bigr)^q\int \chi_{(q),T}(t)
e^{ith^{-1}(A_\eta -\tau)}\bigl((A_\eta  -\tau )^2+\gamma ^2)^{-q} \,dt
\end{aligned}
\label{19-4-40}
\end{multline}
with $\chi_{(q)}(t)=  (-\partial_t +\gamma ^2T^2h^{-2})^q\chi (t)$. Let us pick up $\gamma = \rho_0$.

Let us multiply (\ref{19-4-40}) from the left by $\varphi'_\nu$ which is equal to $1$ in $\epsilon T$-vicinity of $\supp \varphi_\nu$. Then modulo negligible we can multiply it from the right by $\tilde{\varphi}_\nu$ which is equal to 0 in $\epsilon T$-vicinity of $\supp \varphi_\nu$:
\begin{multline}
F_{t\to h^{-1}\tau}\bigl(\chi_T(t)  \varphi'_\nu\,  e^{ith^{-1}A_\eta } \bigr)
\equiv \\[3pt]
(2\pi )^{-1}  \bigl(\frac{h}{T}\bigr)^q\int \chi_{(q),T}(t)
\varphi'_\nu \, e^{ith^{-1}(A_\eta  -\tau)}
\bigl((A_\eta  -\tau )^2+\gamma ^2\bigr)^{-q} \tilde{\varphi}_\nu\,dt
\label{19-4-41}
\end{multline}
and therefore
\begin{multline}
F_{t\to h^{-1}\tau}\bigl(\chi_T(t) \varphi'_\nu\,e^{ith^{-1}A_\eta}\bigr) B \varphi_\nu\equiv \\[3pt]
(2\pi )^{-1}  \bigl(\frac{h}{T}\bigr)^q\int \chi_{(q),T}(t)
\varphi'_\nu\, e^{ith^{-1}(A_\eta  -\tau)}\tilde{\varphi}_\nu
\bigl((A_\eta  -\tau )^2+\gamma ^2\bigr)^{-q} \tilde{\varphi}'_\nu \cdot
\underbrace{\tilde{\varphi}_\nu B  \varphi_\nu }_{\text{sandwich}} \,dt
\label{19-4-42}
\end{multline}
where the origin of both factors $\tilde{\varphi}'_\nu$ is clear: the first one
is again due to the propagation  and the second one is just taken of
$\tilde{\varphi}_\nu$.

Let us note that  $\chi_{(q)}$ is a  function bounded by
$C\bigl(1+\varepsilon T/h\bigr)^{2q}$. Further, note that we already proved that the operator norm of the ``sandwich'' does not exceed
$C\gamma ^{lp}|\log h|^{-p\sigma}\bigl(1+\varepsilon T/h\bigr)^{-s}$ and therefore (\ref{19-4-22}) for term with $p=1$ follows from estimate
\begin{equation}
\sum _\nu
\| \tilde{\varphi}'_\nu
\bigl((A_\eta   -\tau )^2+ \gamma ^2\bigr )^{-q} \tilde{\varphi}'_\nu \|_1
\le Ch^{-r}\gamma ^{1-2q}
\label{19-4-43}
\end{equation}
where $\|.\|_1$ means the trace norm.

\medskip\noindent
(c) Let us prove (\ref{19-4-43}). Because operator is positive, its left-hand expression does not exceed
\begin{equation}
\Tr \bigl(\bigl((A_\eta  -\tau )^2+ \gamma ^2\bigr )^{-q} \psi \bigr).
\label{19-4-44}
\end{equation}
Let us consider $\gamma$-admissible partition; recall that
$\gamma \ge \varepsilon$. Because in $B(z,\gamma)$ with $z=(x,\xi)$ variation of symbol of $A_\eta$ does not exceed $C\gamma$, we conclude that (\ref{19-4-44}) does not exceed
\begin{equation*}
Ch^{-r} \int_\cD   \Tr' \bigl((a(z) -\tau)^2 +\gamma^2)^{-s} \,dz =
C\mu^{-r}h^r \int_\cD \sum_J \bigl((\omega_J-\tau)^2 +\gamma^2)^{-s} \,dz
\end{equation*}
where $a(z)$ is a (matrix) symbol of $A_\eta$  and $\Tr'$ means a (matrix) trace, $\omega_J$ are eigenvalues of $a$ (in order) and integrals are taken over bounded domain in  $\cD \subset \bR^{2r}_{x'',\xi''}$.

Then (\ref{19-4-43}) follows from
\begin{equation}
\int (|\omega _j(z) -\tau | +\gamma)^{-q} \,dz \le C\gamma^{1-2q}
\label{19-4-45}
\end{equation}
as $q>\frac{1}{2}$ which, in turn, follows from the microhyperbolicity. Really, without any loss of the generality one can assume that  
$\ell = \partial _{z_1}$. Then uniformly with respect to  $z'=(z_2,\ldots,z_d)$
$\epsilon \le \partial_{z_1}\omega_j \le C$ which instantly yields that even if we restrict ourselves by integral over $z_1$, (\ref{19-4-45}) still holds.

So, \emph{for a term with $p=1$ estimate \textup{(\ref{19-4-22})} is proven\/}. For $p\ge 2$ we need to modify our trace trick. Consider $(p+1)$-dimensional integrals representing corresponding terms
\begin{multline}
(2\pi )^{-1} h^{-p} \int_{\{t\ge t_1+\ldots+t_p\}} \chi_T(t)
e^{i h^{-1}\bigl((t-t_1-\ldots-t_p)A_\eta -t\tau\bigr)}
B  \times \\
e^{it_1h^{-1}A_\eta } B \cdots B  e^{it_kh^{-1}A_\eta }\,dt_1\cdots\,dt_p\,dt.
\label{19-4-46}
\end{multline}
Substituting
\begin{equation*}
e^{i h^{-1}\bigl((t-t_1-\ldots-t_p)A_\eta -t\tau\bigr)}
=(-h^2\partial^2_t + \gamma^2) (A_\eta^2 +\gamma^2)^{-1}
e^{i h^{-1}\bigl((t-t_1-\ldots-t_p)A_\eta -t\tau\bigr)}
\end{equation*}
and integrating one time by parts with respect to $t$ we gain factor\newline
$(h/T)^2\bigl((A_\eta -\eta)^2 + \gamma^2\bigr)^{-1}$ with $\chi$, replaced  by
$\chi_{(1)}=(-\partial^2_t + h^{-2}T^2\gamma^2) \chi $,  but we also get  new terms
\begin{multline}
ch^{1-p}  \int \chi_T(t_1+\ldots+t_p)
\bigl((A_\eta -\eta)^2 + \gamma^2\bigr)^{-1}\times\\
e^{-ih^{-1}t\tau } B e^{it_1h^{-1}A_\eta }B \cdots 
B  e^{it_ph^{-1}A_\eta }\,dt_1\cdots\,dt_p
\label{19-4-47}
\end{multline}
and may be
\begin{multline}
ch^{2-p}  \int \chi_T(t_1+\ldots+t_{p-1})
\bigl((A_\eta -\eta)^2 + \gamma^2\bigr)^{-1}\times\\
e^{-ih^{-1}t\tau } B e^{it_1h^{-1}A_\eta }B \cdots 
B  e^{it_{p-1}h^{-1}A_\eta }\,dt_1\cdots\,dt_{p-1} 
\label{19-4-48}
\end{multline}
with some other functions $\chi$. These new terms are just lesser-dimensional integrals of the same type. As before we frame operator exponents by $\varphi_*$ functions before applying this trick.

\medskip\noindent
(d) Now let $l=1$, $\sigma >1$. We need to consider $R_m$ and generated by its expression which is $(m+1)$-dimensional integral
\begin{multline}
(2\pi )^{-1} h^{-m} \int_{\{t\ge t_1+\ldots+t_m\}} \chi_T(t)
e^{i h^{-1}\bigl((t-t_1-\ldots-t_m)A_\eta -t\tau\bigr)}
B  \times \\
e^{it_1h^{-1}A_\eta } B \cdots B  e^{it_mh^{-1}A }\,dt_1\cdots\,dt_m\,dt.
\label{19-4-49}
\end{multline}
We treat it in the same way as before, but in the end we
apply estimate
\begin{equation}
\| (A_\eta+i\gamma )^{-1}\psi \|_1 \le Ch^{-r}\left\{\begin{aligned}
\gamma ^{1-q}\qquad &\text{as\ \ }q>1,\\
|\log \gamma|\qquad &\text{as\ \ }q=1,
\end{aligned}\right.
\label{19-4-50}
\end{equation}
which is proven the same way (but simpler) as (\ref{19-4-43});
we see that (\ref{19-4-49}) does not exceed  $Ch^{1-r}
\bigl(\eta |\log h|^{-\sigma} Th^{-1}\bigr)^m |\log h| =
Ch^{1-r}  |\log h|^{m(1-\sigma)+1} $ as $T\ge h^{1-\delta}$;
for large $m$ it is less than the right-hand expression in (\ref{19-4-24}).

Further, one can treat $m$-dimensional integral (\ref{19-4-48})  in the same way.
So, \emph{estimates \textup{(\ref{19-4-22})} and  \textup{(\ref{19-4-24})} are proven completely\/}.
\end{proof}

\section{Short-term estimates}
\label{sect-19-4-4}

So, let us consider dynamics generated by our reduced operator as 
\begin{equation}
T:h^{1-\delta}\le T\le T_*\Def C\mu h |\log h|
\label{19-4-51}
\end{equation}
with arbitrarily small exponent $\delta >0$.
Note that we need to consider only $T \ge \bar{T} \Def \epsilon_0\mu^{-1}$ because for $|t|\le \bar{T}$ standard theory takes place. Another restriction from below will appear later.

Recall that as $T\ll 1$ we do not need to distinguish microhyperbolicity and $\fN$-microhyperbolicity.

\begin{proposition}\label{prop-19-4-9} 
Let microhyperbolicity condition be fulfilled.  Let $\phi$ be supported in the small vicinity of $0$.

\medskip\noindent
(i) Let $l>1$. Then for $T\in [Ch^{1-\delta}, T_*]$ with arbitrarily small exponent $\delta >0$
\begin{align}
&|\phi (hD_t) \chi_T(t) ( \Gamma U\psi )|\le
C h^{-d} \bigl(\frac {h}{T}\bigr)^{l-1} |\log \frac{h} {T}|^{-\sigma}
\bigl(1+\frac{T\varepsilon} {h} \bigr)^{-s},\label{19-4-52}\\[3pt]
&|F_{t\to h^{-1}\tau} \chi_T(t) ( \Gamma U\psi )|\le
C h^{1-d} \bigl(\frac{h} {T}\bigr)^{l-1} |\log \frac{h} {T}|^{-\sigma}
\bigl(1+\frac{T\varepsilon} {h} \bigr)^{-s}
\label{19-4-53}
\end{align}
as $|\tau|\le \epsilon$ with arbitrarily large exponent $s$.

\medskip\noindent
(ii) Let $l=1,\sigma \ge 2$. Then for $T\in [Ch^{1-\delta}, T_*]$ with arbitrarily small exponent $\delta >0$ estimates
\begin{align}
&|\phi (hD_t) \chi_T(t) ( \Gamma U\psi )|\le\label{19-4-54}\\
&\qquad\qquad C h^{-d} \bigl(\frac{h} {T}\bigr)^{l-1} |\log \frac{h} {T}|^{-\sigma}
\bigl(1+\frac {T\varepsilon}{h} \bigr)^{-s}+Ch^{-d}|\log \frac{h} {T}|^{-s}, \notag\\[3pt]
&|F_{t\to h^{-1}\tau} \chi_T(t) ( \Gamma U\psi )|\le\label{19-4-55}\\
&\qquad\qquad C h^{1-d}\bigl( \frac{h} {T}\bigr)^{l-1} |\log \frac{h} {T}|^{-\sigma}
\bigl(1+\frac {T\varepsilon}{h} \bigr)^{-s}+Ch^{1-d} |\log \frac{h} {T}|^{-s}
\notag
\end{align}
as $|\tau|\le\epsilon $ with arbitrarily large exponent $s$.
\end{proposition}

\begin{proof}
As in the previous Subsection~\ref{sect-19-4-3}, we can consider $\chi $ supported in $[\frac {1}{2},1]$ and in what follows $t,T_k$ are positive.

Consider $\eta $-mollification $\cA_\eta$ of operator $A_\cT$ with respect to $x'',\xi''$  with $\eta$ satisfying (\ref{19-4-28}). Then the  error in the operator (i.e. the difference between $A_\eta$ and $A$)  will not exceed
$C\eta^l|\log \eta|^{-s}$ in the sense that
\begin{equation}
\| (A_\eta -A) v\|\le
C\eta^l|\log \eta|^{-s}\Bigl( \|Av\| +\|v\|\Bigr),
\label{19-4-56}
\end{equation}
and then we just follow the proof proposition \ref{prop-19-4-8} with the following modifications:

\medskip\noindent
(a) we consider $\mu^{-1}h$-pseudo-differential operators with respect to $x''$
with ``matrix'' symbols $a$; furthermore, factor  $(\mu h)^{-r}$ comes from the ``matrix'' trace.

\medskip\noindent
(b) Between different copies of $\cB$ we place
not only $\varphi_{(\lambda)\nu j}$  but also
$Q'_j$ where $Q'_j$ are cut-offs with respect to $(x',\xi')$ keeping all in zone $\{ |x'|+ |\xi'|\le C\mu^{-1}\}$ (and also $a_{\fm}$ in vicinity of
$\tau_{\fm}$ for all ${\fm} \in \fM$).

\medskip\noindent
(c) We remember that
$\psi_\cT = \psi (x',\mu^{-1}D'; x'',\mu^{-1}hD'')$ with ${\sF}^{l,\sigma}$-symbol and we need to mollify   $\psi$ as well. Then there will be either
$\psi_{\cT\eta}$ or $(\psi_{\cT\eta}-\psi_\cT)$ factor; in the latter case it is possible that it is the only factor containing difference. In the former case we
follow Subsection~\ref{sect-19-4-2}, in the latter case we just note that the norm of $(\psi_{\cT\eta}-\psi_\cT)$ is $O(\eta^l |\log \eta|^{-\sigma})$ exactly as in $(A_{\cT\eta}-A_\cT)$ but this time it is not accompanied by a rather large factor $T/h$ and this makes things much more comfortable.
\end{proof}

Then summing over partition of $[-T_*,T_*]\setminus [-T'_*,T'_*]$ with $T'_*=Ch^{1-\delta}$    we get immediately inequality
\begin{equation}
|F_{t\to h^{-1}\tau} \bigl(\bar{\chi}_{T_*}(t)- \bar{\chi}_{T'_*}(t)\bigr)
\Gamma (U \psi) |\le Ch^{1-d}.
\label{19-4-57}
\end{equation}

Really, the left-hand expression here does not exceed the right-hand expression of (\ref{19-4-52}) integrated over $\frac{dT}{T}$ from $T'_*$ to $T_*$ which, in turn, does not exceed $Ch^{1-d}$; contribution of extra term (\ref{19-4-54}) as $l=1$ does not exceed $Ch^{1-d}$ as well. On the other hand, we know from rescaling of the standard results that (for $\mu \le \mu^*_2$)
\begin{equation}
|F_{t\to h^{-1}\tau}  \bar{\chi}_T(t) \Gamma(U \psi) |\le Ch^{1-d}
\label{19-4-58}
\end{equation}
as $T=T'_*$ and therefore it holds for $T=T_*$ and it follows from Section~\ref{sect-19-2}, that (\ref{19-4-58}) holds for $T=T^*$. Recall that under microhyperbolicity assumption $T^*=\epsilon$ and under $\fN$-microhyperbolicity assumption $T^*=\epsilon \mu$ as \underline{either} $\#\fN=1$ \underline{or}  $(l,\sigma)\succeq (2,0)$ and
$T^*=\epsilon \mu^{l-1}|\log \mu|^\sigma$ otherwise.

Then standard the Tauberian arguments imply immediately

\begin{corollary} \label{cor-19-4-10} 
Under microhyperbolicity condition  as $\mu \le h^{\delta-1}$ estimate \textup{(\ref{19-2-78})} holds.
\end{corollary}

In the  case $h^{\delta -1}\le \mu \le \mu^*_2$ we will need analysis of the next Subsection~\ref{sect-19-4-5} to prove this statement. Anyway, we need it to get more explicit formula than (\ref{19-2-78}) provides.

\section{Calculations}
\label{sect-19-4-5}
In this Subsection we still under assumption (\ref{19-3-5}) will replace rather implicit Tauberian expression formula  by a more explicit one and we finish the proof the crucial inequality
\begin{equation}
|F_{t\to h^{-1}\tau}\Gamma  \bar{\chi}_{T_*}(t) U \psi |\le
Ch^{1-2r}+C\mu^rh^{1-r}.
\label{19-4-59}
\end{equation}
Our tool will be method of successive approximations on a rather short interval
$[-T_*,T_*]$.

\subsection{Preliminary remarks}
\label{sect-19-4-5-1}

Now even if we have a remainder estimate proven in some cases (see corollary~\ref{cor-19-4-10}) but  its principal part is given by rather implicit Tauberian expression 
\begin{equation}
h^{-1}\int _{-\infty}^0 \Bigl( F_{t\to  h^{-1}\tau}
\bar{\chi}_T(t) \Tr \bigl(e^{ih^{-1}tA_\cT }\psi_\cT \bigr) \Bigr)\, d\tau
\label{19-4-60}
\end{equation}
with arbitrary $T\ge T_*=C\varepsilon^{-1}h|\log h|$ and $\psi_\cT=\cT^*\psi\cT$. 

Surely, microhyperbolicity condition is fulfilled only at levels close to $0$ but as usual we do not need it otherwise. Really, let us decompose expression (\ref{19-4-60}) into sum of
\begin{align}
&\int _{-\infty}^0 \bar{\phi}_L (\tau)\Bigl( F_{t\to  h^{-1}\tau}
\bar{\chi}_T(t) \Tr \bigl(e^{ih^{-1}t\cA }\psi_\cT \bigr) \Bigr)\, d\tau \label{19-4-61}\\
\shortintertext{and}
&\int _{-\infty}^\infty \phi_L(\tau)\Bigl( F_{t\to  h^{-1}\tau}
\bar{\chi}_T(t) \Tr \bigl(e^{ih^{-1}t\cA }\psi_\cT\Bigr) \Bigr)\, d\tau
\label{19-4-62}
\end{align}
where $\bar{\phi}$, $\phi$ are admissible functions, $\bar{\phi}$ is supported in $[-1,1]$ and equal $1$ in $[-\frac{1}{2},\frac{1}{2}]$,
$\phi=\uptheta (-\tau) \bigl(1-\bar{\phi}(\tau)\bigr)$, $L$ is a small constant.
Now we can use microhyperbolicity in the first term and replace $T=T^*$ by any
$T\in [T_*,T^*]$. On the other hand, one can  rewrite the second expression as
\begin{equation}
\Bigl(\phi_L(hD_t)\bigl(\bar{\chi}_T(t) \Tr \bigl(e^{ih^{-1}t\cA }\psi_\cT \bigr) \bigr)\bigr)\Bigr|_{t=0}
\label{19-4-63}
\end{equation}
and then one can replace here $T=T^*$ by any $T\in [T',T^*]$ where $T'=CL^{-1}h|\log \mu |\ll T_*$.
One can see easily that expression (\ref{19-4-63}) is of magnitude
$h^{-d} + \mu^r h^{-r}$ and therefore we expect  the principal part of asymptotics  of the same magnitude.
We assume that microhyperbolicity is in the direction $\partial_{\xi_1}$.

Recall that
\begin{multline}
\int_{-\infty}^0
h^{-1}\Bigl(F_{t\to h^{-1}\tau }\chi_T(t)
\Gamma \bigl(U\psi \bigr)\Bigr)\,d\tau =\\
T^{-1} \Bigl(F_{t\to h^{-1}\tau }\check{\chi}_T (t)
\Gamma \bigl(U\psi \bigr)\Bigr)\Bigr|_{\tau=0}
\label{19-4-64}
\end{multline}
with $\check{\chi}(t)=i t^{-1}\chi (t)$; this formula as usual  plays very important role below. 

\begin{remark}\label{rem-19-4-11} 
It follows from propositions \ref{prop-19-4-9}, \ref{prop-19-5-10}--\ref{prop-19-6-20} that in their framework expression (\ref{19-4-64}) does not exceed
$Ch^{-d} \bigl(h/T\bigr)^l|\log \bigl(h/T\bigr)|^{-\sigma}$
and thus
\begin{multline}
\int_{-\infty}^0
h^{-1}|\Bigl(F_{t\to h^{-1}\tau }(\bar{\chi}_{T_*}(t) - \bar{\chi}_T(t))
\Gamma \bigl(U\psi \bigr)\Bigr)\,d\tau |\le \\
C\bigl(h^{-d}+\mu^rh^{-r}\bigr) \bigl(\frac{h}{T}\bigr)^l
|\log \bigl(\frac{h}{T}\bigr)|^{-\sigma}.
\label{19-4-65}
\end{multline}
In particular for $T=\bar{T}_*\Def \varepsilon^{-1}h$ the right-hand expression does not exceed the remainder estimates we want to prove (in any of our cases). We can even reduce $T$ if the smoothness allows. 

In particular, as $T=\epsilon \mu^{-1}$ (and therefore $\varepsilon=\mu h$) we get a remainder estimate  $Ch^{-d}(\mu h)^l|\log \mu|^{-\sigma}$ which is a bit better than in theorem \ref{thm-19-2-16}.
\end{remark}

\subsection{Successive approximation method}
\label{sect-19-4-5-2}

We apply method of successive approximations to calculate both terms (\ref{19-4-61}) and (\ref{19-4-62}). We consider our operator 
$A_\cT (x',hD';x'',hD_{x''})$ as a matrix operator and then we take
$A_\cT (x',hD',y'',hD_{x''})$ for unperturbed operator $\bar{A}_\cT$ (simplifying it later). 

Then one can see easily that each next term of successive approximations gets an extra factor $\bar{\cG}^\pm R$ or $\cG ^\pm R$ where  $\bar{\cG}^\pm$ and 
$\cG ^\pm$ are forward and backward parametrices of $hD_t-\bar{A}_\cT$ and 
$hD_t -A_\cT $ respectively and by Duhamel principle their operator norms on interval $[-T,T]$ do not exceed $Th^{-1}$.

\subsection{Mollified by $\tau$ asymptotics}
\label{sect-19-4-5-3}
Each next term of successive approximations when plugged there adds a factor $\|R\|$ to its estimate from above.

Really, it follows from the fact that
$t^{-1}\bar{\cG}^\pm$ and $t^{-1}\cG ^\pm$ are operators with norms $h^{-1}$ and factor $t$ either annuls the corresponding term, or is replaced by its commutator with $\phi_L (hD_t)$ thus releasing factor $h$. So, with expression (\ref{19-4-63}) any perturbation with norm $O(h^\delta)$ would be good enough to have only a bounded number of successive approximations to be considered, but $(x_j-y_j)$ is even better: as usual
\begin{equation}
\bar{\cG }^\pm (x_j-y_j)=(x_j-y_j)\bar{\cG}^\pm  -
\bar{\cG }^\pm [\bar{A}, x_j] \bar{\cG}^\pm
\label{19-4-66}
\end{equation}
and similarly for $\cG ^\pm$, and as explained, these parametrices plugged
into expression (\ref{19-4-63}) do not increase an upper estimate, while the norm of $[\bar{A}, x_j]$ is $O(\mu ^{-1}h)$.

Further, each perturbation of the form $(A_\cT -\cA)$ adds a factor
$\mu^{-l}|\log h|^{-\sigma}$ to the corresponding term in (\ref{19-4-63});
we take $\bar{\cA}$ as unperturbed operator in this Subsubsection.

So any term of the successive approximation but the  first one,
leads to term in (\ref{19-4-63}) not exceeding
$\bigl(\mu^{-1}h +\mu^{-l}|\log h|^{-\sigma}\bigr) \times h^{-d}$ which, in turn, does not exceed the  remainder estimate.  Therefore

\begin{claim}\label{19-4-67}
Modulo remainder estimate one can replace in (\ref{19-4-63}) operator $A_\cT$ by $\bar{\cA}$ leading to expression
\begin{multline}
(2\pi)^{-r}\mu ^r h^{-r} \times\\
\int
\int {\bar\phi}_L(\tau)\, d_\tau \biggl(\Tr ' \Bigl(
\uptheta \bigl(\tau -\cA (x'',\xi '') \bigr)
\psi_\cT (x'',\xi'')\Bigr)\biggr)\,dx''d\xi''
\label{19-4-68}
\end{multline}
where both $\cA (x'',\xi '')$ and $\psi_\cT (x'',\xi '')$ are considered as  ``matrices'' i.e. as operators in auxiliary space $\bH=\sL^2(\bR^r)$,
$\uptheta \bigl(\tau -\cA (x'',\xi '') \bigr)$ is its spectral projector  and $\Tr'=\Tr_\bH$ is a ``matrix'' trace.
\end{claim}

\subsection{Unmollified asymptotics. I}
\label{sect-19-4-5-4}

Unfortunately things are not that good for  expression (\ref{19-4-61}) because
parametrices are ``worth'' of $T/h$ and therefore in view of (\ref{19-4-66}) factor $(x_j-y_j)$ is ``worth'' of $T^2h^{-2}\times \mu ^{-1}h\asymp
\varepsilon ^{-2}\mu ^{-1}h |\log h|^2\asymp \mu h|\log h|^2$ (as $T=T_*$) which is not that small,  especially for large $\mu$.

In $2\D$-case considered in the previous Chapter~\ref{book_new-sect-18} of \cite{futurebook} we could always arrange (after an appropriate symplectic map\footnote{\label{foot-19-25} Which after quantization gives us a metaplectic transformation.}) to have $(x_j-y_j)$ accomplished by an extra factor
$(\varepsilon + \varepsilon^{l-1}|\log \varepsilon|^{-\sigma})+\varepsilon^2$ because we considered a scalar symbol. 

Here we have essentially a matrix symbol and this trick does not work (however it will work for a superstrong magnetic field). A bit larger $\varepsilon$ does not help much as $\mu$ is close to $h^{-1}$.
Further, perturbation $(A_\cT-\cA)$ is ``worth'' of
$Th^{-1} \times \mu^{-l}|\log \mu|^{-\sigma }\asymp \varepsilon^{-1}\mu^{-l} |\log h|^{1-\sigma}$ (again as $T=T_*$). In $2\D$-case we just removed such term at the very beginning.

However, using methods of Subsection~\ref{sect-19-4-4} we will be able to insert a factor $\bigl(h/T\bigr)^l |\log (h/T)|^{-\sigma}$ into estimates of successive the approximation terms which will take care of few copies of factor $Th^{-1}$.

So, our goal now is to evaluate properly
\begin{equation}
F_{t\to h^{-1}\tau }\chi_T(t) \Gamma \bigl(U\psi \bigr) =
F_{t\to h^{-1}\tau }\chi_T(t) \Gamma \bigl(U_\cT\psi_{\cT}\bigr)
\label{19-4-69}
\end{equation}
as $|\tau|\le \epsilon$.

Recall the successive approximation method is the following procedure: rewriting down
\begin{equation}
(hD_t-A_\cT)U^\pm _\cT = \mp i h\updelta(t)\updelta (x-y)
\label{19-4-70}
\end{equation}
with $U^\pm = U\uptheta (\pm t)$ and thus
\begin{gather}
U^\pm _\cT = \mp i h \cG^\pm \updelta(t)\updelta (x-y)
\label{19-4-71}\\
\shortintertext{as}
(hD_t-\bar{A}_\cT)U^\pm _\cT = \mp i h\updelta(t)\updelta (x-y) +
(A_\cT - \bar{A}_\cT) U^\pm
\label{19-4-72}
\end{gather}
with $U^\pm = U\theta (\pm t)$ and thus
\begin{equation}
U^\pm _\cT = \mp i h \bar{\cG}^\pm \updelta(t)\updelta (x-y) +
\bar{\cG}^\pm(A_\cT - \bar{A}_\cT) U^\pm
\label{19-4-73}
\end{equation}
we can iterate the last equation few times:
\begin{multline}
U^\pm _\cT = \mp i h \sum_{0\le k\le m-1}
\Bigl(\bar{\cG}^\pm(A_\cT - \bar{A}_\cT) \Bigr)^k \bar{\cG}^\pm
\updelta(t)\updelta (x-y)  \\
\mp i h \Bigl(\bar{\cG}^\pm(A_\cT - \bar{A}_\cT) \Bigr)^m \cG^\pm
\updelta(t)\updelta (x-y);
\label{19-4-74}
\end{multline}
however our latest technique allows usually to take $m=1$. Applying $\psi_\cT$ from the right we get that
\begin{multline}
U^\pm _\cT = \mp i h \sum_{0\le k\le m-1}
\Bigl(\bar{\cG}^\pm(A_\cT - \bar{A}_\cT) \Bigr)^k \bar{\cG}^\pm
\updelta(t)\cK_\psi (x,y)  \\
\mp i h \Bigl(\bar{\cG}^\pm(A_\cT - \bar{A}_\cT) \Bigr)^m \cG^\pm
\updelta(t)\cK_\psi (x,y)
\label{19-4-75}
\end{multline}
where $\cK_\psi (x,y)$ is the Schwartz kernel of $\psi_\cT$. 

\begin{claim}\label{19-4-76}
From now on $\Gamma$ means that we set both Schwartz kernel arguments equal to $y'$ (the same we used in the successive approximations) and then integrate.
\end{claim}

\begin{proposition}\label{prop-19-4-12} 
In the framework of proposition~\ref{prop-19-4-9} for $T\in [T', T^*]$, 
$T'\Def h^{1-\delta}$,
\begin{multline}
T^{-1} |F_{t\to h^{-1}\tau }\chi_T(t) \Gamma
\Bigl({\bar\cG}^\pm(A_\cT - \bar{A}_\cT) \cG^\pm
\delta(t)\cK_\psi\Bigr)|\le \\
C\mu^{-1}h^{1-d}  \bigl(\frac{h}{T}\bigr)^{l-2}|\log h|^{-\sigma}.
\label{19-4-77}
\end{multline}
Therefore
\begin{multline}
 |\int_{-\infty}^0 F_{t\to h^{-1}\tau }
\bigl(\bar{\chi}_{T^*}(t)-\bar{\chi}_{T'}(t) \bigr)\Gamma
\Bigl(\bar{\cG}^\pm(A_\cT - \bar{A}_\cT) \cG^\pm
\updelta(t)\cK_\psi\Bigr)\, d\tau |\le \\
C\left\{\begin{aligned}
&\mu^{-1}h^{1-d}  \qquad\qquad\qquad &\text{as\ \ }(l,\sigma)\succeq (2,1),\\
&\mu^{1-l}h^{1-d}|\log h|^{-\sigma} &\text{as\ \ }(l,\sigma)\prec (2,1)
\end{aligned}\right.
\label{19-4-78}
\end{multline}
and thus does not exceed remainder estimate of theorem \ref{thm-19-4-17} below.
\end{proposition}

\begin{proof} Rewriting
\begin{equation}
A_\cT - \bar{A}_\cT =
\sum_j (x_j-y_j)R_j,\qquad R_j= R_j (x',\mu^{-1}hD'; x'',y'',\mu^{-1}hD'')
\label{19-4-79}
\end{equation}
we see that
\begin{multline}
\Bigl(\bar{\cG}^\pm(A_\cT - \bar{A}_\cT) \cG^\pm
\updelta(t)\cK_\psi \Bigr) =\\
\sum_j \Bigl(\bar{\cG}^\pm R_j \cG^\pm
\updelta(t)\cK_{\psi,j} \Bigr)
-i\mu^{-1}h  \sum_j \Bigl(\bar{\cG}^\pm R_j
\cG^\pm A_\cT ^{(j)} \cG^\pm
\updelta(t)\cK_\psi\Bigr)
\label{19-4-80}
\end{multline}
where we used standard commutation relation (\ref{19-4-66}) and
$\cK_{\psi j}=(x_j-y_j)\cK_\psi$ is the Schwartz kernel of $[x_j,\psi]$. Note that symbols of $\bar{A}_\cT^{(j)}$, $R_j$ belong to $\sF^{l-1,\sigma}$
while symbol of $[x_j,\psi]$ belongs to $\mu^{-1}h\sF^{l-1,\sigma}$.

Therefore due to the Duhamel formula the left-hand expression of (\ref{19-4-77}) is an absolute value of the sum of
\begin{multline}
-i\mu^{-1}h T^{-1} F_{t\to h^{-1}\tau} \chi_T(t) \Tr
\sum_j \Bigl(\bar{\cG}^\pm R_j
\cG^\pm A_\cT ^{(j)} \cG^\pm\psi_\cT\Bigr)=\\
\shoveleft{-i\mu^{-1}h^{-1}T^{-1} F_{t\to h^{-1}\tau} \chi_T(t)  \times}\\[3pt]
\sum_j\Tr\Bigl(\int_{\{t_1+t_2\le t\}}
e^{ih^{-1}(t-t_1-t_2)\bar{A}_\cT} R_j e^{ih^{-1}t_1A _\cT}A_\cT ^{(j)}
  e^{ih^{-1}t_2A_\cT}\psi_\cT \,dt_1\,dt_2\Bigr)
\label{19-4-81}
\end{multline}
and
\begin{multline}
T^{-1} F_{t\to h^{-1}\tau} \chi_T(t) \Tr
\sum_j \Bigl(\bar{\cG}^\pm R_j
\cG^\pm [x_j,\psi_\cT]\Bigr)=\\
T^{-1} F_{t\to h^{-1}\tau} \chi_T(t) \sum_j\Tr
\Bigl(\int_{\{t_1 \le t\}}
e^{ih^{-1}(t-t_1)\bar{A}_\cT} R_j e^{ih^{-1}t_1A _\cT}[x_j,\psi_\cT ] \,dt_1\Bigr).
\label{19-4-82}
\end{multline}

For $\eta = ChT^{-1}|\log h|$ let us make $\eta$-mollification $A_{\cT\,\eta}$ of (symbol) of $A_\cT$ and also of $\psi$; thus we will get also
$\bar{A}_{\cT\,\eta}$. Note that the left-hand expression of (\ref{19-4-77}) with $A_\cT$, $\bar{A}_\cT$ replaced by $A_{\cT\,\eta}$, $\bar{A}_{\cT\,\eta}$
is negligible again due to proposition \ref{prop-19-4-3}. Really, both operators
$A_{\cT\,\eta}$, $\bar{A}_{\cT\,\eta}$ are microhyperbolic in the same direction.

Consider the difference between expressions (\ref{19-4-80}) for $\eta$-mollified and original operators.

To do this we apply formula (\ref{19-4-35})--(\ref{19-4-36}) assuming first that $l>1$ and skipping the negligible remainder. Then the difference in  (\ref{19-4-81}) will be operator
\begin{equation*}
-i\mu^{-1}h^{-1-k} F_{t\to h^{-1}\tau} \chi_T(t) \sum_j\Tr 
\end{equation*}
applied to the sum of the following terms (with the constant coefficients)
\begin{equation*}
\mu^{-1}h^{-1-k-k_1-k_2}\int_{\Delta_k}
e^{ih^{-1}(t-t_1-t_2-\ldots-t_k)\bar{A}_{\cT\eta}}
\cC_1 e^{ih^{-1}t_1\bar{A}_{\cT\eta}}\cC_2 \cdots \cC_k
e^{ih^{-1}t_k\bar{A}_{\cT\eta}}\psi'
\end{equation*}
where operators $\cC_k$ are in a some order:
\begin{enumerate}[leftmargin=*,label=(\alph*)]
\item one copy of $\bar{A}_{\cT\eta}^{(j)}$ or
$(\bar{A}_{\cT\eta}^{(j)}-A_{\cT\eta}^{(j)})$,
\item one copy of  $R_{j\eta}$ or $(R_{j\eta}-R_j)$ and
\item $(\bar{A}_\cT- \bar{A}_{\cT\eta})$ constitute the rest; 
\end{enumerate}also
$\psi '= \psi_{\cT\eta}$ or $\psi'=(\psi_\cT-\psi_{\cT\,\eta})$
and among all the factors  there must be at least one factor which is the difference between mollified and non-mollified operators.

Acting as in the proof proposition \ref{prop-19-4-8} we can also gain the factor $h^{2q}T^{-2q}\bigl((\bar{A}_{\cT\,\eta}-\tau) ^2 +\gamma^2\bigr)^{-q}$ replacing $\chi$ by $(-\partial_t^2 + \gamma^2T^2h^{-2})^q\chi$. Further, $2\D$-integral can also appear instead of $3\D$ one (one extra integral comes from the Fourier transform).

Note that 
\begin{enumerate}[leftmargin=*,label=(\alph*)]
\item each  factor $(\bar{A}_\cT- \bar{A}_{\cT\eta})$ or
$(A_\cT- A_{\cT\eta})$  is worth of $\eta ^l |\log h|^{-\sigma}$ and is accompanied by a factor not exceeding $Th^{-1}$ while 
\item each factor
$(\bar{A}_\cT^{(j)}- \bar{A}_{\cT\rho}^{(j)})$ or $(R_{j\eta}-R_j)$
is worth of $\eta ^{l-1} |\log h|^{-\sigma}$ and  $(\psi_\cT-\psi_{\cT\eta})$
is worth of $\eta ^l |\log h|^{-\sigma}$.
\end{enumerate}
Then using estimate (\ref{19-4-50}) with $\gamma =h/T$ one can derive easily that the difference does not exceed
$C\mu^{-1} h^{1-d}\bigl(h/T\bigr)^{l-2} |\log h|^{l-\sigma}$
and we arrive to the estimate different from the required one by an extra  factor $|\log h|^l$.

There are many methods to get rid of this factor; let us apply the cheapest one. First, let we apply $m$ iterations in successive approximations and note that each extra iteration brings the factor
$\mu^{-1}h \times T^2h^{-2}\le \mu h \le |\log h|^{-1}$ and use the  same arguments as above; we conclude that instead of $\cG^\pm $ we need to consider
$\bigl(\bar{\cG}^\pm(A_\cT - \bar{A}_\cT)\bigr)^k\bar{\cG}^\pm$; so we got rid of $\cG^\pm $ but paid for this by a bit more complicated expression. This expression includes also $[x_j-x_k, R_j]$, $[x_j-x_k, A^{(j)}]$ etc but one can see easily that the corresponding terms would satisfy required estimate.

Also, because $l>1$ we see from the above mollification arguments that all the terms with more than one factor of $(\bar{A}_\cT- \bar{A}_{\cT\eta})$,
$(\bar{A}_\cT^{(j)}- \bar{A}_{\cT\eta}^{(j)})$ and
$(R_j-R_{j\eta})$  satisfy required estimate and thus we need to consider only terms with exactly one such factor. Further, terms containing factor
$(\psi_\cT-\psi_{\cT\eta})$ satisfy this estimate for sure.

Furthermore, if we replace in such factor $\eta$ by $\rho_0=h/T$ we will get terms which can be properly estimated and thus in the terms under consideration we can replace $(\bar{A}_\cT- \bar{A}_{\cT\eta})$ by
$(\bar{A}_{\cT\rho_0}- \bar{A}_{\cT\eta})=\sum _n \cB_n$, also replace
$(\bar{A}_\cT^{(j)}- \bar{A}_{\cT\eta}^{(j)})$ by
$(\bar{A}_{\cT\rho_0}^{(j)}- \bar{A}_{\cT\eta}^{(j)})=\sum _n \cB'_n$, and replace $(R_j-R_{j\eta})$ by $(R_{j\rho_0}-R_{j\eta})= \sum _n \cB''_n$
where $\cB_n \Def (\bar{A}_{\cT\rho_n}- \bar{A}_{\cT\rho_{n+1}})$
etc where again $\rho_n =2^n\rho_0$.

Now we can just repeat arguments of the proof proposition \ref{prop-19-4-8}, related to $\mu^{-1}T$-partition with respect to $x_1$ assuming that microhyperbolicity direction is $\partial_{\xi_1}$.

On the other hand, if $l=1$, using (\ref{19-4-49}) and related integration by part trick (but without $\mu^{-1}T$-partition), one can see easily that the far terms in the successive approximation satisfy estimate in question and thus we got rid of $\cG^\pm$ leaving only $\bar{\cG}^\pm$.

We can apply mollification and formulae (\ref{19-4-35})--(\ref{19-4-36}) skipping terms generated by (\ref{19-4-36});  then we arrive to an estimate which contains in comparison with the desired one extra factor $|\log h|$ and all the terms save the same as above surely satisfy desired estimate. Now we just apply the same $\mu^{-1}T$-partition with respect to $x_1$ etc as above.

Analysis of expression (\ref{19-4-82}) uses the same technique but is simpler.
\end{proof}

\subsection{Unmollified asymptotics. II}
\label{sect-19-4-5-5}

Now we are in the shortest-range zone $|t|\le \bar{T}\Def h^{1-\delta}$. Then as 
$\mu \le \epsilon h^{\delta-1}$ we can apply rescaling arguments and arrive to estimate (\ref{19-4-59}) which implies the corresponding Tauberian estimate.

However as $\epsilon h^{\delta -1}\le \mu \le \epsilon (h|\log h|)^{-1}$ we still need to estimate \newline $F_{t\to h^{-1}\tau}\bar{\chi}_T \Gamma (U_\cT\psi_\cT)$.  Using a successive approximation method in its rough form we conclude that modulo term not exceeding  $CTh^{-d}\times \mu^{-1}h (T/h)^2 \asymp \mu^{-1}h^{1-d-3\delta}$ we can replace here $U$ by $U^0_\cT$ which is a Schwartz kernel of  $e^{ih^{-1}t \bar{\cA}}$. Then calculations show that for $T=\bar{T}$ and $\bar{\chi}=1$ on $[-\frac{1}{2},\frac{1}{2}]$
\begin{multline}
F_{t\to h^{-1}\tau}\bar{\chi}_T \Gamma (u^0_\cT\psi_\cT)\equiv \\
(2\pi )^{-r}\mu^rh^{1-r}\partial_\tau \int \Tr ' \Bigl(
\uptheta  \bigl(\tau -\cA(x'',\xi'')\bigr)\psi_\cT (x'',\xi'')\Bigr)\,dx''d\xi''
\label{19-4-83}
\end{multline}
modulo negligible term where $\uptheta  \bigl(\tau -\cA (x'',\xi'')\bigr)$ is the spectral projector of $\cA$ considered as operator in $\bH$, $\psi(x'',\xi'')$ is also operator in $\bH$ and $\Tr'$ is a trace in $\bH$.

One can see easily that due to the microhyperbolicity condition  for small increment $d\tau$ eigenvalue $\omega_J(x'',\xi'')$ of matrix $\cA (x'',\xi'')$ belongs to $[\tau, \tau+d\tau]$ only for
$(x'',\xi'')$ belonging to the set of measure not exceeding $Cd\tau$.
Since there are $\asymp (\mu h)^{-r}$ of such eigenvalues we conclude that expression (\ref{19-4-83}) does not exceed $Ch^{1-d}$. 

This and the previous arguments imply immediately estimate (\ref{19-4-57}) with $T=\bar{T}$ and thus the same estimate (\ref{19-4-57}) with $T=T_*$ and thus (\ref{19-4-57}) with $T=T^*$. Then the standard Tauberian arguments imply estimate (\ref{19-2-78}). So we have proven

\begin{proposition}\label{prop-19-4-13} 
Under either microhyperbolicity or $\fN$-microhyperbolicity condition  for  
$\mu \le \epsilon (h|\log h|)^{-1}$  estimates \textup{(\ref{19-4-57})} and \textup{(\ref{19-2-78})} hold.
\end{proposition}

Still we need to get more explicit formula:

\begin{proposition}\label{prop-19-4-14} 
Under either microhyperbolicity or $\fN$-microhyperbolicity condition for  
$\mu \le \epsilon (h|\log h|)^{-1}$
\begin{multline}
|\Gamma (\psi \tilde{e}) (\tau)-\\
(2\pi )^{-r}\mu^rh^{-r} \int \Tr ' \Bigl(
\uptheta  \bigl(\tau -A_\cT (x'',\xi'')\bigr)\psi_\cT (x'',\xi'') \Bigr)\,dx''d\xi''|\le \\
C\mu^{-1}h^{1-d}\Bigl( 1+ \varepsilon^{l-2}|\log h|^{-\sigma}\Bigr).
\label{19-4-84}
\end{multline}
\end{proposition}

\begin{proof}
In view of what we have proven already we need to  prove that
\begin{equation}
|h^{-1}\int_{-\infty}^0 F_{t\to h^{-1}\tau } \bar{\chi}_T(t)\Gamma \Bigl(\bar{\cG}^\pm(A_\cT - \bar{A}_\cT )  \cG ^\pm
\updelta(t)\cK_\psi \Bigr)\,d\tau |
\label{19-4-85}
\end{equation}
with $T=\bar{T}$ does not exceed right-hand expression of (\ref{19-4-84}).

Consider expression (\ref{19-4-85}) first. Clearly it does not exceed $C\mu^{-1}h^{-d-2}T^3$ and for $T=\bar{T}$ it is just $C\mu^{-1}h^{1-d-3\delta}$ and therefore the required estimate  holds as $l<2$; so we need to consider case $l\ge 2$ only.

Applying then the brute-force successive approximations we conclude that replacing in (\ref{19-4-85}) $\cG^\pm$ by $\bar{\cG}^\pm$ leads to a much smaller error. Also, note that 
\begin{equation*}
A_\cT - \bar{A}_\cT = \sum_j (x_j-y_j)\bar{A}_{\cT (j)} +
\sum _{j,k} (x_j-y_j)(x_k-y_k) R'_{jk}
\end{equation*}
with $R'_{jk} \in \sF^{l-2,\sigma}$. Then using commutator equality for all $(x_j-y_j)$ factors, we conclude that the contribution of these terms will be also much smaller. In the end we are left with the term
\begin{equation}
-i\mu^{-1}\sum_j F_{t \to h^{-1}\tau} \bar{\chi}_T(t)
\Tr \Bigl({\bar\cG}^\pm \bar{A}_{\cT (j)} \bar{\cG}^\pm \bar{A}_\cT ^{(j)} \bar{\cG}^\pm \psi_\cT \Bigr)
\label{19-4-86}
\end{equation}
which is reduced to
\begin{multline}
i(2\pi)^{-r}\mu^{r-1}h^{-r} \sum_j F_{t \to h^{-1}\tau} \bar{\chi}_T(t)
\Tr '\int \Bigl(\cG^\pm (x'',\xi'') A_{\cT (j)}\\
(x'',\xi'') \cG^\pm (x'',\xi'') A_\cT ^{(j)} (x'',\xi'')  \cG^\pm (x'',\xi'')  \psi^0 (x'',\xi'')  \Bigr)\, dx''d\xi''
\label{19-4-87}
\end{multline}
where we have just matrix-valued symbols $a(x'',\xi'')=A_\cT (x'',\xi'')$
and  $\cG^\pm (x'',\xi'')$ and we replaced matrix-valued function
$\psi_\cT (x'',\xi'')$ by the scalar-valued one $ \psi^0(x'',\xi'')$; this is legitimate since  $\psi_\cT (x'',\xi'') \equiv \psi^0(x'',\xi'')$ modulo $O(\mu^{-1})$.

Using that $\partial_{x_j}\cG^\pm (x'',\xi'')= \cG^\pm (x'',\xi'')A_{\cT(j)} (x'',\xi'') \partial_{x_j}\cG^\pm (x'',\xi'')$ and
similar identity for $\partial_{\xi_j}$ one can rewrite (\ref{19-4-87}) due to
the trace property and scalar nature of $\psi^0$ as
\begin{multline}
 (2\pi)^{-r} \mu^{r-1}h^{-r}
 F_{t \to h^{-1}\tau} \bar{\chi}_T(t) \times \\[2pt]
\Tr '\int \Bigl(\cG^\pm (x'',\xi'') B (x'',\xi'')   \cG^\pm (x'',\xi'')  \Bigr)\psi^0 (x'',\xi'')  \, dx''d\xi'' +\\[2pt]
\frac{i}{2} (2\pi)^{-r} \mu^{r-1}h^{-r}
 F_{t \to h^{-1}\tau} \bar{\chi}_T(t)
\Tr '\int  \cG^\pm (x'',\xi'')  \psi' (x'',\xi'')  \, dx''d\xi''
\label{19-4-88}
\end{multline}
with $B=-\frac{i}{2}\sum_j A_{\cT (j)}^{(j)}$ and $\psi'=\sum_j\psi^{0(j)}_{(j)}$.We can rewrite the first term here as the convolution with respect to $\tau$ of $L\cdot {\widehat{\bar{\chi}}}(L\tau )$ (with $L=Th^{-1}$) with 
\begin{multline}
(2\pi)^{-r} \mu^{r-1}h^{-r} \Tr' \partial_\tau\times \\
\int  \uptheta \bigl(\tau - A_\cT(x'',\xi'')\bigr)
B(x'',\xi'') \psi^0(x''\xi'')\, dx''d\xi''.
\label{19-4-89}
\end{multline}
Due to the microhyperbolicity and monotonicity of $\uptheta$ with respect to $\tau$ this convolution does not exceed
$C\mu^{r-1}h^{1-r} \times (\mu h)^{-r}\times (1+\mu^{2-l}|\log \mu|^{-\sigma})$
where the last factor is just upper bound of $B(x'',\xi'')$  and  $(\mu h)^{-r}$
comes from the trace $\Tr'$ and what we got does not exceed the right-hand expression of (\ref{19-4-84}).

We can rewrite the second term in (\ref{19-4-89}) as the convolution with respect to $\tau$  of $L\widehat{\bar{\chi}}(L\tau )$ (with $L=Th^{-1}$) with
\begin{equation}
(2\pi)^{-r} \mu^{r-1}h^{-r} \Tr' \int \uptheta 
\bigl(\tau - A_\cT(x'',\xi'')\bigr)\psi '(x''\xi'')\, dx''d\xi''
\label{19-4-90}
\end{equation}
and this convolution  does not exceed the same expression where factor
$(1+\mu^{2-l}|\log \mu|^{-\sigma})$ this time is an upper bound for $\psi'$.
\end{proof}

\subsection{Unmollified asymptotics. III}
\label{sect-19-4-5-6}

Now we can decompose expression
\begin{equation}
(2\pi )^{-r}\mu^rh^{-r} \int \Tr ' 
\uptheta  \bigl(\tau -A_\cT (x'',\xi'')\bigr)\psi_\cT (x'',\xi'') \,dx''d\xi''
\label{19-4-91}
\end{equation}
as
\begin{multline}
(2\pi )^{-r}\mu^rh^{-r} \int \Tr '
\uptheta  \bigl(\tau -\cA_0 (x'',\xi'')\bigr)
\psi ^0 (x'',\xi'')\,dx''d\xi''\ + \\
 h^{-d}\cN^\MW_\corr
 \label{19-4-92}
\end{multline}
with
\begin{multline}
\cN^\MW _\corr\Def\\
(2\pi )^{-r}(\mu h) ^r  \int \Tr ' \Bigl(
\uptheta  \bigl(\tau -A_\cT (x'',\xi'')\bigr)-
\uptheta  \bigl(\tau -\cA_0 (x'',\xi'')\bigr)\Bigr)\times \\
\shoveright{\psi ^0 (x'',\xi'') \,dx''d\xi''+}\\
(2\pi )^{-r}(\mu h)^{r} \int \Tr ' 
\uptheta  \bigl(\tau - A_\cT (x'',\xi'')\bigr)\times\qquad\qquad\qquad\qquad\qquad \\
\bigl(\psi_\cT(x'',\xi'')-\psi ^0 (x'',\xi'')\bigr) \,dx''d\xi''
\label{19-4-93}
\end{multline}
Then making change of variables $\Psi_0$ in the integration and using (\ref{19-2-18}) and (\ref{19-3-39}) for the Jacobian we conclude that
\begin{claim}\label{19-4-94}
The first term in (\ref{19-4-92}) equals to 
\begin{equation*}
h^{-d}\int \cN^{\MW} (x,\tau)\psi (x) \,dx.
\end{equation*}
\end{claim}

Let us  represent expression (\ref{19-4-93}) as
\begin{equation}
\cN^\MW_\corr = \int \cN^\MW_{\corr,x} (\tau)\psi (x) \,dx.
\label{19-4-95}
\end{equation}
While it is clearly the case with the first term in it, in the second term we can rewrite as
\begin{multline}
\psi_\cT(x'',\xi'')-\psi ^0 (x'',\xi'')=\\
\sum_{1\le |\alpha|+|\beta|\le m}
\psi_{\alpha\beta}(x'',\xi'') (x ')^\alpha (\mu^{-1}hD')^\beta+
O(\mu^{-l}|\log \mu|^{-\sigma})
\label{19-4-96}
\end{multline}
with $m=\lfloor (l,\sigma)\rfloor$; then we can just skip the last term (because the output of it will be $O(\mu^{-l}|\log \mu|^{-\sigma}h^{-d})$). Since \begin{equation*}
\psi_{\alpha\beta}\circ \Psi_0 = \sum_{\gamma: 1\le \gamma\le |\alpha|+\beta|}
\rho_{\alpha\beta\gamma}\partial ^\gamma\psi (x)
\end{equation*}
we can always rewrite the second term in the required form as well. What we are lacking is the estimate for these two terms.

\begin{proposition}\label{prop-19-4-15}
In the case of the intermediate magnetic field under the microhyperbolicity condition
\begin{equation}
\cN^\MW _\corr \le
C\left\{\begin{aligned}
&\mu^{-l}|\log \mu |^{-\sigma}\qquad 
&&\text{as\ \ }l<2;\\[2pt]
&\mu^{-2}(\mu h)^{l-2}  |\log (\mu h) |^{-\sigma}+\mu^{-1}h \qquad 
&&\text{as\ \ }l>2
\end{aligned}\right.
\end{equation}
(ii) In particular,  $\cN^\MW_\corr \le C\mu^{-1}h $ \ for \
$(l,\sigma)\succeq (3,0)$.
\end{proposition}

Therefore   $(2,0)\prec (l,\sigma)\prec (3,0)$ is the only case when we cannot skip $h^{-d}\cN^\MW_\corr$ in the final answer.

\begin{proof} 
(i) Assume first that $(l,\sigma)\preceq (2,0)$.

First of all, let us replace $A_\cT$ by $\cA$. Due to the microhyperbolicity and the fact that $(A_\cT-\cA)$ is bounded by $C\mu^{-l}|\log \mu|^{-\sigma}$, such operation causes an admissible error $O(\mu^{-l}|\log \mu|^{-\sigma}h^{-d})$. 

Replacing in the second term in (\ref{19-4-93}) $\cA$ by
$\cA_0$ we get an error $O(\mu^{-2})$. Really,  each factor
$\bigl(\uptheta(\tau-\cA)-\uptheta(\tau-\cA_0)\bigr)$ and $(\psi_\cT-\psi^0)$ contributes the factor $\mu^{-1}$ into the remainder estimate. So we need to consider the second term in (\ref{19-4-93}) with $A_\cT$ replaced by $\cA_0$ and with $(\psi_\cT-\psi^0)$ replaced by $\psi^1$, 
\begin{equation}
\psi^m=\sum_{|\alpha|+|\beta|=m}
\psi_{\alpha\beta}(x'',\xi'') (x ')^\alpha (\mu^{-1}hD')^\beta
\label{19-4-98}
\end{equation} 
but then after taking trace $\Tr'$ we get $0$.

In the first term in (\ref{19-4-93})  $\cA = \cA_0+\cA_1$. Consider $\cA_0$ as an unperturbed operator and $\cA_0+\cA_1$ as a perturbed one. Let us return back to representation (\ref{19-2-78}). 

We claim that
\begin{claim}\label{19-4-99}
As $m>l+1$, the contribution of $m$-th term of successive approximations does not exceed $Ch^{-d}\vartheta (\mu^{-1})$.
\end{claim}
Really, replace $\cA_0$ and $\cA_1$ by their $\eta$-mollifications. Then, in virtue of our previous arguments, replacing $\eta$-mollification of $\cA_0$ by its $2\eta$-mollification we get an error not exceeding
\begin{equation}
Ch^{-d}\vartheta (\eta) \times \bigl(\frac{\mu^{-1}T}{h}\bigr)^{m-1} \times
\bigl(1+ \frac{\eta T}{h}\bigr)^{-s} 
\label{19-4-100}
\end{equation}
and similarly replacing $\eta$-mollification of $\cA_1$ by its $2\eta$-mollification we get an error not exceeding the same expression with an  factor $\vartheta(\eta)$ replaced by $hT^{-1} \vartheta(\eta)\eta^{-1}$; thus we get the same expression (\ref{19-4-100}). 

Summation with respect to $\eta$ returns 
$Ch^{-d}\vartheta \bigl(h/T\bigr) \times \bigl(T/\mu h \bigr)^{m-1}$ and then summation with respect to $T\le C\mu h$ returns this expression as $T=\mu h$ i.e. $Ch^{-d}\vartheta (\mu^{-1})$. On the other hand, as $T\ge C\mu h$ we just do not need to consider successive approximations but the whole difference.

Further, as we have $\eta$-mollified $\cA_0$ and $\cA_1$, the contribution of this term does not exceed 
\begin{equation*}
Ch^{1-d}  T^{-1} \times \bigl(\frac{\mu^{-1}T}{h}\bigr)^{m-1} \times
\bigl(1+ \frac{\eta T}{h}\bigr)^{-s} 
\end{equation*}
and taking $\eta =\epsilon$ we get  $Ch^{-d}\mu^{1-m}$ after summation with respect to $T$. So (\ref{19-4-99}) is proven.

\medskip

Therefore for $l<2$ we can take $m=3$ i.e. take only two terms of the successive approximations. 

Note that the first term of the successive approximations annihilates with the corresponding expression for an unperturbed operator. 

Recall that
\begin{gather*}
\cA_0 =\Bigl(\mu^2 \sum_{j,k} a_{jk} \zeta^\dag_j \zeta_k +V\Bigr)^\w
\intertext{preserves spaces $\sH_n$ and} 
\cA_1 =\Bigl(2\Re \mu^2 \sum_{j,k, p} b_{jkp} \zeta^\dag_j \zeta_k\zeta^p +
\sum_k b_k\zeta_k\Bigr)^\w
\end{gather*}
maps $\sH_n$ into $\sH_{n+1}\oplus \sH_{n-1}$
and therefore after taking trace the second term of successive approximations results in 0 for $\psi_\cT$ replaced by $\psi^0$.

We need also to consider the similar expressions for $\psi_\cT$ replaced by
$\psi^1$  but in this case the first term of successive approximations vanishes and the second term has a desired estimate. 

\medskip\noindent
(ii) On the other hand, for $l>2$, $(l,\sigma)\preceq (3,0)$ we need to remember that $\cA = \cA_0+\cA_1+\cA_2$, $\psi_\cT=\psi^0+ \psi^1+\psi^2$ and also we need to consider four extra terms (the rest definitely is estimated in a desired way) which are\footnote{\label{foot-19-26} As we return back to parametrices, cut-off on interval $[-\hat{T},\hat{T}]$ etc.}
\begin{align}
&ih^{-1}\int_{-\infty}^0 F_{t\to h^{-1}\tau}\bar{\chi}_T(t)
\Tr \Bigl(\sum_{\pm} \mp \cG_0^\pm \cA_1\cG_0^\pm \cA_1\cG_0^\pm \psi^0\Bigr)\,d\tau,
\label{19-4-101}\\
&ih^{-1}\int_{-\infty}^0 F_{t\to h^{-1}\tau}\bar{\chi}_T(t)
\Tr\Bigl(\sum_{\pm} \mp\cG_0^\pm \cA_2\cG_0^\pm  \psi^0\Bigr)\,d\tau,
\label{19-4-102}\\
&ih^{-1}\int_{-\infty}^0 F_{t\to h^{-1}\tau}\bar{\chi}_T(t)
\Tr  \Bigl(\sum_{\pm} \mp\cG_0^\pm  \cA_1 \psi^1\Bigr)\, d\tau,
\label{19-4-103}\\
&ih^{-1}\int_{-\infty}^0 F_{t\to h^{-1}\tau}\bar{\chi}_T(t)
\Tr \Bigl(\sum_{\pm} \mp\cG_0^\pm \psi_2 \Bigr)\, d\tau
\label{19-4-104}
\end{align}
with $T=\hat{T} = \mu h$ and $\cG_0^\pm $ parametrices of $hD_t-\cA_0$. Here $\Tr $ means that we take $x=y$ and then integrate; recall that all operators here depend on $y$.

Really, even as $l=3$ we need to consider to $m=4$ but this term vanishes in virtue the above arguments related to $\sH_n$-scale.

Replacing $\bar{\chi}_T(t)$ by $\chi_T(t)$ and applying the same technique as before we find that such modified expressions (\ref{19-4-101})--(\ref{19-4-104}) do not exceed
$C\mu^{-2}h^{-d}(h/T)^{l-2}|\log (h/T)|^{-\sigma}$. 

Then after summation over $t$-partition we find that the contribution of segment $[\bar{T},\hat{T}]$ is estimated now by
$C\mu^{-2}h^{-d}(h/\bar{T})^{l-2}|\log (h/\bar{T})|^{-\sigma}$. We take
$\bar{T}= \epsilon \mu^{-1}$. This is clearly possible provided 
\begin{equation}
\mu \le h^{\delta -1}.
\label{19-4-105}
\end{equation}
Assuming this, let us consider contribution of $[-\bar{T},\bar{T}]$. Returning from $\cT$ reduction to the initial settings we can rewrite (\ref{19-4-101})--(\ref{19-4-104}) (modulo terms, estimated properly) as the similar expressions but with $\cG_0^\pm$ replaced by parametrices $\bar{G}^\pm$ of operator $\bar{A}$, given by (\ref{19-2-84}):
\begin{multline}
\bar{A}= \sum_{j,k} \bar{g}^{jk}  \bar{P}_j\bar{P}_k+\bar{V},\qquad
\bar{g}^{jk}=g^{jk}(y), \bar{V}=V(y), \\
\bar{P}_j=hD_j-V_j(y)-\sum_k(\partial_kV_j)(y)(x_k-y_k),
\label{19-4-106}
\end{multline}
and $\cA_j$ and $\psi^j$ are replaced by
\begin{gather}
\bar{A}_j =
\sum_{\alpha:|\alpha|\ge 2+j} \mu^{2-|\alpha|} b_\alpha (y) \bar{P}^\alpha +
\sum_{\alpha:|\alpha|\ge j} \mu^{-|\alpha|}b_\alpha (y) \bar{P}^\alpha,
\label{19-4-107}\\
\shortintertext{and}
\bar{\psi}^j=\sum_{\alpha:|\alpha|\ge j} \mu^{-|\alpha|}
\psi_\alpha (y)\bar{P}^\alpha
\label{19-4-108}
\end{gather}
respectively.

Note that the coefficients of all operators are smooth in
$\varepsilon$-scale with $\varepsilon=\mu^{-1}$ which is larger than
$Ch|\log h|$ required by the standard theory which is applicable then (after we scale $x\mapsto \mu x$, $h\mapsto \mu h$, $t\mapsto \mu t$ producing that

\begin{enumerate}[leftmargin=*,label=(\alph*)]
\item 
before rescaling the contributions of intervals 
$[-\epsilon \mu^{-1}, -Ch|\log h|]$
and  $[Ch|\log h|, \epsilon \mu^{-1}]$  are negligible and

\item 
the contribution of interval $[-Ch|\log h|, Ch|\log h|]$ is equal to
\newline $\mu^{-2} \sum_{j\ge 0} \kappa_j (\mu h)^j$
with the coefficients $\kappa_j$ which do not depend on $\mu, h$.
\end{enumerate}
We should not care about terms with $j\ge 1$ here because these terms do not exceed $C\mu^{-1}h^{1-d}$. Therefore,
\begin{equation*}
h^{-d}\cN^\MW_\corr=
\kappa_0\mu^{-2}h^{-d}+O\bigl(\mu^{-2}(\mu h)^{l-2}|\log (\mu h)|^{-\sigma}
+\mu^{-1}h\bigr)
\end{equation*}
(we can alway replace it modulo properly estimated).

Now notice that the ``main part'' of asymptotics is
$h^{-d}\int \cN^\MW (x)\psi(x)\,dx$ which in comparison with theorem \ref{thm-19-2-16} for $\mu \asymp (h|\log h|)^{-\frac{1}{2}}$ implies that $\kappa_0=0$.

Therefore, under assumption (\ref{19-4-105}) statement (ii) of proposition is proven.

\medskip\noindent
(iii) To get rid of assumption (\ref{19-4-105}) let us return to terms (\ref{19-4-101})--(\ref{19-4-104}) with $\bar{\chi}$ replaced by $\chi$ and with
$T\ge \epsilon \mu^{-1}$.

Each of these expressions is equal to $\mu^{-2}h^{-d}f(T/h, \mu h)$ and therefore the above arguments imply estimate
\begin{multline}
\mu h|\log h|\le \epsilon,\quad  T\ge h^{1-\delta},\quad  1\ge \varepsilon \ge C_0\mu^{-1}
\implies \\
|f(\frac{T}{h}, \mu h)|\le C(\frac{h}{T})^{l-2}|\log (\frac{h}{T})|^{-\sigma}.
\label{19-4-109}
\end{multline}
This is a really strange inequality because its assumptions involve $T$ and $h$
separately but its conclusion contains them only as $T/h$. This gives us a certain flexibility. Namely, let us assume that
\begin{equation}
h\le T \le h^{1-\delta}
\label{19-4-110}
\end{equation}
(otherwise everything is fine) and replace $h\mapsto h\lambda$,
$T\mapsto T\lambda$, $\mu \mapsto \mu\lambda^{-1}$ which affects neither $\mu h$ nor $(h/T)$. Assumptions of (\ref{19-4-109}) should be fulfilled now with these modified $h, T, \mu$. To fulfil the second one we take
$\lambda = h^{-1} (hT^{-1})^{1/\delta}$ and this is greater than $1$ due to (\ref{19-4-110}). 

The first of the assumptions of (\ref{19-4-109}) will survive\footnote{\label{foot-19-27} In fact, it will be fulfilled as long as $\mu h |\log (h/T)|\le \epsilon'$.} but the third one is now replaced by a more restrictive condition  $\varepsilon \ge C_0\mu^{-1}\lambda$. To get rid of it let us apply $C\mu^{-1}\lambda$-mollification which leads to the error in $f(T/h, \mu h)$,  not exceeding
\begin{equation*}
C(\mu ^{-1}\lambda )^{l-1}|\log \mu ^{-1}\lambda |^{-\sigma}  \bigl(\frac{T}{h}\bigr)^2 +
C(\mu ^{-1}\lambda )^l |\log \mu ^{-1}\lambda |^{-\sigma} \bigl(\frac{T}{h}\bigr)^3.
\end{equation*}
This is clearly less than the right-hand expression in (\ref{19-4-109}) provided
$T\ge h(\mu h)^{-\delta'}$ with arbitrarily small exponent $\delta'>0$
($\delta $ depends on it). It also implies $1\ge \varepsilon$.

Now inequality (\ref{19-4-109}) is proven under
this humble assumption. Then one can take $\bar{T}=\epsilon \mu^{-1}$ which clearly satisfies it.

The rest of arguments of (ii) hold without any modifications.
\end{proof}

\begin{remark}\label{rem-19-4-16} 
Case $l=2$ is obviously missing in the statement of proposition \ref{prop-19-4-15}.

\medskip\noindent
(i) Using the same arguments as above we find that contribution of zone
$\{ \bar{T}\le |t|\le \hat{T}\}$ does not exceed
$C\mu^{-2}h^{-d}|\log (\mu h)|^{1-\sigma}$ as  $\sigma >1$ and
$C\mu^{-2}h^{-d}|\log h|^{1-\sigma}$ as $\sigma <1$ and
$C\mu^{-2}h^{-d} |\log \bigl(|\log (\mu h)|/|\log h|\bigr) |$ as $\sigma=1$.

On the other hand, contribution of zone $\{|t|\le \bar{T}\}$ will be smaller as
$\sigma \le 1$ and
$\kappa_0\mu^{-2}h^{-d}+O\bigl(\mu^{-2}h^{-d}|\log (\mu h)|^{1-\sigma}\bigr)$
as $\sigma >1$ and the same arguments as above show that $\kappa_0=0$.

Therefore in comparison with the case $l\ne 2$ an extra factor $|\log \mu h|$ appears as $\sigma\ne 1$; as $\sigma=1$ it is 
$|\log \mu h|\cdot |\log \bigl(|\log (\mu h)|/|\log h|\bigr) |$.

\medskip\noindent
(ii) Further, if \underline{either} there are no third-order resonances  \underline{or}  $g^{jk}$ and $F_{jk}$ are constant, then $\cA_1=0$ after proper reduction. We can exploit it properly only for $\sigma\le 0$: one can see easily that then $\cN^\MW_\corr \le C\mu^{-2}|\log h|^{-\sigma}$.
\end{remark}

This concludes the proof theorem~\ref{thm-19-4-17} below.

\section{Main theorem}
\label{sect-19-4-6}

Now combining all the results of this Section we arrive to

\begin{theorem}\label{thm-19-4-17} 
Let assumptions \textup{(\ref{19-1-4})}--\textup{(\ref{19-1-6})}, $\textup{(\ref{19-1-25})}_{1-3}$ with $(\bar{l},\bar{\sigma})\succeq(l,\sigma)\succeq (1,2)$, $(\bar{l},\bar{\sigma})\succeq (2,1)$ and \textup{(\ref{19-2-38})}  be fulfilled. Let 
\begin{equation}
\mu^*_1 \Def C(h |\log h|)^{-\frac{1}{2}}\le \mu \le \mu^*_2\Def \epsilon ' (h |\log h|)^{-1}
\label{19-4-111}
\end{equation}
with sufficiently small constant $\epsilon'>0$.

Then there are two framing approximations\footref{book_new-foot-18-16} (see Chapter~\ref{book_new-sect-18} of \cite{futurebook}) such that 

\medskip\noindent
(i) Let  $\fN$-microhyperbolicity condition (see definition~\ref{def-19-2-5}) be fulfilled. Then if \underline{either}  $(l,\sigma)\succeq (2,0)$ \underline{or} $\#\fN=1$, then 
\begin{multline}
\R_1^\MW \Def\\
|\int \Bigl( \tilde{e}(x,x,0) - h^{-d}\cN^\MW (x,0)- 
h^{-d}\cN^\MW_{1\corr }(x,0)\Bigr)\psi (x)\, dx |\le\\[3pt]
C\mu ^{-1}h^{1-d} + C\mu ^{-l} |\log h|^{-\sigma} h^{-d}
\label{19-4-112}
\end{multline}
and if  $(l,\sigma)\prec (2,0)$ \underline{and} $\#\fN\ge 2$ then
\begin{equation}
\R_1^\MW \le
C\mu ^{l-l} h^{1-d}|\log \mu |^{- \sigma } + 
C\mu ^{-l} |\log h|^{-\sigma} h^{-d};
\label{19-4-113}
\end{equation}

\medskip\noindent
(ii) Let  microhyperbolicity condition (see definition~\ref{def-19-2-4}) be fulfilled. Then
\begin{equation}
\R_1^\MW \le
C h^{1-d}| +  C\mu ^{-l} |\log h|^{-\sigma} h^{-d}.
\label{19-4-114}
\end{equation}
\end{theorem}

\begin{remark}\label{rem-19-4-18} 
(i) Condition (\ref{19-4-111}) means exactly that estimate (\ref{19-4-112})  is stronger than (\ref{19-2-86}), (\ref{19-2-87}).

\noindent
(ii) If $(l,\sigma) \succeq (3,1)$ then the best remainder estimate (\ref{19-2-89}) holds as  $\mu \le \mu^*_2$.
\end{remark}

\section{Special case of constant $g^{jk}$, $F_{jk}$}
\label{sect-19-4-7}

\subsection{Framework; Mid- and long-range propagation}
\label{sect-19-4-7-1}

Consider now the case of constant $g^{jk}$ and $F_{jk}$, $(l,\sigma)\succeq (2,0)$ and $V$ having only non-degenerate critical points. Our goal is to extend results obtained for $V$ which does not have any critical points.

We assume that $l<4$ as $l>3$ will be sufficient for the remainder estimate $O(\mu^{-1}h)$. Note that the contribution of 
$\{x:\ |\nabla V|\le \bar{\nu}\Def C_1\mu^{-\frac{1}{2}}\}$ to the remainder does not exceed $C\bar{\nu}^d \mu h^{1-d}=O(\mu^{-1}h^{-1})$ and therefore we need to consider domain $\{x:\ |\nabla V|\ge \bar{\nu}\}$ only. 

Let us introduce in this domain
\begin{equation}
\varepsilon = C_0\mu^{-1}\nu^{-1}\le \nu,\qquad \nu \Def |\nabla V|.
\label{19-4-115}
\end{equation}
One can see easily that 
\begin{claim}\label{19-4-116}
Mollification error does not exceed $C\mu^{-l}|\log h|^{-\sigma}h^{-d}$
\end{claim}
exactly as in the case when $V$ does not have critical points.

Exactly as in Subsection~\ref{sect-19-4-1} one can prove easily that 

\begin{claim}\label{19-4-117}
At $B(y,\nu (y))$ with $\nu(y)\ge \bar{\nu}$ estimates (\ref{19-4-7}) and (\ref{19-4-8}) hold for $T\in [T_*,T^*]$ with $T^*= \epsilon \mu$ and 
\begin{equation}
T_*= C\varepsilon^{-1}\nu^{-1}h|\log h| \asymp C\mu h|\log h|.
\label{19-4-118}
\end{equation}
\end{claim}

\subsection{Short-range theory}
\label{sect-19-4-7-2}

Short-range theory as in Subsections~\ref{sect-19-4-2}--\ref{sect-19-4-4} is more complicated:

\begin{proposition}\label{prop-19-4-19} 
Let $\psi \in \sC_0^\infty (B(y,\nu(y))$ with $\nu (y)\ge \bar{\nu}$ and let $\phi$ be supported in $L$-vicinity  of $0$, $L\ge \nu^2+\mu h+ h/T$. 

Then for $T\in [Ch^{1-\delta}, T^*]$ with arbitrarily small exponent $\delta >0$ such that $T\ge \nu^{-2}h$
\begin{align}
&|\phi (hD_t) \chi_T(t) ( \Gamma U\psi )|\le \label{19-4-119}\\[2pt]
&\qquad\qquad  C L\nu^{d+1}  h^{-d} 
\bigl(\frac{h} {\nu T}\bigr)^{l-1} |\log \frac{h} {\nu T}|^{-\sigma}
\bigl(1+\frac{h} {\nu T\varepsilon} \bigr)^{-s},\notag\\[3pt]
&|F_{t\to h^{-1}\tau} \chi_T(t) ( \Gamma U\psi )|\le \label{19-4-120}\\
&\qquad\qquad  C L \nu^{d-1}  h^{1-d} 
\bigl(\frac{h} {\nu T}\bigr)^{l-1} |\log \frac{h} {\nu T}|^{-\sigma}
\bigl(1+\frac{h} {\nu T\varepsilon} \bigr)^{-s}\notag
\end{align}
as $\tau\le c$ with arbitrarily large exponent $s$.
\end{proposition}

\begin{proof}[Idea of the proof]
The proof follows the same scheme as proposition~\ref{prop-19-4-9}. Instead of giving the full proof, which we leave to the reader, let us consider the pilot-model for the proof proposition~\ref{prop-19-4-8}.

Consider again $I_1,I_2$ defined by (\ref{19-4-25}) and (\ref{19-4-26}) respectively. Replace $\lambda$ by its $\eta$-mollification with $\eta\le \epsilon \nu h^\delta $. Then an error in $I_1$ is  $O\bigl(h^{-1}T\vartheta(\eta)\nu^d\bigr)$ (just acquires factor $\nu^d$ as a volume of $B(y,\nu(y))$ in comparison with the original pilot model) while an error in $I_2$ is given by the left-hand expression of (\ref{19-4-27}) which does not exceed $Ch^{-1}T\vartheta(\eta)\nu^{d-2}$ now.

Obviously $I_1$ with mollified $\lambda$ does not exceed 
$C(h/\nu T)^k \eta^{-k-1}\vartheta(\eta)\nu^d$ and $I_2$ with mollified $\lambda$ does not exceed 
$C\nu^{-1} (h/\nu T)^k \eta^{-k-1}\vartheta(\eta)\nu^d$ and therefore
\begin{equation*}
I_1+ \nu^2 I_2\le Ch^{-1}T\Bigl(1 + \nu (h/\nu T)^{k+1} \eta^{-k-1}\Bigr)\vartheta(\eta)\nu^d.
\end{equation*}
Taking near optimal $\eta = h/T\nu ^{1-\delta'}$ we conclude that 
\begin{equation*}
I_1+ \nu^2 I_2\le Ch^{-1}T\vartheta(h/\nu T)\nu^d.
\end{equation*}
Note that $\eta \le h^{\delta}\nu$ as $T\ge \max(h^{1-\delta}, h/\nu^2)$.

This would imply for a single $\mu^{-1}h$-operator estimates as in the proof proposition~\ref{prop-19-4-9} but with $T$ replaced by $T\nu$ and with extra factors $\nu^{d+1}$ and $\nu^{d-1}$. However we also need to take into account that the ellipticity is broken for no more than $CL (\mu h)^{-r}$
multiindices $\alpha\in \bZ^{+\,r}$ rather than $(\mu h)^{-r}$ as it was for $\nu\asymp 1$ and this leads to an extra factor $CL$.

Easy details are left to the reader.
\end{proof}

Then summation of the right-hand expression of (\ref{19-4-120}) for 
$T\ge h/\nu^2$ returns $O(\nu^dh^{1-d})$. 

On the other hand note that  
\begin{equation*}
|F_{t\to h^{-1}\tau} \bar{\chi}_T(t) ( \Gamma U\psi )|\le \\
 C  T \nu^d  h^{-d} 
\end{equation*}
and plugging $T= h/\nu^2$ into the right-hand expression we get 
$C  \nu^{d-2}  h^{1-d}$. 

\begin{corollary} \label{cor-19-4-20} 
 In the framework of proposition~\ref{prop-19-4-19} as either 
$\mu \le h^{\delta-1}$ or $\nu \le h^\delta$

\medskip\noindent
(i) Estimate holds
\begin{equation}
|F_{t\to h^{-1}\tau} \bar{\chi}_T(t) ( \Gamma U\psi )|\le C \nu^{d-2}  h^{1-d} 
\label{19-4-121}
\end{equation}
for $T\le T^*=\epsilon \mu $;

\medskip\noindent
(ii) Therefore contribution of $B(y,\nu(y))$ to the Tauberian remainder estimate does not exceed $C\mu^{-1} \nu^{d-2}h^{1-d}$.
\end{corollary}

\subsection{Calculations}
\label{sect-19-4-7-3}
Following Subsection~\ref{sect-19-4-5} we arrive to expression (\ref{19-4-83}) where now 
$\cA=\cA_0 +\mu^{-2}\cA_2 + \cA'$ with $\cA'=O(\varepsilon^l |\log \varepsilon|^{-\sigma})$ and applying the method of the successive approximations we conclude first that modulo $o(h^{1-d})$ we can replace here $\cA$ by $\cA_0$ and then  (\ref{19-4-121}) holds without restriction 
``either $\mu \le h^{\delta-1}$ or $\nu \le h^\delta$''.

Summation with respect to a partition of unity results in the remainder estimate $O(\mu^{-1}h^{1-d})$ under assumption (\ref{19-2-26}). Therefore

\begin{claim}\label{19-4-122}
Under assumption (\ref{19-2-96}) Tauberian remainder is $O(\mu^{-1}h^{1-d})$.
\end{claim}

Further, let us consider expression (\ref{19-4-83}) integrated with respect to $\tau$. Let us skip  $\cA'$; one can see easily that the contribution of $B(y,\nu(y))$ to  the error does not exceed 
$C\bigl( h^{-d}+\mu h^{1-d} \nu^{-2}\bigr) \nu ^d \varepsilon^l |\log \varepsilon|^{-\sigma}$ and summation with respect to partition of unity results in $O(h^{-d}\varepsilon^l |\log \varepsilon|^{-\sigma})$. 

Furthermore, skipping term $\mu^{-2}\cA_2$ results in the correction term which after replacing a Riemann sum by an integral as $\nu^2\ge \mu h$ results in 
the expression $\mu^{-2}h^{-d} \int \kappa_0(y)\,dy$; contribution of $B(y,\nu(y))$  into an error does not exceed 
\begin{equation*}
C\mu^{-2}h^{-d} (\mu h /\nu^2)^{l-2} |\log (\mu h/\nu^2)|^{-\sigma}\nu^d 
\end{equation*}
and summation with respect to partition of unity results in 
\begin{equation}
C\mu^{-2}h^{-d} (\mu h )^{l-2} |\log (\mu h)|^{-\sigma}.
\label{19-4-123}
\end{equation}

Meanwhile skipping term $\mu^{-2}\cA_2$ as $\nu ^2\le \mu h$ results in an error; one can prove easily that the contribution of $B(y,\nu(y))$  into this error does not exceed  $C\mu^{-1}h^{1-d} \nu^{d-2}$. Treating all other terms (arising when we replace $\psi$ by $\psi^0$) in the similar way we arrive to a correction term not exceeding (\ref{19-4-123}) as ``final'' $\kappa_0$ must be $0$.

\subsection{Main theorem}
\label{sect-19-4-7-4}

Then we arrive to

\begin{theorem}\label{thm-19-4-21}
Let assumptions \textup{(\ref{19-1-4})}--\textup{(\ref{19-1-6})}, $\textup{(\ref{19-1-25})}_3$ with $(l,\sigma)\succeq (1,2)$ and \textup{(\ref{19-2-38})}  be fulfilled. Let $g^{jk}$, $F_{jk}$ be constant. 

Then there are two framing approximations\footref{book_new-foot-18-16} (see Chapter~\ref{book_new-sect-18} of \cite{futurebook}) such that the following statements are true:

\medskip \noindent
(i) Under microhyperbolicity assumption  $|\nabla V|\ge \epsilon$ estimate  \textup{(\ref{19-4-112})} holds;

\medskip\noindent
(ii) Under assumptions $(l,\sigma)\succeq (2,0)$ and \textup{(\ref{19-2-96})} estimate  \textup{(\ref{19-4-112})} holds;

\medskip\noindent
(iii) We can skip correction term unless $(2,0)\prec (l,\sigma)\prec (3,0)$ in which case it does not exceed \textup{(\ref{19-4-123})}.

\medskip\noindent
(iv) We have remainder estimate $O(\mu^{-1}h^{1-d})$ as 
$\mu \ge (h|\log h|)^{-\frac{1}{2}}$ provided $(l,\sigma)\succeq (3,1)$.
\end{theorem}

\chapter{Strong magnetic field}
\label{sect-19-5}

\section{Framework and special cases}
\label{sect-19-5-1}
Let now 
\begin{equation}
\mu^*_2\Def \epsilon (h|\log h|)^{-1}\le \mu \le \mu^*_3\Def \epsilon h^{-1}.
\label{19-5-1}
\end{equation}

Then we need to modify our arguments in several way. First of all now, as 
$(\mu^{-1}h|\log h|)^{\frac{1}{2}}\ge \mu^{-1}$ we must assume that
\begin{equation}
\varepsilon \ge C(\mu^{-1}h|\log h|)^{\frac{1}{2}}.
\label{19-5-2}
\end{equation}
Then under $\fN$-microhyperbolicity condition with $\#\fN=1$ we are done and we have immediately the first clause (\underline{either}) of our final theorem~\ref{thm-19-5-2} below which differs from theorem~\ref{thm-19-4-17}(i) only by a choice of $\varepsilon$ and thus by a mollification error. 

Furthermore, in the case of the constant $g^{jk}$, $F_{jk}$ we are done as well and we have immediately our final theorem~\ref{thm-19-5-3} which differs from theorem~\ref{thm-19-4-21} in the same way. 

However under either microhyperbolicity condition with $\#\fM \ge 2$ or $\fN$-microhyperbolicity condition with $\#\fN\ge 2$ we were dealing with the partition of energies and then variables $(x',\xi')$ are no more microlocal as logarithmic uncertainty principle for them fails and our approach here needs to be modified.

Still, we are also done if we \underline{either} use microhyperbolicity assumption and there are no $2$-nd order resonances \underline{or} we use $\fN$-microhyperbolicity assumption and there are no $2$-nd and $3$-rd order resonances. 

Really, in these cases we need to consider operators 
\begin{equation}
a_j= f_j(x'',\mu^{-1}hD'')(h^2D_j^2 + \mu^2 x_j^2)
\label{19-5-3}
\end{equation}
with $\fm=\{j\}$ and $\fn=\{j\}$ respectively and instead we can consider commuting operators 
$(h^2D_j^2 + \mu^2 x_j^2)$ and we do not need to have logarithmic uncertainty principle for these variables to establish the following proposition similar to proposition~\ref{prop-19-2-10}(i),(ii) with small constant $\nu T$:

\begin{proposition}\label{prop-19-5-1}
Let \underline{either} there are no $2$-nd order resonances and $T^*=\epsilon$ 
\underline{or} there are no $2$-nd and $3$-rd order resonances and 
$T^*=\epsilon \mu$. Let $\mu \le \epsilon_0 h^{-1}$.

Let $B_j= (h^2D_j^2 + \mu^2 x_j^2)$ and $\mathbf{B}=(B_1,\ldots, B_d)$. Let $\phi_1,\phi_2$ be two functions with 
\begin{gather}
\dist (\supp(\phi_1),\supp(\phi_2))\ge (T/T^* +\epsilon_1)
\label{19-5-4}\\
\intertext{the following estimate holds:}
\|\bigl(\phi_2 (\mathbf{B})\bigr)_x U\,
^t\! \bigl(\phi _1 (\mathbf{B} )\bigr)_y \| \le Ch^s. 
\label{19-5-5}
\end{gather}
\end{proposition}
Then we have immediately clause \underline{or} of our final theorem~\ref{thm-19-5-2}.

Thus we arrive to our first main theorem in the case of a strong magnetic field:

\begin{theorem}\label{thm-19-5-2} 
Let assumptions \textup{(\ref{19-1-4})}--\textup{(\ref{19-1-6})}, $\textup{(\ref{19-1-25})}_{1-3}$ with \linebreak $(\bar{l},\bar{\sigma})\succeq(l,\sigma)\succeq (1,2)$, $(\bar{l},\bar{\sigma})\succeq (2,1)$ and \textup{(\ref{19-2-38})}  be fulfilled.  Let $\mu$ satisfy \textup{(\ref{19-5-1})}. 

Then there are two framing approximations\footref{book_new-foot-18-16} (see Chapter~\ref{book_new-sect-18} of \cite{futurebook}) such that 

\medskip\noindent
(i) Let   $\fN$-microhyperbolicity condition (see definition~\ref{def-19-2-5}) be fulfilled and \underline{either} $\#\fN=1$ \underline{or} there be no $2$-nd and $3$-rd order resonances. Then
\begin{multline}
\R_1^\MW \Def\\
|\int \Bigl( \tilde{e}(x,x,0) - h^{-d}\cN^\MW (x,0)- 
h^{-d}\cN^\MW_{1\corr }(x,0)\Bigr)\psi (x)\, dx |\le\\[3pt]
C\mu ^{-1}h^{1-d} +
C\mu ^{-\frac{l}{2}} |\log h|^{\frac{l}{2}-\sigma} h^{-d+\frac{l}{2}}
\label{19-5-6}
\end{multline}

\medskip\noindent
(ii) Let  microhyperbolicity condition (see definition~\ref{def-19-2-4}) be fulfilled and there are be no $2$-nd order resonances. Then
\begin{equation}
\R_1^\MW \le
C h^{1-d} +  C\mu ^{-\frac{l}{2}} |\log h|^{\frac{l}{2}-\sigma} h^{-d+\frac{l}{2}}.
\label{19-5-7}
\end{equation}
\end{theorem}

\begin{theorem}\label{thm-19-5-3}
Let assumptions \textup{(\ref{19-1-4})}--\textup{(\ref{19-1-6})}, $\textup{(\ref{19-1-25})}_{3}$ with $(l,\sigma)\succeq (1,2)$ and \textup{(\ref{19-2-38})}  be fulfilled.  Let $g^{jk}$, $F_{jk}$ be constant. 

Then there are two framing approximations\footref{book_new-foot-18-16} (see Chapter~\ref{book_new-sect-18} of \cite{futurebook}) such that the following statements are true:

\medskip \noindent
(i) Under microhyperbolicity assumption  $|\nabla V|\ge \epsilon$ estimate  \textup{(\ref{19-5-6})} holds;

\medskip\noindent
(ii) Under assumptions $(l,\sigma)\succeq (2,0)$ and \textup{(\ref{19-2-96})} estimate  \textup{(\ref{19-5-6})} holds.
\end{theorem}

\section{Partition of energy space} 
\label{sect-19-5-2}

Consider the general case now. To  prove a statement similar to proposition \ref{prop-19-4-5}  under either microhyperbolicity condition with $\#\fM\ge 2$ or $\#\fN$-microhyperbolicity condition  with $\#\fN\ge 2$  one needs to consider evolution of $a_{\fm}$ or $a_{\fn}$ and the main trouble is that we can microlocalize them only in the box of size $\mu ^2\varepsilon ^2 \asymp \mu h  |\log h |$ which is by no means smaller than $\epsilon$. We need instead \emph{operator localization\/}\index{operator localization}: something like proposition~\ref{prop-19-5-1} but with $\mathbf{B}=(B_1,\ldots,B_m)$ with $B_\fm$ operators ``close'' to $a^0_\fm-\tau_\fm$, $m=\#\fM$\,\footnote{\label{foot-19-28} Or similar estimate with $\fm$ and $\fM$ replaced by $\fn$ and $\fN$.}. 

Therefore we need to define function of operators and thus operators need to be self-adjoint and ``almost commute''.

Let  $B$ be an operator of the form 
\begin{equation}
B= \sum _{i,k\in \fm} \bigl(
\beta^{\w}_{j,k}(x'',\mu^{-1}hD'') \, Z_j Z_k^*,\qquad 
Z_j= hD_j+i\mu x_j
\label{19-5-8}
\end{equation}
where $\beta_{jk}\in \sF^{2,0}$, $\beta_{jk}=\beta_{kj}^\dag$  and the complex sesquilinear  form $\sum_{j,k} \beta_{jk} \zeta_j\zeta_k^\dag$ is positive definite; here $\fm$ is just a subset in $\{1,\ldots,d\}$. 

This will be sufficient to study propagation if \underline{either} $T^*=\epsilon$ \underline{or} $T^*=\epsilon \mu$ and there are no $3$-rd order resonances; however if $T^*=\epsilon \mu$ and there are $3$-rd order resonances we consider 
\begin{multline}
B= \sum _{j,k\in \fn} \bigl(
\beta^{\w}_{jk}(x'',\mu^{-1}hD'')\, Z_jZ_k +\\
\mu^{-1} \sum _{j,k,m\in \fn} \Bigl(\beta^{\w}_{jkm}(x'',\mu^{-1}hD'')\,Z_jZ_kZ_m^*+
\beta^{\dag\, \w}_{jkm}(x'',\mu^{-1}hD'')\,Z_m Z_k^*Z_j^*\Bigr) +\\
C_0 \mu^{-2} (\sum _{j\in \fn}  Z_j^*Z_j)^2
\tag*{$\textup{(\ref{19-5-8})}'$}\label{19-5-8-'}
\end{multline}
where $\beta_{jkm}\in \sF^{1,0}$ and the last term is added to make the symbol positive.

Then this operator is self-adjoint and 
\begin{equation}
\phi (B)= \int e^{it  B}\hat{\phi} (t ) \,dt  
\label{19-5-9}
\end{equation}
where $\hat{\phi}$ is a Fourier transform.

Let temporarily $A=\alpha^\w(x'',\mu^{-1}hD'')$ or 
\begin{equation}
A=\sum_j \Bigl(\alpha^\w(x'',\mu^{-1}hD'')Z_j+ 
\alpha^{\dag\w}(x'',\mu^{-1}hD'') Z_j^*\Bigr)
\label{19-5-10}
\end{equation}
with $\alpha\in \sF^{1,0}$. Let us consider $[\phi(B),A]$. Note first that 
\begin{multline*}
e^{it B}Ae^{-it B} = \sum_{0\le k\le m-1} \frac{1}{k!} (it )^k
\Ad_B^k (A) + \\
\frac{1}{(m-1)!} (it)^m \int_0^1 (1-z)^{m-1} e^{itz B} \Ad_B^m (A) e^{-itz B}\,dz
\end{multline*}
where
\begin{equation}
\Ad_B^0 A\Def A,\quad
\Ad_B^m A\Def \bigl[B,\Ad_B^{m-1}A\bigr].
\label{19-5-11}
\end{equation} 
and therefore
\begin{multline*}
[e^{i\tau B},A]= \sum_{1\le k\le m-1} \frac{1}{k!} (it)^k
\Ad_B^k (A) e^{it B} +\\
\frac{1}{(m-1)!}(it)^m 
\int_0^1 (1-z)^{m-1}e^{itz B} \Ad_B^m (A) e^{it(1-z)B}\,dz
\end{multline*}
and 
\begin{multline}
[\phi(B),A]= \sum_{1\le k\le m-1} \frac{1}{k!} 
\Ad_B^k (A) \phi^{(k)}(B) +\\
\frac{1}{m!}\int_{-\infty}^\infty  \int_0^1   e^{itz B} \Ad_B^m (A) e^{it(1-z )B}\widehat{\phi^{(m)}}(t) \,dz dt.
\label{19-5-12}
\end{multline}
Consider repeated commutators $\Ad_B^k (A)$. Note first that if $B$ is of the form (\ref{19-5-8}), $A$ is defined by (\ref{19-5-10}) and if in the commutators  we consider only derivatives with respect to $(x',\xi')$ then obviously these shortened commutators would be also of the form (\ref{19-5-10}) but with norms 
$C_1 (C_0\mu h)^k$. 

On the other hand, in virtue of proposition~\ref{book_new-prop-1-A-1} of cite{futurebook} we can find uniformly smooth  function $\phi =\phi_n$, $\phi_n\ge \frac{1}{2}$ on the given interval and such that  $|\phi^{(k)}|\le C_0^{k+1} n^k $ as $k\le m+2\le n$. Then terms in the sum in the right-hand expression in (\ref{19-5-12}) would not exceed  
$C_1(C_0  \mu h)^k n^k k^{-k} $ which is $O(h^s)$ as $\mu h\le \epsilon$ and 
\begin{equation}
k\ge C_s|\log h|/|\log \mu h|, \qquad n\le \epsilon_0 k/(\mu h).
\label{19-5-13}
\end{equation}
Let us now consider complete commutators. Then as $\beta_{jk}\in \sF^{2,0}$, $\alpha_j\in \sF^{1,0}$ we conclude that in comparison with $\Ad_B^{m-1}(A)$  each term in  $\Ad_B^m(A)$ acquires either factor $C_0\mu h $ (if there are only differentiations with respect to $(x',\xi')$ or $C_1 \mu^{-1}h\varepsilon^{-1}N$ where  $N= \mu h^{-1}\varepsilon^2 $; it happens as long as $k\le N$, which is obviously the case; recall that $\varepsilon \ge C(\mu^{-1}h|\log h|)^{\frac{1}{2}}$. 

Further, similarly one can estimate an integral term in (\ref{19-5-12}). Furthermore, similar arguments work for operator $B$ defined by \ref{19-5-8-'}.

Therefore we arrive to
\begin{proposition}\label{prop-19-5-4}
As $\varepsilon$ is defined by \textup{(\ref{19-5-3})} (or larger)
\begin{multline}
[\phi(B),A]\equiv 
\sum_{1\le k\le M-1} \frac{1}{k!} \Ad_B^k (A) \phi^{(k)}(B), \\
\text{with\ \ } M= C_1|\log h|/|\log (\mu h)|.
\label{19-5-14}
\end{multline}
\end{proposition}

Similarly one can see easily that 
\begin{equation*}
e^{itB_1}e^{it'B_2}e^{-itB_1}e^{-it'B_2}\sim I+ 
\sum_{p,q}  \frac{1}{p!q!} (it)^p(it')^q L_{pq}
\end{equation*}
with
\begin{multline}
L_{pq}= \sum_{1\le j\le p}\sum_{1\le k\le q}
\frac{p!q!}{j!k!(p-j)!(q-k)!}  B_1^j B_2^k (-B_1)^{p-j}(-B_2)^{q-k}=\\
\left\{\begin{aligned}
&\sum_{1\le k\le q} \frac{q!}{k!(q-k)!}  \Ad^p_{B_1} (B_2^k)(-B_2)^{q-k}= \\
&\sum_{1\le j\le p}\frac{p!}{j!(p-j)!}  B_1^{p-j} \Ad^q_{B_2} ((-B_1)^j);
\end{aligned}\right.
\label{19-5-15}
\end{multline}
we will use the first (the second) expression as $q\le p$ ($p\le q$) respectively.

Then one can prove easily

\begin{proposition}\label{prop-19-5-5}
As $\varepsilon$ is defined by \textup{(\ref{19-5-3})} (or larger)
\begin{multline}
[\phi_1(B_1),\phi_2(B_2)]\equiv  
\sum_{1\le q, p\le M-1} \frac{1}{p!q!} L_{pq} \phi_2^{(q)}(B_2)\phi_1^{(p)}(B_1),\\
\text{with\ \ } M= C_1|\log h|/|\log (\mu h)|.
\label{19-5-16}
\end{multline}
\end{proposition}

\begin{proof}
We leave an easy proof using the same arguments as the proof proposition~\ref{prop-19-5-4} to the reader.
\end{proof}

\begin{corollary}\label{cor-19-5-6}
If $\dist(\supp \phi_1,\supp\phi_k)\ge \epsilon_0$, $i_1=i_k$ and 
$\mu h\le \epsilon $ then
\begin{equation}
\phi_1(B_{i_1}) \phi_2(B_{i_2})\cdots \phi_k(B_{i_k}) \equiv 0.
\label{19-5-17}
\end{equation}
\end{corollary}

\begin{proof}
If $k=3$ multiplying (\ref{19-5-16}) with $B_1$, $B_2$ replaced by $B_{i_1}$, $B_{i_2}$ respectively by $\phi _k(B_{i_1})$ on its right and using 
$\phi_1^{(p)}(B_{i_1})\phi_{i_3}(B_{i_3})=0$ we arrive to the required equality.

In the general case let $\varphi$ be functions such that $\varphi =1$ in vicinity of $\supp \phi_k$ and $\varphi =0$ in vicinity of $\supp \phi_1$. Then in virtue of (\ref{19-5-17}) with $k=3$ we can insert $\varphi (B_{i_1})$ between $\phi_{k-2}(B_{i_{k-2}})$ and $\phi_{k-1}(B_{i_{k-1}})$ and apply induction with respect to $k$.
\end{proof}

This is very important: we can localize with respect to $B_1,\ldots ,B_k$ \emph{simultaneously\/} by $\phi_1(B_1) \phi_2(B_2)\cdots \phi_k(B_k)$.

\section{Propagation}
\label{sect-19-5-3}

Now we can apply the same arguments as in Theorem~\ref{book_new-thm-2-3-1} of \cite{futurebook}. Consider operator 
\begin{equation}
\upchi \bigl(\epsilon ^{-1}(B-s) +  \frac{t}{T}\bigr)
\label{19-5-18}
\end{equation}
and apply it to $U^{\pm}=\uptheta(t) U$. Then
\begin{equation}
(hD_t-A)U^\pm = \mp ih \updelta (t) \updelta (x-y)
\label{19-5-19}
\end{equation}
and applying $\mu^{-1}h$-pseudo-differential operator cutoff we get
\begin{multline}
(hD_t-A)W^\pm   = \mp ih f(x,y)\qquad\text{with}\\[3pt]
W=U \,^t\!\bar{\varphi}_y (B)\,^t\!Q_y,\quad
f= \updelta (t) \updelta (x-y) \,^t\!\bar{\varphi}_y (B)\,^t\!Q_y,
\label{19-5-20}
\end{multline}
where $\bar{\varphi}$ is supported in vicinity of  $\tau \in\bR $ while $\upchi$ is our standard function and then
\begin{equation}
|T \Re ih^{-1}\Bigl(\Bigl[(hD_t - A),
\upchi \bigl(\epsilon ^{-1}(B-\tau) +\frac{t}{T}\bigr) \Bigr]W^\pm , W^\pm\Bigr)|\le Ch^{2s}.
\label{19-5-21}
\end{equation}

Transforming this inequality we get
\begin{multline}
|\Bigl(\upchi' \bigl(\epsilon ^{-1}(B-\tau) + \frac{t}{T} \bigr) W^\pm , W^\pm \Bigr)|\le \\
T\epsilon^{-1} |\Bigl(\Bigl[A, \upchi \bigl(\epsilon ^{-1}(B-\tau ) + 
C_0 \frac{t}{T}\bigr) \Bigr] W^\pm , W^\pm\Bigr) | + Ch^{2s}
\label{19-5-22}
\end{multline}
where $\upchi '$ is derivative of $\upchi$.

Approximating commutator by Poisson brackets (multiplied by $(-ih)$) we get 
\begin{equation*}
\epsilon^{-1} T \{a,b\} \upchi' \bigl(\epsilon ^{-1}(B-\tau) +  \frac{t}{T}\bigr)
\end{equation*}
with the first factor being an operator with a norm not exceeding $\frac{1}{2}$ provided 
\begin{multline}
T\le T^*\Def \\
\left\{\begin{aligned} 
&\epsilon_1 \qquad &&\text{as\ \ } b=a_\fm,\\
&\epsilon_1\mu  \qquad &&\text{as\ \ } b=a_\fn+O(\mu^{-2}),\quad (l,\sigma)\succeq (2,0),\\
&\epsilon_1\mu^{l-1}|\log \mu|^{-\sigma} &&\text{as\ \ } b=a_\fn+O(\mu^{-2}),\quad (l,\sigma)\prec (2,0)
\end{aligned}\right.
\label{19-5-23}
\end{multline}
where $\epsilon_1$ is small enough constant.

Therefore, taking $\upchi$ to be a primitive of function $\varphi ^2$, considering next terms of the commutator and assuming that 
\begin{equation}
\|\varphi\bigl(\epsilon ^{-1}(B-\tau) +  \frac{t}{T}-\epsilon_2\bigr)W^\pm \| 
\le Ch^{s-\delta}
\label{19-5-24}
\end{equation}
with a small exponent $\delta>0$, we get
\begin{multline}
\|\varphi \bigl(\epsilon ^{-1}(B-\tau ) + \frac{t}{T}\bigr )W^\pm \|^2 \le \\
 C_1 \biggl( \sum _{1\le m \le M}  \frac{1}{m!}  \eta ^m
\| \varphi ^{(m)} \bigl(\epsilon ^{-1}(B-\tau) +  \frac{t}{T}\bigr)W^\pm  \|\biggr)^2 + Ch^{2s}
\label{19-5-25}
\end{multline}
with $\eta=C_0\mu h$.

Plugging  $\frac{1}{m!} \eta ^m \varphi ^{(m)}$ instead of $\varphi $ and
$(M-m)$ instead of $M$ into  inequality (\ref{19-5-25}) we conclude that
\begin{equation}
\Bigl(\frac{1}{m!} \eta ^m
\| \varphi ^{(m)}\bigl(\epsilon ^{-1}(B-\tau) +\frac{t}{T}\bigr)W^\pm \|\Bigr)^2
\label{19-5-26}
\end{equation}
also does not exceed the right-hand expression of (\ref{19-5-25}).

Taking a sum with respect to $m$, $0\le m\le M$, and again increasing $C_0,C_1$ we get that expression
\begin{equation}
\sum_{0\le m\le M} \Bigl({\frac{1}{m!}} \eta ^m\|
\varphi ^{(m)}\bigl(\epsilon ^{-1}(B-\tau ) +
\frac{t}{T}\bigr)W^\pm \|\Bigr)^2
\label{19-5-27}
\end{equation}
does not exceed itself multiplied by $C_1 \eta $ plus $Ch^{2s}$. Taking $\eta$ small enough we conclude that (\ref{19-5-24}) holds with $\epsilon_2=0$ and $\delta=0$. Consequently increasing $s$ by $\delta$ and doubling $\epsilon_2$ we get rid off assumption (\ref{19-5-24}). Thus
\begin{equation}
\|\varphi  \bigl(\epsilon ^{-1}(B-\tau) + C_0 \mu ^{-1} t\bigr)
U \,^t\!\bar{\varphi}_y  (B) \,^t\!Q_yg \|
\le Ch^s.
\label{19-5-28}
\end{equation}

Further, we can plug $-(B-\tau)$ instead of $(B-\tau)$ into (\ref{19-5-28}). Therefore we have proven 

\begin{proposition}\label{prop-19-5-7} 
Let $\bar{\varphi}$ be supported in $\epsilon$-vicinity of $\tau $ and $\varphi =1$  in $2\epsilon$-vicinity of $\tau$.  Then for $T\le T^*$ defined by \textup{(\ref{19-5-23})}
\begin{equation}
\|\bigl(I-\varphi (B)\bigr) U \,^t\!\bar{\varphi}_y(B) \,^t\!Q_y \|
\le Ch^s.
\label{19-5-29}
\end{equation}
\end{proposition}

Using corollary~\ref{cor-19-5-6} which allows us to apply simultaneous cut-off by
\begin{equation}
\Phi (B_1,\dots,B_\nu)\Def
\varphi _1(B_1)\cdot \varphi _2(B_2)\cdots \varphi _\nu(B_\nu)
\label{19-5-30}
\end{equation}
we arrive to

\begin{corollary}\label{cor-19-5-8} 
Let functions $\bar{\varphi}_j$ be supported in $\epsilon$-vicinities of $\tau _j$ and let $\varphi_j =1$  in $2\epsilon$-vicinities of $\tau_j$. 
Then for $T\le T^*$ defined by \textup{(\ref{19-5-23})}
\begin{equation}
\|\bigl(1-\Phi(B_1,\dots,B_\nu)\bigr ) U
\,^t\!\bar{\Phi}_y(B_1,\dots, B_\nu) \,^t\!Q_y \| \le Ch^s.
\label{19-5-31}
\end{equation}
\end{corollary}

Recall that we consider  set of $B_j$ coinciding either with 
$\{a_\fm\}_{\fm \in \fM}$ or with $\{a_\fn+O(\mu^{-2})\}_{\fn \in \fN}$.

\section{Propagation  and microhyperbolicity}
\label{sect-19-5-4}

Consider now the proof that singularities leave diagonal -- again as theorem~\ref{book_new-thm-2-3-1} of \cite{futurebook} with 
\begin{equation}
\phi = \langle \ell , x-y\rangle \mp \epsilon tT^{-1}.
\label{19-5-32}
\end{equation}
We plug-in 
\begin{equation}
W =U  \,^t\!\bigl(\Phi (B_1,\dots,B_\nu )\bigr)_y \,^t\!Q_y
\label{19-5-33}
\end{equation}
localizing $x$ near $y$ and $B_j$ near $\tau_j$ for all $j$, and take
$\ell = \ell (y,\boldtau)$. 

To recover the proof we  need to estimate from below quadratic form 
\begin{equation*}
(\mu ^2 \sum _{j,k} p_j^\w K_{jk} p_k^\w w, w)
\end{equation*}
with $w = \upchi(\phi )^\w W$, $W=\bar{\Phi}(B_1,\ldots,B_\nu)U$ and $K_{jk} = \ell (g^{jk})$; one can assume that $V=-1$.

Due to proposition \ref{cor-19-5-8} one can replace $W$ by $\Phi(B_1,\ldots,B_\nu)W $ with a bit more wide supports of $\varphi_j$ than those of $\bar{\varphi}_j$.

Let us introduce $\mu ^{-1}$-admissible partition $\psi_\nu (x)$ and on each element of it apply gauge transformation making $|V_j|\le C\mu^{-1}$. Rescaling
$x\mapsto \mu (x-{\bar x}_\nu)$ we get $(\mu h)$-pseudo-differential operators $\mu p_j^\w$ and
$\mu ^2 \sum _{j,k} p_j^\w K_{jk} p_k^\w$ while $\Phi(B)$ is a legitimate
$(\mu h)$-pseudo-differential operator as well.

Then  for $\mu h\le \epsilon_1$
\begin{equation*}
\Phi(B_1,\ldots,B_\nu )^*
\Bigr(\mu ^2 \sum _{j,k} p_j^\w K_{jk} p_k^\w - \epsilon_0 +C \mu h\Bigl) \Phi(B_1,\dots,B_\nu )
\end{equation*}
is non-negative operator and
\begin{multline*}
\Phi(B_1,\ldots,B_\nu )^*  \mu ^2 \sum _{j,k} p_j^\w K_{jk} p_k^\w \Phi(B_1,\ldots,B_\nu ) \ge \\
\frac{1}{2} \epsilon_0 \Phi(B_1,\ldots,B_\nu )^* \Phi(B_1,\ldots,B_\nu )
\end{multline*}
in the operator sense and therefore
\begin{equation}
\Bigl(\mu ^2 \sum _{j,k} p_j^\w K_{jk} p_k^\w w, w \Bigr) \ge
\frac{1}{2} \epsilon _0 \|w\|^2 -Ch^s;
\label{19-5-34}
\end{equation}
the rest of the proof needs no modifications.

Now applying the same approach as in the proof proposition \ref{prop-19-2-13} we arrive to

\begin{proposition}\label{prop-19-5-9}
Let  \underline{either} $\fM$-microhyperbolicity assumption \underline{or} $\fN$-microhyperbolicity assumption be fulfilled.

Then estimate \textup{(\ref{19-2-63})} holds for $T\in [T_*,T^*]$ with
$T^*$  defined by \textup{(\ref{19-5-23})} and  
$T_* =C\varepsilon ^{-1}h |\log h| \le
\epsilon _1 ( \mu h |\log h|)^{\frac{1}{2}}$.
\end{proposition}

\begin{proof}
An easy proof we leave to the reader.
\end{proof}

\section{Short-range estimates}
\label{sect-19-5-5}

In this strong magnetic field case we can apply the same arguments as in Subsection~\ref{sect-19-4-4} combined with corollary~\ref{cor-19-5-8} and extend proposition \ref{prop-19-4-9}:

\begin{proposition}\label{prop-19-5-10} 
Both statements (i), (ii) of proposition \ref{prop-19-4-9} remain true in the case of the strong magnetic field as well.
\end{proposition}

\begin{proof} 
Just repeating proof proposition \ref{prop-19-4-9} without any significant difference. We use that in the operator sense
$|x'|+|\xi'|\le C\mu^{-1}$ still even if in the microlocal sense one should put $C(\mu ^{-1}h|\log h|)^{\frac{1}{2}}$ instead.

Easy details are left to the reader.
\end{proof}

\section{Calculations}
\label{sect-19-5-6}

In the strong magnetic field case one can apply the same arguments in the successive approximation method remembering that now
$\varepsilon \ge C(\mu^{-1}h|\log h|)^{\frac{1}{2}}$ and
$\bar{T}_0\le h\varepsilon^{-1}\le (\mu h |\log h|^{-1})^{\frac{1}{2}}$.

Then propositions \ref{prop-19-4-12}--\ref{prop-19-4-14} and their proofs remain true with this minor modification:

\begin{proposition}\label{prop-19-5-11} 
Under microhyperbolicity condition for  $\mu \le \epsilon h^{-1}$

\medskip\noindent
(i)  Estimates \textup{(\ref{19-4-77})} and \textup{(\ref{19-4-78})} hold;

\medskip\noindent
(ii)  Estimates \textup{(\ref{19-4-57})}, \textup{(\ref{19-2-78})} and \textup{(\ref{19-4-84})} hold.
\end{proposition}

Furthermore, arguments of proposition~\ref{prop-19-4-15} related to
$T\ge \bar{T}=\epsilon\mu^{-1}$ do not change either. However, we need to reconsider contribution of the segment $[-\bar{T},\bar{T}]$ as $l>2$. Again, we need to consider only the sum of terms (\ref{19-4-101})--(\ref{19-4-104}).

We already know from the proof proposition \ref{prop-19-4-15}(ii) that under assumptions $\mu h|\log h\le \epsilon$ and $\varepsilon \ge C\mu^{-1}$  this sum does not exceed  $C\mu^{-2}h^{-d}(\mu h)^{l-2}|\log (\mu h|)^{-\sigma}$.

On the other hand, all terms (\ref{19-4-101})--(\ref{19-4-104}) with $T=\bar{T}$ and their sum have the form $\mu^{-2}h^{-d}f(\mu h)$ and then repeating arguments of part (iii) of the proof proposition \ref{prop-19-4-15} with $\lambda = \mu^K h^{K-1}$ we arrive to the same estimate
$C\mu^{-2}h^{-d}(\mu h)^{l-2}|\log (\mu h|)^{-\sigma}$ for the sum as
$\mu h\le \epsilon '$. Furthermore, remark \ref{rem-19-4-16} remains true as well:

\begin{proposition}\label{prop-19-5-12} 
Statements of proposition \ref{prop-19-4-15} and remark \ref{rem-19-4-16} remain true as well.
\end{proposition}

\section{Main theorems}
\label{sec-19-5-7}

Thus we arrive to the following

\begin{theorem}\label{thm-19-5-13} 
Let assumptions \textup{(\ref{19-1-4})}--\textup{(\ref{19-1-6})}, $\textup{(\ref{19-1-25})}_{1-3}$ with \linebreak $(\bar{l},\bar{\sigma})\succeq(l,\sigma)\succeq (1,2)$, $(\bar{l},\bar{\sigma})\succeq (2,1)$ and \textup{(\ref{19-2-38})}  be fulfilled. 

 Let 
\begin{equation}
 \mu^*_2\Def\epsilon (h|\log h|)^{-1} \le \mu \le \mu^*_3\Def 
\epsilon ' h^{-1}
\label{19-5-35}
\end{equation}
with sufficiently small constant $\epsilon'>0$. Then for two framing approximations

\medskip\noindent
(i) Under  $\fN$-microhyperbolicity condition (see definition~\ref{def-19-2-5})  estimate \textup{(\ref{19-5-6})} holds;

\medskip\noindent
(ii) Under  microhyperbolicity condition  (see definition~\ref{def-19-2-4}) estimate \textup{(\ref{19-5-7})} holds.
\end{theorem}

\begin{remark}\label{rem-19-5-14}
(i) Obviously theorem~\ref{19-5-13} generalizes theorem~\ref{thm-19-5-2}; recall that our second main theorem is theorem~\ref{thm-19-5-3};

\medskip\noindent
(ii) Theorems~\ref{thm-19-5-2} and \ref{19-5-3} hold for very strong magnetic field (\ref{19-6-1}) as well but then we can still relax definition of $\fN$-hyperbolicity (see condition (\ref{19-6-13}) below).
\end{remark}

\chapter{Very strong and superstrong magnetic field}
\label{sect-19-6}

\section{Framework}
\label{sect-19-6-1}
The remaining case $\mu \ge \mu_3^*=\epsilon_1 h^{-1}$ is split into two subcases:
\begin{align}
\mu_3^*=\epsilon_1 h^{-1} \le &\mu \le \mu_4^*=C_0 h^{-1},\label{19-6-1}\\
& \mu \ge \mu_4^*=C_0 h^{-1}\label{19-6-2}
\end{align}
which we refer as \index{magnetic field!very strong}\emph{very strong\/} and \index{magnetic field!superstrong}\emph{superstrong\/} magnetic field cases respectively.

In the former subcase we know that the spatial speed of the propagation is $O(\mu^{-1})$ but we have not estimated the propagation speed with respect to the energy partition and thus we cannot use microhyperbolicity condition except $\fN$-microhyperbolicity condition with $\#\fN=1$. Further, this condition should be modified to accommodate spectral gaps.

In the latter subcase (\ref{19-6-2}) we have  nothing about propagation at all. Recall that in this case we need to consider the generalized Schr\"odinger-Pauli operator (\ref{19-1-13}).

\section{Hermitian decomposition}
\label{sect-19-6-2}

In both cases, however, we can apply a very usefully decomposition %(\ref{0-49})
\begin{gather}
U_\cT (x,y,t)=
\sum _{\alpha \in \bZ^{+r}} U_{\alpha \beta} (x'',y'',t)
\Upsilon _\alpha (x') \Upsilon _\beta (y')
\label{19-6-3}%{3-16}{0-49}
\shortintertext{with}
\Upsilon _\alpha (x')\Def \mu ^{\frac{1}{2}r}h^{-\frac{1}{2}r} 
\upsilon _\alpha\bigl( \mu ^{\frac{1}{2}}h^{-\frac{1}{2}} x'\bigr), \qquad 
\upsilon _\alpha (x)= \prod_{1\le j\le r} \upsilon _{\alpha_j}(x_j).
\label{19-6-4}
\end{gather}
Recall that  $\upsilon _{\alpha_j}$ are $1$-dimensional Hermite functions.

Note that
\begin{multline}
\cA _0 \Bigl(u_\alpha(x'')\Upsilon_\alpha (x')\Bigr)= \\
\sum_{\beta: |\beta|=|\alpha|,\,  |\beta-\alpha|=0,2}
\Bigl(\mu h  \cA^0_{\beta\alpha}(x'',\mu^{-1}hD'') +\updelta_{\alpha\beta}\,q_0(x'',\mu^{-1}hD'')\Bigr)u_\alpha\cdot  \Upsilon_\beta(x')
\label{19-6-5}
\end{multline}
where
\begin{equation}
\cA^0_{\beta\alpha}= \left\{\begin{aligned}
&\sum _j (2\alpha_j+1) b_{jj}\qquad &\text{as\ \ }\beta=\alpha,\\
&2\sqrt{\alpha_j\beta_k}\, b_{jk}\qquad &\text{as\ \ }
\beta_i-\alpha_i=\updelta_{ik}-\updelta_{ij}.
\end{aligned}\right.
\label{19-6-6}
\end{equation}
Further,
\begin{claim}\label{19-6-7}
$\bH_n \Def\bigl\{\sum_{\alpha:|\alpha|=n}v_\alpha (x'')\Upsilon_\alpha(x')\bigr\}$
are invariant subspaces of both $\cA_0$ and $\cA^0$
\end{claim}
and
\begin{equation}
\Spec \bigl(\cA^0_{\alpha\beta}(x'',\xi'')\bigr)_{\alpha,\beta}\Bigr|_{\bH_n} =\Bigl\{ \sum_j  (2\alpha_j +1)f_j,\  |\alpha |=n\Bigr\}.
\label{19-6-8}
\end{equation}

\section{Propagation}
\label{sect-19-6-3}
\subsection{Special case of constant $g^{jk},F_{jk}$}
\label{sect-19-6-3-1}

This is far the easiest case. First, we can assume without any loss of the generality that
\begin{equation}
\cA = \cA_0.
\label{19-6-9}
\end{equation}
Really,  no cubic terms appear from the ``kinetic'' part of Hamiltonian
$\sum_j g^{jk}P_jP_k$ and the first-order terms appearing from the potential are eliminated in the process of reduction. 

Further, we can skip an error $O(\varepsilon ^l |\log \varepsilon|^{-\sigma})$ and for $(l,\sigma)\succeq (2,0)$ we can skip quadratic terms appearing from the potential because in the operator sense they are less than $\mu ^{-2}(\cA^0 +1)$ and since our reduction in this case is global, we can just take approximate operator as $\cT \cA_0 \cT^*$.

Furthermore, we can assume without  any loss of the generality that $b_{jk}=f_j\updelta_{jk}$ and then
\begin{equation}
\cA \Bigl(u_\alpha (x'')\Upsilon _\alpha (x')\Bigr) =
\Bigl(\sum _j (2\alpha_j+1)f_j\mu h + q_0\Bigr)u_\alpha (x'')\cdot
\Upsilon _\alpha (x').
\label{19-6-10}
\end{equation}
Now we are dealing with $r$-dimensional scalar $\mu^{-1}h$-pseudo-differential operators and the proofs of the following two statements are obvious:

\begin{proposition}\label{prop-19-6-1} 
Let $g^{jk},F_{jk}$ be constant and one of conditions \textup{(\ref{19-6-1})}, \textup{(\ref{19-6-2})} be fulfilled. Let $Q=q(x'',\mu^{-1}hD'')$. Assume that on the support of $q$ an ellipticity condition 
\begin{equation}
|\sum _i (2\alpha _i+1)f_i \mu h +V-\tau| \ge \epsilon _1\qquad
\forall \alpha \in \bZ^{+r}.
\label{19-6-11}
\end{equation}
is fulfilled.

Then for $T\ge T_*\Def Ch|\log \mu|$ 
\begin{equation}
|F_{t\to h^{-1}\tau} \Bigl(\bar{\chi}_T (t) U_\cT \,^t\!Q_y \Bigr)|\le CT\mu^{-s}.
\label{19-6-12}
\end{equation}
\end{proposition}

\begin{proof}
An easy proof  based on our standard elliptic arguments is left to the reader.
\end{proof}

\begin{proposition}\label{prop-19-6-2} 
Let $g^{jk},F_{jk}$ be constant and one of conditions \textup{(\ref{19-6-1})}, \textup{(\ref{19-6-2})} be fulfilled. Let $Q=q(x'',\mu^{-1}hD'')$. 

\medskip\noindent
(i) Assume that on the support of $q$ a microhyperbolicity  condition  
\begin{equation}
|\sum _i (2\alpha _i+1)f_i\mu h +V-\tau| + |\nabla V|\ge \epsilon _1\qquad
\forall \alpha \in \bZ^{+r}.
\label{19-6-13}
\end{equation}
is fulfilled.
Then for $T\in [T_*,T^*]$
\begin{equation}
|F_{t\to h^{-1}\tau} \Bigl(\chi_T(t) \Gamma ''U_\cT \,^t\!Q_y\Bigr)|\le C\mu^{-s}.
\label{19-6-14}
\end{equation}

\medskip\noindent
(ii) Assume that on the support of $q$   condition  
\begin{equation}
|\sum _i (2\alpha _i+1)f_i\mu h +V-\tau| \le \epsilon  \implies  |\nabla V|\asymp \epsilon _1\nu \qquad
\forall \alpha \in \bZ^{+r}
\label{19-6-15}
\end{equation}
with $C_0 (\mu^{-1}h|\log h|)^{\frac{1}{2}}\le \nu \le \epsilon$. Then for $T\in [T_*,T^*]$ \textup{(\ref{19-6-14})} holds where now $T_*= C_0\nu^{-1}\mu h|\log \mu|$.
\end{proposition}

\begin{proof}
An easy proof  based on our standard microhyperbolicity arguments  is left to the reader.
\end{proof}

\subsection{Very strong magnetic field: general operators}
\label{sect-19-6-3-2}

Assume that  condition (\ref{19-6-1}) holds.

Consider operator $A_\cT$ as a $\mu^{-1}h$-pseudo-differential operator with respect to $x'$ with a ``matrix'' symbol with values in $\cL (\bH,\bH)$ where 
$\bH =\sL^2(\bR^r_{x'})$ and $\cL (\bH_1,\bH_2)$ is the space of linear operators from $\bH_2$ to $\bH_1$. 

Obviously, operator  $\cA$ does not belong to this class, but we can also replace $A_\cT$ by $\cA_0 + Q_1^*(A_\cT -\cA_0)Q_1$ where
$Q_1= I- \phi (\varepsilon^{-1} x',\varepsilon^{-1}\mu^{-1}hD')$ and $\phi $
is supported in $2C_1$-vicinity  and equal 1 in $C_1$-vicinity of
$0\in \bR^d_{x'',\xi''}$. Applied to $U$, $Q_1$ produces negligible output\footnote{\label{foot-19-29} After applying $F_{t\to h^{-1}\tau }\bar{\chi}_T (t) $ with $\tau \le c$ and $Ch|\log h|\le T \le ch^{-s}$.}.

\begin{proposition}\label{prop-19-6-3} 
Let condition \textup{(\ref{19-6-1})} be fulfilled. Assume that on the support of $q$ the ellipticity condition \textup{(\ref{19-6-11})} is fulfilled and  let $Q=q(x'',\mu^{-1}hD'')$.
Then for $T\ge Ch|\log h|$, $\tau \in [-\epsilon',\epsilon']$ (with small enough constant $\varepsilon'>0$) estimate \textup{(\ref{19-6-12})} holds.
\end{proposition}

Furthermore

\begin{proposition}\label{prop-19-6-4} 
Let condition \textup{(\ref{19-6-1})} be fulfilled. Assume that on the support of $q$ for a symbol $\cA_0$ the standard microhyperbolicity assumption 
\begin{equation}
\blangle\bigl(\ell \cA_0(x'',\xi'')\bigr) v, v \brangle  \ge
\epsilon_0 \1 v\1 ^2 - C\1 \cA_0 v\1^2\qquad \forall v\in \bH
\label{19-6-16}
\end{equation}
is fulfilled  and let $Q=q(x'',\mu^{-1}hD'')$. Then for $T\in [T'_*,T^*]$ estimate \textup{(\ref{19-6-14})} holds.
\end{proposition}

\begin{proof}
An easy proof  based on  the standard arguments of the proof theorem~\ref{book_new-thm-2-3-1} of \cite{futurebook} is left to the reader.
\end{proof}

\begin{remark}\label{rem-19-6-5}
Note that we need the single direction $\ell$. We just cannot localize  in $(\mu^2 a^0_{{\fn}_1},\dots, \mu^2 a^0_{{\fn}_\nu})$ even after our reduction. This is really frustrating because  Poisson brackets
$\bigl\{\mu^2 a^0_{{\fn}_j}, \mu^2 a^0_{{\fn}_k}\bigr\}$ are $O(\mu ^3 h |p|^4)$ which is $O(\mu h^3 |\log h|^2)$ even in the microlocal sense.
\end{remark}

\subsection{Superstrong magnetic field case. Reduction}
\label{sect-19-6-3-3}

Assume that  condition (\ref{19-6-2}) holds. Then we should take into account that unless  $F_{jk},g^{jk}$ are constant (the case we already considered),  generally variations of  $\mu h f_j $ are of magnitude $\mu h \gg 1$. To overcome this difficulty we need to  modify our assumptions.

First, we need to replace the Schr\"odinger operator  by the generalized Schr\"odinger-Pauli operator (\ref{19-1-13}).

Ellipticity condition would mean that 
\begin{equation}
|\sum_j \bigl(\fz_j- (2\alpha_j+1)\bigr)f_j|\ge \epsilon\qquad \forall\alpha\in \bZ^{+\,r}
\label{19-6-17}
\end{equation}
but we concentrate mainly on the opposite case
\begin{gather}
\fz_j=2\alpha_j+1\qquad\text{as\ \ }\alpha\in \fA \subset \bZ^{+\,r}, \label{19-6-18}\\
|\sum_j (\fz_j- 2\alpha_j-1 )f_j|\ge \epsilon\qquad 
\forall \alpha\in \bZ^{+\,r}\setminus \fA
\label{19-6-19}
\end{gather}
and $C_0=C_0(\epsilon)$ in condition (\ref{19-6-2}). Then obviously $\fA$ is a finite subset.

We still need to deal with the fact that we have only local reduction and thus $A_\cT$ differs from $\cA$ by operator  $\cA'$ with symbol belonging to
$\omega \sF^{0,0}$ with
\begin{equation}
\omega = \varepsilon^l|\log \mu |^{-\sigma}+
\mu h \varepsilon^{\bar{l}}|\log \mu |^{-\bar{\sigma}}.
\label{19-6-20}
\end{equation}

\begin{proposition}\label{prop-19-6-6} 
Let $b\in \sF^{0,0}$ be supported in $\{|x'|+|\xi'|\le c\varepsilon\}$. Then operator
\begin{gather}
b^\w_{\alpha\beta}: \bK \ni v \to
\blangle b^\w (v\otimes \Upsilon _\beta ) , \Upsilon _\alpha \brangle \in \bK, \qquad \bK\Def \sL^2(\bR^r_{x''})
\label{19-6-21}\\
\shortintertext{has norm}
\| b^\w_{\alpha\beta}\| \le Ce^{-\epsilon |\alpha -\beta|}.
\label{19-6-22}
\end{gather}
\end{proposition}

\begin{proof} 
Obviously, it is sufficient to prove this proposition as $x'\in \bR^1$ (with 
$x''\in \bR^r$). Then since
\begin{equation}
H_j \Upsilon_\beta =  (2\beta_j+1) \mu h\Upsilon_\beta,\qquad
H_j\Def \mu h(|\zeta_j|^2)^\w = \mu^2x_j^2+h^2D_j^2
\label{19-6-23}
\end{equation}
we conclude that
\begin{equation}
b^\w_{\alpha\beta} =(\alpha_j-\beta_j)^{-1} \bigl\{|\zeta_j|^2,b\bigr\}^\w_{\alpha\beta}
\label{19-6-24}
\end{equation}
and continuing this process we conclude that
$b^\w_{\alpha\beta} =(\alpha_j-\beta_j)^{-n} (b_{(n)})^\w_{\alpha\beta}$ with
$b_{(0)}=b$ and $b_{(n)}= \bigl\{|\zeta_j|^2,b_{(n-1)}\bigr\}$. 

One can see easily $b_{(n)}\in C^n n! \sF^{0,0}$ uniformly with respect to $n$\,\footnote {\label{foot-19-30} With simultaneous decay of the number of derivatives checked, but  everything works as long as $n\le \epsilon|\log \varepsilon|$.}
operator norm of $b_{(n)}^\w$ in $\cL (\bH\otimes \bK,\bH\otimes \bK)$ does not exceed $C^n n!$ and then operator norm of $b^\w_{\alpha\beta}$ in
$\cL (\bK,\bK)$ does not exceed  $|\alpha_j-\beta_j|^{-n} C^n n!$. It reaches minimum  $e^{-\epsilon |\alpha_j-\beta_j|}$ as $n =C^{-1}|\alpha_j-\beta_j|$.
\end{proof}

\subsubsection{No $2$-nd order resonances case.}\label{sect-19-6-3-2-1} 
Assume first that all  $f_j$ are disjoint. Then
\begin{gather}
\cA_0 =\sum_{1\le j \le r} b_j(x'',\mu^{-1}hD'')
\bigl(\mu^2 x_j^2 + h^2 D_j^2\bigr) +b_0(x'',\mu^{-1}hD'')
\label{19-6-25}
\shortintertext{and}
\cA_0 \bigl(v(x'')\Upsilon_\alpha (x')\bigr)=(W_\alpha v)\Upsilon_\alpha (x')
\label{19-6-26}
\end{gather}
where $W_\alpha$ is $(\mu^{-1}h)$-pseudo-differential operator with the principal symbol
\begin{equation}
W_{\alpha} = \bigl(V + \sum_j (2\alpha_j +1-\fz_j) f_j \mu h \bigr)\circ \Psi_0.
\label{19-6-27}
\end{equation}

\begin{proposition}\label{prop-19-6-7}
Let ellipticity condition \textup{(\ref{19-6-17})} be fulfilled and let condition \textup{(\ref{19-6-2})} be fulfilled at $\bar{x}$ with $C_0=C_0(\epsilon)$. Then modulo $O(\mu^{-s})$
\begin{equation}
F_{t\to h^{-1}\tau} \bar{\chi}_T (t)U_{\alpha\beta}\equiv 0
\qquad \forall \tau:|\tau|\le \epsilon' \mu h\qquad\text{in\ \ }\Omega\times\Omega
\label{19-6-28}
\end{equation}
where $\Omega=B(\bar{z}'', C_0\varepsilon)\subset \bR^{2r}_{x'',\xi''}$ is a domain  in which reduction is done, $\bar{z}''=(\bar{x}'',\bar{x}'')=\Psi_0^{-1}(\bar{x})$.
\end{proposition}

\begin{proof}
Proof due to the standard elliptic arguments is left to the reader.
\end{proof}

\begin{proposition}\label{prop-19-6-8} 
Let conditions  \textup{(\ref{19-6-18})}--\textup{(\ref{19-6-19})} be fulfilled  and let condition \textup{(\ref{19-6-2})} be fulfilled at $\bar{x}$ with $C_0=C_0(\epsilon)$.  Assume that all $f_j$ are disjoint. Let $T=\mu ^m$ with large enough $m$.   Then

\medskip\noindent
(i) As $\alpha\notin \fA$
\begin{multline}
F_{t\to h^{-1}\tau} \bar{\chi}_T (t)U_{\alpha\beta}\equiv
\sum _{\alpha'\in \fA} E_{\alpha\alpha', x}
F_{t\to h^{-1}\tau} \bar{\chi}_T  (t) U_{\alpha'\beta} \\ \forall \tau:|\tau|\le \epsilon' \mu h\qquad\text{in\ \ }\Omega\times\Omega
\label{19-6-29}
\end{multline}
where $E_{\alpha\alpha'} = E_{\alpha\alpha'} (\tau , x'',\mu^{-1}hD'')$ are pseudo-differential operators with the norms
\begin{equation}
\| E_\alpha\| \le C (\mu h)^{-1}
\Bigl((\mu^{-1}h)^{\frac{1}{2}\updelta (\alpha)}
+e^{-|\alpha -\bar{\alpha}|}\omega \Bigr)
\label{19-6-30}
\end{equation}
and a similar dual equation holds; here 
\begin{equation}
\updelta (\alpha) \Def  \min_{\alpha'\in \fA} |\alpha-\alpha'|.
\label{19-6-31}
\end{equation}
(ii) As $\alpha\notin \fA$ and $\beta\ne \notin \fA$
\begin{multline}
F_{t\to h^{-1}\tau} \bar{\chi}_T (t)U_{\alpha\beta}\equiv
\sum_{\alpha',\beta'\in \fA} E_{\alpha\alpha', x} F_{t\to h^{-1}\tau} \bar{\chi}_T  (t) U_{\alpha'\beta'} \,^t\!E_{\beta\beta', y}\\[2pt] 
\forall \tau:|\tau|\le \epsilon' \mu h\qquad\text{in\ \ }\Omega\times\Omega;
\label{19-6-32}
\end{multline}
(iii) Further, $\mathbf{U}\Def \{U_{\alpha \beta }\}_{\alpha ,\beta \in \fA}$ satisfies
\begin{multline}
F_{t\to h^{-1}\tau} \bar{\chi}_T (t)
\Bigl( \bigl(hD_t - W_{\alpha}(x'',\mu^{-1}hD'')\bigr) U_{\alpha\beta}- \\
\shoveright{\sum_{\alpha'\in \fA} W'_{\alpha\alpha'} ( x'',\mu^{-1}hD'')U_{\alpha'\beta}\Bigr) \equiv 0}\\
\forall \tau:|\tau|\le \epsilon' \mu h\qquad \text{in\ \ }\Omega\times\Omega
\label{19-6-33}
\end{multline}
and adjoint equation with respect to $y''$ where $W_{\alpha}$ is defined by \textup{(\ref{19-6-27})} and $W'_{\alpha\alpha'}$ are $\mu^{-1}h$-pseudo-differential operator with  symbol belonging to
$(\omega +h^2)\sF^{0,0}$.
\end{proposition}

\begin{proof} 
(a) Plugging representation (\ref{19-6-3}) into equation
$\bigl(hD_t-A_{\cT, x}\bigr)У=0$ and using (\ref{19-6-26})--(\ref{19-6-27}), form of  $\cA$ and applying  proposition \ref{prop-19-6-6} to $A_\cT - \cA$ we arrive to the system
\begin{equation}
\bigl(hD_t -W_\alpha-W'_\alpha \bigr)U_{\alpha\beta}  \equiv
\sum_{\gamma\ne\alpha } B_{\alpha\gamma}\,  U_{\gamma \beta}
\qquad\text{in\ \ }\Omega\times \bR^{2r}
\label{19-6-34}
\end{equation}
with $W'_\alpha$ not including $hD_t$. Here and below all operators are acting with respect to $x''$. We also get an adjoint system
\begin{equation}
U_{\alpha\beta}\,^t\!\bigl(hD_t -W_\beta-W'_\beta\bigr)  \equiv
\sum_{\gamma\ne\\beta}   U_{\alpha\gamma}\,^t\!B_{\gamma\beta}
\qquad\text{in\ \ } \bR^{2r}\times \Omega
\tag*{$\textup{(\ref*{19-6-34})}'$}\label{19-6-34-'}
\end{equation}
with operators acting with respect to $y''$.

Applying $F_{t\to h^{-1}\tau}\bar{\chi}_T(t)$ with $T$ described above, we arrive to the system
\begin{equation}
\bigl(\tau -W_\alpha -W'_\alpha\bigr)v_{\alpha\beta}- \equiv
\sum_{\gamma\ne \alpha} B_{\alpha\gamma}\,  v_{\gamma \beta}
\qquad\text{in\ \ }\Omega\times \bR^{2r},
\label{19-6-35}
\end{equation}
with $v_{\alpha\beta}=F_{t\to h^{-1}\tau}\bar{\chi}_T(t)U_{\alpha\beta}$
and a similar dual system.

Due to assumptions (\ref{19-6-18})--(\ref{19-6-19})  all operators 
$\bigl(\tau -W_\alpha - W'_\alpha\bigr)$ with 
$\alpha \notin \fA$ are elliptic with inverse operators norms of magnitude
$(\mu h\updelta(\alpha))^{-1}$ as  $|\tau|\le \epsilon \mu h$.

Note that even in the worst case $B_{\alpha\beta}$ are operators with the symbols belonging to $\mu h(\mu^{-1}h|\alpha|)^{\frac{1}{2}}\sF^{0,1}$. Then by the successive approximations  
\begin{equation}
v_{\alpha \beta}= \sum _{\alpha'\in \fA}E_{\alpha \alpha'}v_{\alpha'\beta}\qquad \text{as\ \ } \alpha\notin \fA
\label{19-6-36}
\end{equation}
where each $E_{\alpha\alpha'} $ is the sum of terms of the following type:
\begin{equation}
(\tau - W_\alpha )^{-1}B'_{\alpha \alpha_1}(\tau - W_{\alpha _1})^{-1}
B'_{\alpha _1\alpha_2}\cdots \\(\tau - W_{\alpha _k})^{-1} B'_{\alpha_k \alpha'}
\label{19-6-37}
\end{equation}
with $k\ge 0$ where $B'_{\gamma\gamma'}$ are  operators with symbols
belonging to $\mu h(\mu^{-1}h|\alpha|)^{\frac{1}{2}}\sF^{0,1}$. For the sake of simplicity we included $W'_\alpha$ into $W_\alpha$.

Similarly,
\begin{equation}
v_{\alpha\beta}= \sum_{\beta'\in \fA} v_{\alpha\beta' }\,^t\!E_{\beta\beta'}\qquad \text{as\ \ } \alpha\notin \fA
\tag*{$\textup{(\ref{19-6-36})}'$}\label{19-6-36-'}
\end{equation}
Then statements (i),(ii) follow then from  formulae (\ref{19-6-36}), \ref{19-6-36-'} and (\ref{19-6-37}).
 
\medskip\noindent
(b) Plugging (\ref{19-6-27}), (\ref{19-6-3}) into equation $(hD_t-A_{\cT\, x})u=0$ and a dual  equation with respect to $y$, we arrive to an equation similar to (\ref{19-6-33}) (as well as an adjoint equation with respect to $y$):
\begin{equation}
\bigl(\tau - W_{\alpha} \bigr)v_{\alpha\beta} \equiv \sum \cW_{\alpha\alpha'}
v_{\alpha' \beta} \qquad\text{in\ \ }\Omega\times \bR^{2r}
\label{19-6-38}
\end{equation}
where $\cW_{\alpha\alpha'}$ is a sum of terms of (\ref{19-6-37}) type with an extra factor  on their left:
\begin{equation}
(\mu^{-1}h)^{\frac{p}{2}} B'_{\alpha\alpha_1}
(\tau - W_{\alpha_1} )^{-1}
B'_{\alpha _1\alpha_2}\cdots (\tau - W_{\alpha _k})^{-1} 
B'_{\alpha_k \alpha'}
\label{19-6-39}
\end{equation}
with $B'_{\bar{\alpha}\alpha}$ of the same type as above; $p=0$ here. Note that
\begin{equation}
(\tau - W_\gamma )^{-1}=(W_{\alpha'} - W_\gamma)^{-1}-
(\tau - W_\gamma )^{-1}(\tau - W_{\alpha'} )(W_{\alpha'} - W_\gamma)^{-1}
\label{19-6-40}
\end{equation}
in $\Omega\times \bR^{2r}$.

We apply this formula to $(\tau-W_{\alpha_k})^{-1}$. Then we drag
$(\tau - W_{\alpha'} )$ to the right. If it perishes at commuting, we get the same expression with the same $k$ but with $p$ replaced by $p+1$ because in the commuting we gain $\varepsilon^{-1}\mu^{-1}h$ factor. Continuing this process with  $(\tau-W_{\alpha_{k-1}})$ etc we arrive to
\begin{multline}
\bigl(\tau - W_{\alpha} \bigr)v_{\alpha\beta} \equiv 
\sum_{\alpha'\in \fA} \Bigl(W''_{\alpha \alpha'} +
\cW ''_{\alpha\alpha'} \,  (\tau - W_{\alpha'}) + \cW '''_{\alpha\alpha'}\Bigr) v_{\alpha' \beta}\\
\text{in\ \ }\Omega\times \bR^{2r}
\label{19-6-41}
\end{multline}
where $W''_{\alpha\alpha'}$ does not depend on $\tau $, 
$\cW '''_{\alpha\alpha'}$ is the sum of products of (\ref{19-6-39}) type with $p$ replaced by $p+1$,
$ \cW ''_{\alpha\alpha'}$ is the sum of products of (\ref{19-6-39}) type with $p$ replaced by $p+2$ and $k$ replaced by $k-1$. Considering system (\ref{19-6-41}) as a $(\#\fA)\times(\#\fA)$-matrix equation and multiplying by $(I+\cW'' )$ one can rewrite this system as (\ref{19-6-38}) with $W_{\alpha}$ replaced by $W_{\alpha}+W''_{\alpha\alpha}$ and $\cW'''_{\alpha\alpha'}$ of the same type with $p$ replaced by $p+1$.

Continuing this process we arrive to negligible $\cW'''_{\alpha\alpha'}$. Statement (iii) is proven.
\end{proof}

\begin{remark}\label{rem-19-6-9}  
(i) It follows from our construction that an added correction modulo $(\omega + \mu ^{-1}h)\sF^{0,0}$  is equal to
\begin{equation}
W'_{\alpha\alpha'} = \sum _{\beta:|\alpha-\beta|=|\alpha-\beta|=1}
B_{\alpha\beta}(W_\beta)^{-1}B_{\beta\alpha'}
\label{19-6-42}
\end{equation}
with $B_{\alpha\beta}$ appearing exclusively from cubic terms in $\cA$.

\medskip\noindent
(ii) Further, note that these terms are $b_{m;jk}Z_m^* Z_jZ_k$ and
$b_{m;jk}^*Z_m Z_j^*Z_k^*$ with $f_m$ not disjoint from $f_j+f_k$. Conditions
(\ref{19-6-18})--(\ref{19-6-19}) and $\#\fA=1$, $\alpha\in \fA$ imply that  $\alpha_m=0$ \underline{and} either $\bar{\alpha}_j=0$ or $\bar{\alpha}_k=0$ if $j\ne k$ and $\alpha_j=0,1$ as $j=k$. Therefore cubic terms, applied to $v(x'')\Upsilon_{\bar{\alpha}} (x')$, produce $0$, and therefore $W'_{\alpha\alpha'} \in (\omega + \mu ^{-1}h)\sF^{0,0}$.

However, correction
$W'_{\alpha\alpha'} \equiv  h^2 \sum_j \kappa_j Z_j^*Z_j \;\mod O(\omega +\mu^{-1}h )$ has been already generated in the process of reduction in Section~\ref{sect-19-3}. Still we see that it  $0$ in the case of constant $g^{jk},F_{jk}$.
\end{remark}

\subsubsection{$2$-nd order resonances case.}\label{sect-19-6-3-4}

Consider now a more general case when some of $f_j$ are not disjoint; assume that $f_j$ is disjoint from the rest of eigenvalues as $j=1,\ldots,p$ and is not as $j=p+1,\ldots,r$. Then
\begin{multline}
\cA_0= \sum_{1\le j\le p} f_j(x'',\mu^{-1}hD'')\bigl(\mu^2x_j^2+h^2D_j^2\bigr) + \\[2pt] 
\sum_{p+1\le j,k \le r}b_{jk}(x'',\mu^{-1}hD'')
\bigl(Z_j^*Z_k +\mu h \delta_{jk}\bigr) +b_0(x'',\mu^{-1}hD'').
\label{19-6-43} 
\end{multline}

Assume that
\begin{equation}
\alpha_j=0 \quad\text{as\ \ } j=p+1,\ldots,r\qquad\forall \alpha\in \fA.
\label{19-6-44}
\end{equation}
Note that conditions (\ref{19-6-18})--(\ref{19-6-19}) and $\#\fA=1$  yield (\ref{19-6-44}).

Due to (\ref{19-6-18})--(\ref{19-6-19}) and (\ref{19-6-44}) $\cA_0$ restricted to $v(x'')\Upsilon_{\bar{\alpha}}(x')$ equals
\begin{equation}
\cA_0= 
\sum_{1\le j\le p} f_j(x'',\mu^{-1}hD'')\bigl(2\bar{\alpha}_j+1\bigr)\mu h  + b_0(x'',\mu^{-1}hD'');
\label{19-6-45} 
\end{equation}
then arguments of the proof proposition \ref{prop-19-6-7} still work but instead of the individual subspaces $\sH_\alpha\Def \bigl\{v(x'')\Upsilon_\alpha(x')\bigr\}$ for $\alpha=(\alpha ';\alpha'') =(\alpha_1,\dots, \alpha_p;\alpha_{p+1},\dots,\alpha_r)$ one should consider
subspaces 
\begin{equation*}
\sH_{\alpha ',n}\Def \bigl\{\sum _{|\alpha ''|=n} v_{(\alpha';\alpha'')} (x'')\,\Upsilon_{(\alpha';\alpha'')}(x')\bigr\}.
\end{equation*}

Therefore we have

\begin{proposition}\label{prop-19-6-10} 
Propositions \ref{prop-19-6-7} and \ref{prop-19-6-8} holds even if $f_i$ are not necessary disjoint but \textup{(\ref{19-6-44})} is fulfilled.
\end{proposition}

\subsection{Superstrong magnetic field case. Propagation}
\label{sect-19-6-3-5}

So under assumption $\#\fA=1$ basically we reduced our operator to a single $r$-dimensional $\mu^{-1}h$-pseudo-differential operator $W_{\bar{\alpha}}$ with a principal symbol $V\circ \Psi_0$ while under assumption (\ref{19-6-44}) we have a matrix-operator with the diagonal principal part instead. Then in the framework of (\ref{19-6-18})--(\ref{19-6-19}) we can impose an ellipticity 
\begin{gather}
|V|\ge \epsilon_0\label{19-6-46} \\
\shortintertext{or a microhyperbolicity}
|V|+|\nabla V| \ge \epsilon_0\label{19-6-47} \\
\intertext{or as $(l,\sigma)\succeq (2,0)$ a non-degeneracy assumption}
|V|+|\nabla V| + |\det \Hess V|\ge \epsilon_0.\label{19-6-48} 
\end{gather}

\begin{proposition}\label{prop-19-6-11} 
Let conditions \textup{(\ref{19-6-2})}, \textup{(\ref{19-6-18})}--\textup{(\ref{19-6-19})}, and \textup{(\ref{19-6-44})} be fulfilled. Then under assumption \textup{(\ref{19-6-46})} estimate \textup{(\ref{19-6-28})} holds.
\end{proposition}

\begin{proof}
The standard elliptic arguments applied to equation (\ref{19-6-33}) yield that under condition (\ref{19-6-46} ) estimate (\ref{19-6-28}) holds as
$\tau \le \epsilon'$, $\alpha,\beta\in \fA$. Then in virtue of propositions~\ref{prop-19-6-8} and \ref{prop-19-6-10}  this is true for all $\alpha,\beta$. Easy details we leave to the reader.
\end{proof}

\begin{proposition}\label{prop-19-6-12} 
Let conditions \textup{(\ref{19-6-2})},  \textup{(\ref{19-6-18})}--\textup{(\ref{19-6-19})}, and \textup{(\ref{19-6-44})} be fulfilled.

Then for $T_*\le T \le T^*=\epsilon_0 T$, $T\ge C \varepsilon^{-1} h|\log \mu |$ and large enough constant $M$
\begin{equation}
|F_{t\to h^{-1}\tau }
\bar{\chi} _T(t) Q'_xU_{\alpha\beta} (x,y,t) \,^t\!Q''_y |
\le C\mu^{-s}\qquad \forall \tau \le c.
\label{19-6-49} 
\end{equation}
as $\tau \le \epsilon'$ and distance between supports of (symbols of) $Q'$ and $Q''$ is at least $C\mu^{-1}T$.
\end{proposition}

\begin{proof}
The standard propagation arguments applied to equation (\ref{19-6-33}) yield that (\ref{19-6-49} ) holds as $\alpha,\beta\in \fA$. Then in virtue of propositions~\ref{prop-19-6-8} and \ref{prop-19-6-10}  this is true for all $\alpha,\beta$. Easy details we leave to the reader.
\end{proof}

\begin{proposition}\label{prop-19-6-13} 
Let conditions \textup{(\ref{19-6-2})} and \textup{(\ref{19-6-18})}--\textup{(\ref{19-6-19})}, \textup{(\ref{19-6-44})} be fulfilled. Further, let microhyperbolicity condition \textup{(\ref{19-6-47})} be fulfilled.

Then for $T_*\le T \le T^*=\epsilon_0 T$, $T\ge C\varepsilon^{-1} h|\log \mu|$ and large enough constant $M$
\begin{equation}
|F_{t\to h^{-1}\tau }
\bar{\chi} _T(t) Q'_xU_{\alpha\beta} (x,y,t) \,^t\!Q''_y |
\le C\mu^{-s}\qquad \forall \tau \le c.
\label{19-6-50} 
\end{equation}
as $\tau \le \epsilon'$ and distance between the most distant points of supports of (symbols of) $Q'$ and $Q''$ is at most $C\epsilon \mu^{-1}T$.
\end{proposition}

\begin{proof}
The standard propagation arguments applied to equation (\ref{19-6-33}) yield that (\ref{19-6-50} ) holds as $\alpha,\beta\in \fA$. Then in virtue of propositions~\ref{prop-19-6-8} and \ref{prop-19-6-10}  this is true for all $\alpha,\beta$. Easy details we leave to the reader.
\end{proof}

\begin{corollary}\label{cor-19-6-14} 
(i) In the framework of proposition~\ref{prop-19-6-12} 
\begin{equation}
|F_{t\to h^{-1}\tau }
\bar{\chi} _T(t) \psi'(x) U (x,y,t) \psi''(y) |
\le C\mu^{-s}\qquad \forall \tau \le c
\label{19-6-51} 
\end{equation}
as $\tau \le \epsilon'$ and distance between supports of (symbols of) $\psi'$ and $\psi''$ is at least $C\mu^{-1}T$ under additional assumption 
\begin{equation}
T\ge C_0(\mu h|\log h|)^{\frac{1}{2}};
\label{19-6-52} 
\end{equation}
(ii) In the framework of proposition~\ref{prop-19-6-13} 
\begin{equation}
|F_{t\to h^{-1}\tau }
\bar{\chi} _T(t) \psi'(x) U (x,y,t) \psi''(y) |
\le C\mu^{-s}\qquad \forall \tau \le c
\label{19-6-53} 
\end{equation}
as $\tau \le \epsilon'$ and distance between the most distant points of supports of (symbols of) $Q'$ and $Q''$ is at most $C\epsilon \mu^{-1}T$ under additional assumption \textup{(\ref{19-6-52})}.
\end{corollary}

\begin{corollary}\label{cor-19-6-15} 
In the framework of proposition~\ref{prop-19-6-13} 
\begin{equation}
|F_{t\to h^{-1}\tau }
\bar{\chi} _T(t) \Gamma \psi'(x) U (x,y,t) \psi''(y) |
\le C\mu^{-s}\qquad \forall \tau \le c.
\label{19-6-54} 
\end{equation}
\end{corollary}

\begin{remark}\label{rem-19-6-16} 
(i) One can prove easily that under assumptions  (\ref{19-6-2}), (\ref{19-6-18})--(\ref{19-6-19}) and (\ref{19-6-44}) as $(l,\sigma)\succeq (2,0)$ singularities propagate along trajectories
\begin{equation}
\frac {dx_j} {dt}  = \mu^{-1}\sum _k \upphi^{jk} \partial_{x_k} V 
\label{19-6-55} 
\end{equation}
($\tau:|\tau|\le \epsilon'$)
where $(\upphi^{jk})=(F_{jk})^{-1}$. In particular, $V_{\bar{\alpha}}$ is an integral; 

\medskip\noindent
(ii) For $d=2$ even condition $(l,\sigma)\succeq (2,0)$ is not necessary.
\end{remark}

\begin{remark}\label{rem-19-6-17} 
As $g^{jk}=\const, F_{jk}=\const $ all our results remain true even if there are $2$-nd  order resonances and (\ref{19-6-44}) fails.
\end{remark}

We leave to the reader the following

\begin{problem}\label{Problem-19-6-18} 
Formulate and prove similar results as $|\nabla V|\asymp \nu \ge C\varepsilon$ (instead of $|\nabla V|\asymp 1$).
\end{problem}

\section{Short-range estimates}
\label{sect-19-6-4}

\subsection{Very strong magnetic field}
\label{sect-19-6-4-1}

As before the  case $\epsilon_1 h^{-1}\le \mu \le C_0 h^{-1}$ is just a variation of the strong magnetic field case but with the mandatory assumption $\#\fN=1$ and the microhyperbolicity condition required only for  $\tau=\sum_j (2\alpha_j+1)\mu h f_j$ with $\alpha\in \bZ^{+\,r}$. Using the same arguments as in the proof proposition \ref{prop-19-4-9}  we arrive to

\begin{proposition}\label{prop-19-6-19} 
Let $\epsilon h^{-1}\le \mu \le C_0 h^{-1}$,
$\varepsilon \ge  (\mu^{-1}h|\log h|)^{\frac{1}{2}}$, and
$T_*=C\varepsilon ^{-1}h|\log h|$.

Let microhyperbolicity condition \textup{(\ref{19-6-13})} be fulfilled and let $T\in [h^{1-\delta}, T_*]$  with an arbitrarily small exponent $\delta >0$. Then both statements (i),(ii) of proposition \ref{prop-19-4-9} remain true.
\end{proposition}

\begin{proof}
Easy details are left to the reader.
\end{proof}

\subsection{Superstrong magnetic field}
\label{sect-19-6-4-2}

Let us assume now that  $\mu \ge \mu_4^*=C_0h^{-1}$ and we consider generalized Schr\"odinger-Pauli operator (\ref{19-1-13}).

We need to assume now that \underline{either} $g^{jk},F_{jk}$ are constant  \underline{or} conditions (\ref{19-6-18})--(\ref{19-6-19}), (\ref{19-6-44}) are fulfilled. Exactly the same method as before leads us to similar results:

\begin{proposition}\label{prop-19-6-20} 
Let $\mu\ge C_0 h^{-1} $, $\varepsilon \ge (\mu^{-1}h|\log h|)^{\frac{1}{2}}$, and  $T_*=C\varepsilon ^{-1}h|\log h|$. Let conditions \textup{(\ref{19-6-18})}--\textup{(\ref{19-6-19})}, \textup{(\ref{19-6-44})}, and microhyperbolicity condition \textup{(\ref{19-6-47})} be fulfilled.

\medskip\noindent
(i) Let $l>1$. Then for $T\in [C\mu^\delta h, T_*]$ with an arbitrarily small exponent $\delta >0$
\begin{align}
& |\phi (hD_t) \chi_T(t) ( \Gamma U\psi_y )|\ \le
C \mu^r h^{-r}\bigl( \frac{h}{T} \bigr)^{l-1}
\bigl(1+ \frac{h}{T\varepsilon} \bigr)^{-s}|\log \frac{h}{T} |^{-\sigma} \label{19-6-56}\\
\shortintertext{and}
&|F_{t\to h^{-1}\tau} \chi_T(t) ( \Gamma U\psi_y )|\le
C \mu^r h^{1-r}\bigl(\frac{h}{T} \bigr)^{l-1}
\bigl(1+\frac{h}{T\varepsilon } \bigr)^{-s}
|\log \frac{h}{T} |^{-\sigma}\label{19-6-57}
\end{align}
as $|\tau|\le\epsilon$ with arbitrarily large exponent $s$.

\medskip\noindent
(ii) Let $l=1,\sigma\ge 2$. Then for $T\in [Ch^{1-\delta}, T_*]$ with an arbitrarily small exponent $\delta >0$ and small enough constant $T_*$
\begin{multline}
 |\phi (hD_t) \chi_T(t) ( \Gamma U\psi_y)| \le\\
C \mu ^r h^{-r}|\log \frac{h}{T} |^{-\sigma}
\bigl(1+{\frac {T\varepsilon} h}\bigr)^{-s}+
C \mu^r h^{-r}|\log \frac{h}{T} |^{-s},
\label{19-6-58}
\end{multline}
and
\begin{multline}
|F_{t\to h^{-1}\tau} \chi_T(t) ( \Gamma U\psi_y )|\le\\
C\mu^r h^{1-r}|\log \frac{h}{T} |^{-\sigma}
\bigl(1+{\frac {T\varepsilon} h}\bigr)^{-s}+
C \mu ^r h^{1-r}|\log \frac{h}{T} |^{-s}
\label{19-6-59}
\end{multline}
as $|\tau|\le\epsilon$ with arbitrarily large exponent $s$.
\end{proposition}

\begin{remark}\label{rem-19-6-21} 
As $g^{jk}=\const, F_{jk}=\const $ all our results remain true even if there are $2$-nd  order resonances and (\ref{19-6-44}) fails.
\end{remark}

\begin{problem}\label{problem-19-6-22} 
Formulate and prove similar results as $|\nabla V|\asymp \nu \ge C\varepsilon$ (instead of $|\nabla V|\asymp 1$).
\end{problem}

\section{Calculations}
\label{sect-19-6-5}

Now propositions~\ref{prop-19-6-19} and~\ref{prop-19-6-20} imply immediately

\begin{proposition}\label{prop-19-6-23} 
In the frameworks of propositions~\ref{prop-19-6-19} and~\ref{prop-19-6-20} 
the following estimates hold:
\begin{equation}
|F_{t\to h^{-1}\tau }\bar{\chi}_T(t)\Gamma (U\psi )| \le C\mu^r h^{1-r},
\label{19-6-60}
\end{equation}
\begin{multline}
T^{-1} |F_{t\to h^{-1}\tau }\chi_T(t) \Gamma
\bigl(\bar{\cG}^\pm(A_\cT - \bar{A}_\cT) \cG^\pm
\updelta(t)\cK_\psi\bigr)|\le \\
C\mu^{r-1}h^{1-r}  \bigl(\frac{h}{T}\bigr)^{l-2}|\log \mu |^{-\sigma}
\label{19-6-61}
\end{multline}
and
\begin{multline}
|\Gamma (\psi \tilde{e}) (\tau)- h^{-1} \int_{-\infty}^\tau
\bigl(F_{t \to h^{-1}\tau'} \bar{\chi}_T(t) \Gamma (\psi u)\bigr)\,d\tau'| \le C\mu^{r-1}h^{1-r}\\
\forall \tau:|\tau|\le \epsilon
\label{19-6-62}
\end{multline}
with $T\in [T_* , T^* ]$, $T_*=Ch\mu^\delta$, $T^*\Def \epsilon \mu$, $|\tau|\le \epsilon$.
\end{proposition}

In turn, proposition~\ref{prop-19-6-23} immediately implies

\begin{proposition}\label{prop-16-6-24} 
In the frameworks of propositions~\ref{prop-19-6-19} and~\ref{prop-19-6-20}
\begin{multline}
|h^{-1} \int_{-\infty}^\tau
\bigl(F_{t \to h^{-1}\tau'} \bar{\chi}_T(t) \Gamma (\psi u)\bigr)\,d\tau' -\\
(2\pi)^{-r}\mu^{r-1}h^{-r}\Tr ' 
\uptheta \bigl(\tau - \cA(x'',\xi'')\bigr)\psi^0(x'',\xi'')\, dx''d\xi''|\le \\[2pt]
C\mu^{r-1}h^{1-r} + C\mu^r h^{-r}\varepsilon ^l|\log \mu|^{-\sigma}\qquad \forall \tau:|\tau|\le \epsilon.
\label{19-6-63}
\end{multline}
\end{proposition}

\section{Main theorems}
\label{sect-19-6-6}

Finally we immediately arrive to two theorems~\ref{thm-19-6-25} and \ref{thm-19-6-26} below:

\begin{theorem}\label{thm-19-6-25} 
Let assumptions \textup{(\ref{19-1-4})}--\textup{(\ref{19-1-6})}, $\textup{(\ref{19-1-25})}_{1-3}$ with \linebreak $(\bar{l},\bar{\sigma})\succeq(l,\sigma)\succeq (1,2)$, $(\bar{l},\bar{\sigma})\succeq (2,1)$ be fulfilled. Let 
\begin{equation}
 \mu^*_3\Def\epsilon h^{-1} \le \mu \le \mu^*_4\Def  C_0 h^{-1}
\label{19-6-64}
\end{equation}
with sufficiently small constant $\epsilon'>0$. Assume that the microhyperbolicity condition~\textup{(\ref{19-6-16})} is fulfilled.

Then there are two framing approximations\footref{book_new-foot-18-16} (see Chapter~\ref{book_new-sect-18} of \cite{futurebook}) such that:
\begin{multline}
\R^\MW  \Def |\int \Bigl( {\tilde e}(x,x,0) - 
h^{-d}\cN^\MW (x,0)\Bigr)\psi (x)\, dx |\le\\
C\mu ^{r-1}h^{1-r} +
C\mu ^r h ^{-r}(\mu ^{-1} h)^{\frac{l}{2}} |\log \mu |^{\frac{l}{2}-\sigma}+
C\mu ^{r+1}h^{1-r} (\mu^{-1} h)^{\frac{\bar{l}}{2}}
|\log \mu |^{\frac {\bar{l}} {2}-\bar{\sigma}}.
\label{19-6-65}
\end{multline}
\end{theorem}

\begin{theorem}\label{thm-19-6-26} 
Let assumptions \textup{(\ref{19-1-4})}--\textup{(\ref{19-1-6})}, $\textup{(\ref{19-1-25})}_{1-3}$ with \linebreak $(\bar{l},\bar{\sigma})\succeq(l,\sigma)\succeq (1,2)$, $(\bar{l},\bar{\sigma})\succeq (2,1)$ be fulfilled. Let  
\begin{equation}
 \mu^*_4\Def C_0 h^{-1} \le \mu 
\label{19-6-66}
\end{equation}
with sufficiently small constant $\epsilon'>0$. Consider generalized Schr\"odinger-Pauli operator. Assume that \textup{(\ref{19-6-18})}--\textup{(\ref{19-6-19})}  and \underline{either} microhyperbolicity condition \textup{(\ref{19-6-47})} \underline{or} nondegeneracy condition \textup{(\ref{19-6-48})} are fulfilled.

Then there are two framing approximations\footref{book_new-foot-18-16} (see Chapter~\ref{book_new-sect-18} of \cite{futurebook}) such that:
\begin{multline}
\R^\MW_2 \Def |\int \Bigl( {\tilde e}(x,x,0) - \\
h^{-d}\cN^\MW (x,0)-
h^{-d}\cN_{2\corr} ^\MW (x,0)\Bigr)\psi (x)\, dx |\le\\
C\mu ^{r-1}h^{1-r} +
C\mu ^r h ^{-r}(\mu ^{-1} h)^{\frac{l}{2}} |\log \mu |^{\frac{l}{2}-\sigma}+
C\mu ^{r+1}h^{1-r} (\mu^{-1} h)^{\frac{\bar{l}}{2}}
|\log \mu |^{\frac {\bar{l}} {2}-\bar{\sigma}}
\label{19-6-67}
\end{multline}
with
\begin{multline}
h^{-d}\cN_{2\corr}^\MW (x,\tau) \Def  \\
(2\pi )^{-r} \mu ^rh^{2-r}
\sum _{\alpha \in \fA\subset \bZ^{+r}}
\updelta \Bigl(\tau - \sum_j (2\alpha_j +1) f_j\mu h -V\Bigr) W'_\alpha
\label{19-6-68}
\end{multline}
with $W'_\alpha=O(1)$.
\end{theorem}

\begin{remark}\label{rem-19-6-27}
The exact formula for $W'_\alpha$ is rather complicated; however $W'_\alpha=0$ provided $g^{jk}$, $F_{jk}$ are constant.
\end{remark}
 
\begin{Problem}\label{Problem-19-6-28}
Derive explicit formula for  $W'_\alpha$.
\end{Problem}

Finally, consider  spectral gaps:

\begin{theorem}\label{thm-19-6-29} 
Let assumptions \textup{(\ref{19-1-4})}--\textup{(\ref{19-1-6})}, $\textup{(\ref{19-1-25})}_{1-3}$ with \linebreak $(\bar{l},\bar{\sigma})\succeq(l,\sigma)\succeq (1,2)$, $(\bar{l},\bar{\sigma})\succeq (2,1)$ be fulfilled.  Consider two cases:

\medskip\noindent
(i) Let assumptions  \textup{(\ref{19-6-64})} and \textup{(\ref{19-6-11})} be fulfilled. 

\medskip\noindent
(ii) Let assumptions  \textup{(\ref{19-6-66})}, \textup{(\ref{19-6-18})}--\textup{(\ref{19-6-19})} and either $\#\fA=0$ or \textup{(\ref{19-6-46})} be fulfilled. 

\medskip
Then in both cases (i), (ii)  there are two framing approximations\footref{book_new-foot-18-16} (see Chapter~\ref{book_new-sect-18} of \cite{futurebook}) such that 
\begin{equation}
\R^\MW \le C\mu ^{-s}
\label{19-6-69}
\end{equation}
with arbitrarily large $s$.
\end{theorem}

\chapter{Simple degeneration}
\label{sect-19-7}

Now we want to consider the case of the degeneration but completely different from one when $g^{jk}$, $F_{jk}$ constant and $V$ having non-degenerate critical points.

\section{Framework}
\label{sect-19-7-1}

Assume now that
\begin{claim}\label{19-7-1}
There are no $2$-nd order resonances i.e. eigenvalues $f_j$ of $(F^j_k)$ are simple
\end{claim}
albeit the microhyperbolicity condition (see definition \ref{def-19-2-4}) is violated i.e. there exists $\boldlambda=(\lambda_1,\ldots,\lambda_r)$ with 
\begin{gather}
\lambda_1\ge 0,\ldots,\lambda_r\ge 0, \quad \lambda_1+\ldots+\lambda_r=1
\label{19-7-2}
\shortintertext{such that}
\sum_k\lambda_k \nabla \log (-V/f_k)(\bar{x})=0;
\label{19-7-3}
\end{gather}
let $\bK$ be the set of $\boldlambda$, satisfying (\ref{19-7-2}).

Further, let us assume that\begin{phantomequation}\label{19-7-4}\end{phantomequation}
\begin{equation}
\rank \bigl\{\nabla \log (-V/f_1)(\bar{x}),\ldots, \nabla \log (-V/f_r)(\bar{x})\bigr\}\ge m
\tag*{$\textup{(\ref*{19-7-4})}_m$}\label{19-7-4-m}
\end{equation}
with $1\le m\le r-1$.

Note that if $m=r-1$ (case of the \emph{simple degeneration}\index{degeneration!simple} then $\boldlambda=(\lambda_1,\ldots,\lambda_r)\in \bK$ in (\ref{19-7-3}) is defined uniquely; more generally, $\boldlambda\in \bK$ is defined uniquely modulo elements of $(r-m-1)$-dimensional subspace $\bL$ and
\begin{equation}
|\sum_k\lambda_k \nabla \log (-V/f_k)(\bar{x})|\asymp \nu \Def 
\dist (\boldlambda, \bL')\qquad\forall \boldlambda\in \bK 
\label{19-7-5}
\end{equation}
where affine subspace $\bL'$ is defined by (\ref{19-7-3}) and $\lambda_1+\ldots+\lambda_r=1$. 

We consider only the smooth case i.e. assume that
\begin{equation}
g^{jk}, V_j, V\in \sC^K
\label{19-7-6}
\end{equation}
with sufficiently large $K=K(d)$. Then (\ref{19-7-5}) remains true in  $\epsilon \nu$-vicinity of $x$.

\section{Weak magnetic case}
\label{sect-19-7-2}

Using standard reduction we reach 
\begin{equation}
\cA =\mu h \sum_j f_j(x'',\mu^{-1}hD'') \bigl(x_j^2+ D_j^2\bigr) + O(\mu^{-l})
\label{19-7-7}
\end{equation}
with $l=2,1$ if there are no $3$-rd order resonances and in the general case respectively. 

We leave to the reader the following

\begin{Problem}\label{Problem-19-7-1}
(i) Prove by our standard method that one can localize  with respect to 
$\mu h\bigl(x_j^2+ D_j^2\bigr)$ modulo 
$O(\mu^{1-l} T + C\mu h|\log h|)$. 

\medskip\noindent
(ii) Prove by the method of Section~\ref{sect-19-5} that one can localize with respect to $\mu h\bigl(x_j^2+ D_j^2\bigr)$ modulo  $O(\mu^{1-l} T + C\mu h)$.
\end{Problem}

Assume that localization is done in vicinity of $\boldtau$. Consider symbol
\begin{equation}
b(x'',\xi'';\boldtau) = \sum_j f_j(x'',\xi'')\tau_j.
\label{19-7-8}
\end{equation}

Consider point  $(\bar{x''},\bar{\xi}'')$. Then if in its $\gamma$-vicinity 
\begin{gather}
|\nabla b |\asymp \nu\gamma^{-1}, \quad |b|\lesssim \nu, \qquad \nu = \gamma^2, 
\label{19-7-9}\\
\gamma \ge C_0\max(\mu^{-1} , (\mu^{-1}h|\log h|)^{\frac{1}{2}} ); 
\label{19-7-10}
\end{gather}
then the shift for time $T$ is $\asymp \mu^{-1}\nu T$ and it is observable as 
$\mu^{-1}\nu T \times \gamma \ge C\mu^{-1}h|\log h|$ and one can take $T= \epsilon_0\mu^{-1}$ as long as 
\begin{equation}
\gamma \ge \bar{\gamma}\Def
C_0 \max\bigl((\mu h|\log h|)^{\frac{1}{2}} ,\mu^{-l}\bigr) ;
\label{19-7-11}
\end{equation}
here the last requirement ($\gamma \ge C_0\mu^{-l}$) is needed to counter the irreducible terms of order $l+2$ but it contributes only 
$C\mu h^{1-d}\times \mu^{-l}=C\mu^{1-l}h^{1-d}$ into the final remainder estimate.  Then one can take 
\begin{equation}
T^*\asymp \mu^{l-1}\gamma
\label{19-7-12}
\end{equation}
and $T_* =\epsilon_0\mu^{-1}$ and then contribution of $\gamma$-vicinity to Tauberian remainder does not exceed $C\gamma^{2d} h^{1-d} T^{*\,-1}\asymp C\gamma^{2d-1} \mu^{1-l} h^{1-d}$ and then contribution of zone $\{\gamma \ge \bar{\gamma}\}$ to such remainder does not exceed 
\begin{equation}
C\mu^{1-l}h^{1-d} \int_{\{\gamma \ge \bar{\gamma}\}} \gamma^{-1}\,dx''d\xi''d\boldtau,
\label{19-7-13}
\end{equation}
while contribution of zone $\{\gamma \le \bar{\gamma}\}$ does not exceed 
\begin{equation}
C\mu h^{1-d}\mes \bigl(\{\gamma \le \bar{\gamma}\}\bigr)
\label{19-7-14}
\end{equation}
and the total remainder does not exceed (\ref{19-7-13})+(\ref{19-7-14}):
\begin{equation}
C\mu^{1-l}h^{1-d} \int_{\gamma\ge \bar{\gamma}}  \gamma^{-1}\,d M(\gamma) + 
C\mu h^{1-d}M(\bar{\gamma})
\label{19-7-15}
\end{equation}
with
\begin{multline}
M (\gamma )\Def \mes_{3r} \bigl(\{(x'',\xi'',\boldtau): |\nabla_{x'',\xi''}  b(x'',\xi'',\boldtau) | \asymp \gamma\}\bigr)= \\
\int M(x'',\xi'';\gamma)\,dx''d\xi'',
\label{19-7-16}
\end{multline}
\begin{equation}
M(x'',\xi'';\gamma) =\mes _r \bigl(\{(\boldtau): |\nabla_{x'',\xi''}  b(x'',\xi'',\boldtau)|\asymp \gamma\}\bigr).
\label{19-7-17}
\end{equation}

Note that 
\begin{claim}\label{19-7-18}
Under assumption \ref{19-7-4-m} $M(x'',\xi'';\gamma)= O(\gamma^m)$.
\end{claim}

\begin{Problem}\label{Problem-19-7-2}
Justify the following: We have $(3r-1)$-dimensional variable $(x,\boldlambda)$ (due to restriction $\lambda_1+\ldots+\lambda_r=1$) and (\ref{19-7-3}) imposes $2r$ restrictions; thus it happens generically on a variety of dimension $(r-1)$. However \ref{19-7-4-m}  with $m=(r-2)$ imposes $(2r-1)$ further restrictions more and thus never happens generically.
\end{Problem}

So let us assume that assumption $\textup{(\ref{19-7-4})}_{r-1}$ is fulfilled.  Then there exists a single  $\boldlambda =\boldlambda (x)$ minimizing 
\begin{equation}
|\sum_k \lambda_k \nabla \log (-V/f_k)|^2.
\label{19-7-19}
\end{equation}
Let us assume that 
\begin{phantomequation}\label{19-7-20}\end{phantomequation}
\begin{multline}
\nabla (\sum \lambda_i \log (-V/f_i)) =0 \implies\\
\rank \bigl\{ \sum _i \partial _{j} (\lambda_i  \partial_k \log (-V/f_i)) \bigr\}_{j,k=1,\ldots,2r}=n.
\tag*{$\textup{(\ref*{19-7-20})}_n$}\label{19-7-20-n}
\end{multline}
Note that 
\begin{claim}\label{19-7-21}
Under assumptions $\textup{(\ref{19-7-4})}_{r-1}$ and \ref{19-7-20-n} $M(\gamma)= O(\gamma^{m+n})$.
\end{claim}

\begin{Problem}\label{Problem-19-7-3}
Prove that generically (\ref{19-7-9}) is fulfilled with $n=r+1$ and therefore $M(\gamma)=\gamma^{2r}$.
\end{Problem}

Therefore under assumption $\textup{(\ref{19-7-4})}_{r-1}$ we conclude that
\begin{equation}
\R^\T \le C\mu^{-1}h^{1-d} + C \mu (\mu h |\log h|)^{\frac{1}{2}(r-1)} h^{1-d}
\label{19-7-22}
\end{equation}
and under assumptions $\textup{(\ref{19-7-4})}_{r-1}$ and $\textup{(\ref{19-7-20})}_{r+1}$
\begin{equation}
\R^\T \le 
C\mu^{1-l}h^{1-d} + C \mu (\mu h|\log h|)^r h^{1-d}.
\label{19-7-23}
\end{equation}

Further the same estimates hold for $\R^\W_{(\infty)}$.

\begin{remark}\label{rem-19-7-4}
Actually summation with respect to the partition as $r=2$ and only the first assumption is fulfilled returns $\mu^{1-l}h^{1-d}\log \mu$ instead of the first term.

However, as $m=1$ we can consider propagation in the direction of $\nu$ increasing and we would be able to take $T^*\asymp \mu \nu^{1-\delta}$ if there are no $3$-rd order resonances and $T^*\asymp  \nu^{1-\delta}$ otherwise (we used similar arguments multiple times, for example in the analysis of the exterior zone for $3\D$ magnetic Schr\"odinger) and then we would arrive to (\ref{19-7-22}) again.
\end{remark}

\begin{remark}\label{rem-19-7-5}
Using rescaling technique one can get rid off logarithmic factors arriving to
\begin{gather}
\R^\T \le C\mu^{-1}h^{1-d} + C  (\mu h )^{\frac{1}{2}(r+1)} h^{-d}
\tag*{$\textup{(\ref*{19-7-22})}^*$}\label{19-7-22-*}\\
\shortintertext{and} 
\R^\T \le 
C\mu^{1-l}h^{1-d} + C  (\mu h)^{r+1} h^{-d}.
\tag*{$\textup{(\ref*{19-7-23})}^*$}\label{19-7-23-*}
\end{gather}
respectively.
\end{remark}

Therefore we arrive to the same estimates for $\R^\W_{(\infty)}$. 

\begin{theorem}\label{thm-19-7-6}
Let assumptions \textup{(\ref{19-1-4})}--\textup{(\ref{19-1-6})}, \textup{(\ref{19-7-6})} and \textup{(\ref{19-2-38})}  be fulfilled.

Then under assumption $\textup{(\ref{19-7-4})}_{r-1}$
\begin{equation}
\R^\MW \le C\mu^{-1}h^{1-d} + C  (\mu h )^{\frac{1}{2}(r+1)} h^{-d}
\label{19-7-22-24}\\
\end{equation}
and under assumptions $\textup{(\ref{19-7-4})}_{r-1}$ and $\textup{(\ref{19-7-20})}_{r+1}$
\begin{equation}
\R^\MW \le 
C\mu^{1-l}h^{1-d} + C  (\mu h)^{r+1} h^{-d}.
\label{19-7-25}
\end{equation}
respectively.
\end{theorem}

\begin{proof}
We leave to the reader the transition from $\R^\W_{(\infty)}$ to $\R^\MW$.
\end{proof}

\section{Intermediate magnetic field}
\label{sect-19-7-3}

Let us assume now that \underline{either}  $\textup{(\ref{19-7-4})}_{r-1}$ is fulfilled and $(\mu h)^{(r+1)/2} \ge C\mu^{-1}h$ which is equivalent to 
\begin{equation}
 C_0 h^{-(r-1)/(r+3)}\le \mu \le \epsilon h^{-1}
\label{19-7-26}
\end{equation} 
\underline{or}  both $\textup{(\ref{19-7-4})}_{r-1}$ and $\textup{(\ref{19-7-20})}_{r+1}$ are fulfilled and $(\mu h)^{(r+1)} \ge C\mu^{-1}h$ 
which is equivalent to 
\begin{equation}
 C_0 h^{-r/(r+2)}\le \mu \le \epsilon h^{-1}.
\label{19-7-27}
\end{equation} 

\subsection{Tauberian remainder: the pilot-model}
\label{sect-19-7-3-1}

Assume temporarily that there are no resonances of order less than $L$ with large enough $L$ except the trivial ones i.e. with the irreducible terms
\begin{equation}
\mu^{2-L} b_\beta (x'',\mu^{-1}hD'') (h^2D_1^2 + \mu^2 x_1^2)^{\beta_1}\cdots (h^2D_r^2 + \mu^2 x_r^2)^{\beta_r};
\label{19-7-28}
\end{equation} 
then in view of (\ref{19-7-26}) or (\ref{19-7-27}) we can skip all the terms corresponding to resonances of order $\ge L$ (with a very small error) and we get a family of $r$-dimensional $\mu^{-1}h$-pseudo-differential operators and we can apply our standard theory to them. 

Under the non-degeneracy assumption $\textup{(\ref{19-7-20})}_{r+1}$ all of these operators satisfy non-degeneracy assumption ``no degenerate critical points'' and therefore the individual remainder estimate would be $O\bigl( (\mu^{-1}h)^{1-r}\bigr)$; then multiplying by the number of operators 
$\asymp (\mu h)^{-r}$ we get $O(\mu^{-1}h^{1-d})$.

\medskip
Meanwhile without assumption $\textup{(\ref{19-7-20})}_{r+1}$ these operators could be highly degenerate and the individual remainder estimate would be $O\bigl( (\mu^{-1}h)^{-r}\bigr)$ and thus we do not expect a better total remainder estimate. 

Let us prove that this is also a total remainder estimate. According to Subsection~\ref{book_new-sect-5-1-4} of \cite{futurebook} we can introduce an admissible function 
$\ell =\ell(z,\boldlambda)$ with $z=(x'',\xi'')$ such that in $B(z,\ell(z))$ 
\begin{gather}
\ell (z,\boldlambda) \asymp \sum_{\alpha:|\alpha |\le m}
|\nabla^\alpha b(z,\boldlambda)|^{\frac{1} {m+1-|\alpha|}} 
\label{19-7-29}\\
\intertext{and then we redefine}
\ell (z,\boldlambda) \Def \max\bigl(\ell (z,\boldlambda) , \bar{\ell}\bigr), \qquad \bar{\ell}= C_0(\mu^{-1}h)^{\frac{1}{2}(1-\delta)}
\label{19-7-30}
\end{gather}
with arbitrarily large $m$ and arbitrarily small $\delta>0$. 

Then again according to the same Subsection~\ref{book_new-sect-5-1-4} of \cite{futurebook}  the Tauberian remainder for individual $\boldlambda$ does not exceed
\begin{equation}
C\mu^{r-1}h^{1-r-\delta'} \int \ell (z, \boldlambda)^{-2} \,dz + C\mu^r h^{-r} 
\int_{\ell (z, \boldlambda) \le \bar{\ell}}\, dz
\label{19-7-31}
\end{equation}
and therefore the total Tauberian remainder does not exceed
\begin{equation}
C\mu^{r-1}h^{1-r-\delta'} \sum_{\boldlambda}\int \ell (z, \boldlambda)^{-2} \,dz + C\mu^r h^{-r} \sum_{\boldlambda}
\int_{\ell (z, \boldlambda) \le \bar{\ell}}\, dz
\label{19-7-32}
\end{equation}
where $\boldlambda$ runs through the lattice. 

Note that $\gamma \Def \ell_1 \le C\ell^m$ (where $\ell=\ell_m$). Therefore  $\ell \asymp\bar{\ell}$ implies that $\gamma \lesssim \mu h$ provided \begin{equation}
m (r+1)>6.
\label{19-7-33}
\end{equation}
But then for each $z$ there exist no more than $C_0$ values of $\boldlambda$ such that $\ell(z,\boldlambda)\le \bar{\ell}$ and therefore the second term in (\ref{19-7-32}) is $O(\mu^rh^{-r})$. This statement remains valid as we redefine $\bar{\ell}= (\mu h)^{1/(m-1)}$. 

Meanwhile for any given $z$ the number of $\boldlambda$ such that $\ell(z,\boldlambda)\le \ell$ does not exceed 
$C_0(\mu h)^{-r}\gamma^r + C_0\lesssim C_0(\mu h)^{-r}\ell^{(m-1)r}+C_0$ and therefore the first term in (\ref{19-7-32}) does not exceed 
\begin{equation}
C\mu^{-1}h^{1-d-\delta'}\int \ell^{m r-2} \,dz + C\mu^r h^{-r}.
\label{19-7-34}
\end{equation}
Since $m r> 2$ this term is $O(\mu^{-1}h^{1-d-\delta'})$ which is almost exactly what we need. And it is less than $C\mu^r h^{-r}$ unless $\mu\le C h^{-(r-1)/(r+1)-\delta''}$ so we need to consider only this case.

To cover this case we need to make $\delta'=0$. Then we need to consider only elements of $\ell_m$-partition  with $\ell=\ell_m\ge h^{\delta}$. However if we consider elements of $\ell_{m-1}$-partition with $\ell_{m-1}\le h^{\delta_1}$ (but $\ell_m\ge h^\delta$), their contribution to the remainder will be properly estimated as well, so we need  to consider only elements of $\ell_{m-1}$-partition  with $\ell_{m-1}\ge h^{\delta_1}$; continuing this process  we conclude that we need to consider only balls with $\ell_1\ge h^{\delta_{m-1}}$. But this is exactly $\gamma$ and we can apply weak magnetic field approach.

\subsection{Tauberian estimates: the general case}
\label{sect-19-7-3-2}

As far as Tauberian estimate is concerned everything remains the same as 
$\ell^{m+1} \ge \mu^{-L}$ where $\mu^{-L}$ is the magnitude of the highest irreducible ``extra'' terms. Further, let 
\begin{equation}
\ell (z,\boldlambda) \asymp \sum_{\alpha:1\le |\alpha |\le m}
|\nabla^\alpha b(z,\boldlambda)|^{\frac{1} {m+1-|\alpha|}}  +
\min _{s:|s|\le C_0\mu^{-L}}  |b(z,\boldlambda)+s|^{\frac{1} {m+1}}
\tag*{$\textup{(\ref*{19-7-29})}'$}\label{19-7-29-'}
\end{equation}
(in contrast to (\ref{19-7-29}) term with $\alpha=0$ is modified).  

However this would be absorbed by the case $\ell ^m \le \mu h$ provided $\mu^{-L}\le \mu h$. 

\medskip\noindent
(a) As $L=1$ we are are looking at the remainder estimate $O(h^{1-d})$ i.e. we need to consider the case $\mu \ge h^{-(r-1)/(r+1)}$ and therefore 
$\mu \ge h^{-\frac{1}{2}}$ as $r\ge 3$ and $\mu \ge h^{-\frac{1}{3}}$ as $r=2$. However condition of absorption is $\mu ^{-1}\le \mu h$; so as $r\ge 3$ everything is fine: we arrive to the remainder estimate $O\bigl(h^{1-2r}+ \mu^rh^{-r}\bigr)$.

As $r=2$ we also arrive to this remainder estimate  unless $h^{-\frac{1}{3}}\le \mu \le h^{-\frac{1}{2}}$ in which case the remainder estimate is 
\begin{multline}
Ch^{-3}+ C\min \bigl( (\mu h)^{\frac{3}{2}},\mu^{-2}\bigr)h^{-4}=\\
Ch^{-3}+ Ch^{-4}\left\{\begin{aligned}
&(\mu h)^{\frac{3}{2}}\qquad &&\text{as\ \ } \mu \le h^{-\frac{3}{7}},\\
&\mu^{-2}\qquad &&\text{as\ \ } \mu \ge h^{-\frac{3}{7}}.
\end{aligned}\right.
\label{19-7-35}
\end{multline}
(b) As $L=2$ we are are looking at remainder estimate $O(\mu^{-1}h^{1-d})$ and therefore we need to consider case of (\ref{19-7-26}) i.e. $\mu \ge h^{-\frac{1}{3}}$ as $r\ge 3$ and $\mu \ge h^{-\frac{1}{5}}$ as $r=2$. However condition of absorption is $\mu ^{-2}\le \mu h$; so as $r\ge 3$ everything is fine: we arrive to the remainder estimate $O\bigl(\mu^{-1}h^{1-2r}+ \mu^rh^{-r}\bigr)$.

As $r=2$ the same arguments work for $L=4$. As $L=2$ we also get this remainder estimate unless $h^{-\frac{1}{5}}\le \mu \le h^{-\frac{1}{3}}$ when we get 
\begin{multline}
C\mu^{-1}h^{-3}+ C\min \bigl( (\mu h)^{\frac{3}{2}},\mu^{-4}\bigr)h^{-4}=\\
C\mu^{-1}h^{-3}+ Ch^{-4}\left\{\begin{aligned}
&(\mu h)^{\frac{3}{2}}\qquad &&\text{as\ \ } \mu \le h^{-\frac{3}{11}},\\
&\mu^{-2}\qquad &&\text{as\ \ } \mu \ge h^{-\frac{3}{11}}.
\end{aligned}\right.
\label{19-7-36}
\end{multline}
(c) As $r=2$, $L=3$ contribution of perturbation is $O(\mu^{-6}h^{-3})=O(\mu^{-1}h^{-1})$ as $\mu \ge h^{-\frac{1}{5}}$ and therefore we get a proper remainder estimate.

Therefore we arrive to

\begin{proposition}\label{prop-19-7-7}
Let $\mu \le \epsilon h^{-1}$ and assumption $\textup{(\ref{19-7-4})}_{r-1}$ be fulfilled. Then

\medskip\noindent
(i) As $L=1$ 
\begin{equation}
\R^\T \le Ch^{1-2r}+ C\mu^r h^{-r}
\label{19-7-37}
\end{equation} 
unless $r=2$, $h^{-\frac{1}{3}}\le \mu \le h^{-\frac{1}{2}}$ when $\R^\T$ does not exceed expression \textup{(\ref{19-7-35})};

\medskip\noindent
(ii) As $L=2$ 
\begin{equation}
\R^\T \le C\mu^{-1}h^{1-2r}+ C\mu^r h^{-r}
\label{19-7-38}
\end{equation} 
unless $r=2$, $h^{-\frac{1}{5}}\le \mu \le h^{-\frac{1}{3}}$ when $\R^\T$ does not exceed expression \textup{(\ref{19-7-36})}; furthermore, as $r=2$, $L=3$ estimate \textup{(\ref{19-7-38})} holds.
\end{proposition}

Assume now that extra assumption $\textup{(\ref{19-7-20})}_{r+1}$ is fulfilled. Then 

\medskip\noindent
(d) As $L=1$ we are looking for the remainder $O(h^{1-d})$ i.e. we need to consider a case $(\mu h)^{r+1}\ge h$ i.e. $\mu \ge h^{-\frac{2}{3}}$ and then the contribution to the remainder of the degeneration will be $O(\bigl((\mu^{-1/2})^{2r} h^{-2r}\bigr)=O(h^{1-2r})$.

\medskip\noindent
(e) As $L=2$ we are looking for the remainder $O(\mu^{-1}h^{1-d})$ i.e. we need to consider a case of (\ref{19-7-27}) i.e. $\mu\ge h^{-\frac{1}{2}}$ and then contribution to the remainder of the degeneration will be $O(\bigl((\mu^{-2/2})^{2r} h^{-2r}\bigr)=O(\mu^{-1}h^{1-2r})$.

\begin{proposition}\label{prop-19-7-8}
Let $\mu \le \epsilon h^{-1}$ and assumptions $\textup{(\ref{19-7-4})}_{r-1}$ and $\textup{(\ref{19-7-20})}_{r+1}$ be fulfilled. Then

\medskip\noindent
(i) As   $L=1$ 
\begin{equation}
\R^\T \le Ch^{1-d};
\label{19-7-39}
\end{equation} 
(ii) As   $L=2$ 
\begin{equation}
\R^\T \le C\mu^{-1}h^{1-d}.
\label{19-7-40}
\end{equation} 
\end{proposition}

\subsection{Calculations: the pilot-model}
\label{sect-19-7-3-3}

In the pilot-model case we get the Tauberian expression $h^{-d}\cN^\T$ with 
\begin{multline}
\cN^\T\Def \sum _\alpha \mu^r h^{r} \uptheta \bigl(1 - b_0(x, \alpha) - \mu^{-2}b_1 (x,\alpha) - \ldots \bigr) \times\\
\bigl(\psi_0 (x) + \mu^{-2} \psi_1 (x,\alpha) + \ldots\bigr)f_1\cdots f_r \sqrt{g}\,dx
\label{19-7-41}
\end{multline}
and if we remove all powers of $\mu^{-2}$ then we get the final answer.

As $\mu \ge h^{-\frac{1}{3}}$ we can following our approach  skip instantly  $O(\mu^{-4})$; we can also remove $O(\mu^{-2})$ but add a correction term instead which can be written as $\kappa \mu^{-2} $ and prove that an error does not exceed $C\mu^{-1}h$ provided $r\ge 3$; then we conclude that $\kappa =0$ and therefore as $r\ge 3$ we arrive modulo $O(\mu^{-1}h)$ to $\cN^\MW$. Therefore for the pilot-model operator with $r\ge 3$ we arrive to the desired formula.

For $r=2$ situation is more complicated: we can assume only that $\mu \ge h^{-\frac{1}{5}}$ and then we can skip instantly only $O(\mu^{-6})$ (rather than $O(\mu^{-4})$); we can also  remove $O(\mu^{-4})$ but add a correction term instead which will be given  by (\ref{19-4-93}). 

For $r=2$ and $\mu\ge h^{-\frac{1}{3}}$ we can remove $O(\mu^{-2})$ and add a correction term instead which will be $0$ but only modulo $O(\mu^{-2}(\mu h)^{\frac{1}{2}} h^{-4})= O(\mu^{-\frac{3}{2}}h^{-\frac{7}{2}})$.

On the other hand, under assumption $\textup{(\ref{19-7-20})}_{r+1}$ we can deal with $O(\mu^{-2})$ terms with an error $O(\mu^{-2}(\mu h)^r h^{-2d})$ which is $O(\mu^{-1}h^{1-2d})$ even as $r=2$. We leave all details to the reader.

Thus we arrive to the following theorem:

\begin{theorem}\label{thm-19-7-9}
Let assumptions \textup{(\ref{19-1-4})}--\textup{(\ref{19-1-6})}, \textup{(\ref{19-7-6})} and \textup{(\ref{19-2-38})}  be fulfilled. Let there be no resonances of the $2$-nd and $3$-rd order, and no non-trivial resonances of order less than $L$ (with large enough $L$) and let assumption $\textup{(\ref{19-7-4})}_{r-1}$ be fulfilled.  Then 

\medskip\noindent
(i) As $r\ge 3$, $\mu \le \epsilon (h|\log h|)^{-1}$ 
\begin{equation}
\R^\MW \le C\mu^{-1}h^{1-2r} + C\mu ^r h^{-r};
\label{19-7-42}
\end{equation} 
(ii) As as  $r=2$, $h^{-\frac{1}{5}}\le \mu \le \epsilon (h|\log h|)^{-1}$
\begin{multline}
\R_1^\MW \Def\\
|\int \Bigl( \tilde{e}(x,x,0) - h^{-d}\cN^\MW (x,0)- 
h^{-d}\cN^\MW_{1\corr }(x,0)\Bigr)\psi (x)\, dx |\le\\[3pt]
C\mu ^{-1}h^{1-d} + C\mu ^{2} h^{-2}
\label{19-7-43}
\end{multline}
where $\cN^\MW_{1\corr }(x,0)\Bigr)$ is defined by \textup{(\ref{19-4-93})} and $\cA_\cT=\cA_0 + \mu^{-2}\cA_2$; in particular, as $r=2$, $h^{-\frac{1}{3}}\le \mu \le \epsilon (h|\log h|)^{-1}$
this correction term does not exceed $C\mu^{-\frac{3}{2}}h^{-\frac{7}{2}}$  and therefore one can skip it as $\mu \ge h^{-\frac{3}{7}}$;

\medskip\noindent
(iii) On the other hand, under assumption $\textup{(\ref{19-7-20})}_{r+1}$ 
\begin{equation}
\R^\MW \le C\mu^{-1}h^{1-2r}.
\label{19-7-44}
\end{equation} 
\end{theorem}

\subsection{Calculations: the general case}
\label{sect-19-7-3-4}

Consider now the general case.  Assume first that there are no $3$-rd order resonances. Then both arguments of the proof and statements of theorem~\ref{thm-19-7-9} obviously remain true with the singular exception of (ii) where now we need to take
\begin{equation}
\cA_\cT=\cA_0 + \mu^{-2}\cA_2+\mu^{-3}\cA_3
\label{19-7-45}
\end{equation}
and $\cA_2$ may contain also a non-diagonal part (if there are non-trivial $4$-th order resonances) and $\cA_3$ is purely a non-diagonal term (if there are $5$-th order resonances). 

Now we would like to get rid of all non-diagonal terms in $\mu^{-2}\cA_2$ and $\mu^{-3}\cA_3$. It works as long as $\mu \ge h^{-\frac{1}{3}}$ and $\mu \ge h^{-\frac{1}{4}}$ respectively. Then we will have another correction term. However it will be a bit smaller than $\mu^{-L}(\mu h)^{\frac{1}{2}}h^{-4}$. Namely if we take in account only the leading term in $\psi$, i.e. $\psi_0$, we get a non-diagonal term with the trace $0$ in the second approximation term and the error will be $O(\mu^{-\frac{3}{2}L}h^{-4})$ which is a bit better than  
$O(\mu^{-L}(\mu h)^{\frac{1}{2}}h^{-4})$ as $\mu^{L+1}h\ge 1$.

Modulo this term there will be term coming from $\psi_1\mu^{-1}$ and our standard arguments yield that it will not exceed $O(\mu^{-L-1}h^{-4})$ which is exactly what we got before as $L=2$. However as $L=3$ one can see easily that the first perturbation term would be out-of-diagonal by more than $1$ and we need ``interact'' it with $\psi_2 \mu^{-2}$ to get non-$0$ trace and this term will be  $O(\mu^{-L-2}h^{-4})$. 

Therefore we arrive to

\begin{theorem}\label{thm-19-7-10}
Let assumptions \textup{(\ref{19-1-4})}--\textup{(\ref{19-1-6})}, \textup{(\ref{19-7-6})} and \textup{(\ref{19-2-38})}  be fulfilled.
Let there be no resonances of the $2$-nd and $3$-rd order and let assumption $\textup{(\ref{19-7-4})}_{r-1}$ be fulfilled.  Then 

\medskip\noindent
(i) Statements (i), (ii) of theorem~\ref{thm-19-7-9} hold; statement (ii) remains true without any modifications if there are no non-trivial $4$-th and no $5$-th order resonances;

\medskip\noindent
(ii) In statement (ii) of theorem~\ref{thm-19-7-9} one should take $\cA_\cT$ in the form \textup{(\ref{19-7-45})} but the conclusion about magnitude of the correction term remains true; 

\medskip\noindent
(iii) If $r=2$ and there are non-trivial $4$-th order resonances and we take $\cA_\cT= \cA_0 +\mu^{-2}\cA_2$ and remove from $\cA_2$ non-diagonal terms, estimate 
\begin{equation}
\R_1^\MW \le C\mu^{-1}h^{-3} + C\mu^2 h^{-2} + C\mu^{-3}h^{-4}
\label{19-7-46}
\end{equation}
holds as $\mu \ge h^{-\frac{1}{3}}$; in particular \textup{(\ref{19-7-43})} holds as $\mu \ge h^{-\frac{2}{5}}$;

\medskip\noindent
(iv) If $r=2$ and there are no non-trivial $4$-th order resonances but there are  $5$-th order resonances and we take $\cA_\cT= \cA_0 +\mu^{-2}\cA_2$, estimate 
\begin{equation}
\R_1^\MW \le C\mu^{-1}h^{-3} + C\mu^2 h^{-2} + C\mu^{-\frac{9}{2}}h^{-4}
\label{19-7-47}
\end{equation}
holds as $\mu \ge h^{-\frac{1}{4}}$; in particular \textup{(\ref{19-7-43})} holds as $\mu \ge h^{-\frac{2}{7}}$.
\end{theorem}

Assume now that there are $3$-rd order resonances. Then for $r\ge 3$ the weak magnetic field approach gives remainder estimate $O(h^{1-d})$ as $\mu \le h^{-\frac{1}{2}}$ and an error when we remove $O(\mu^{-1})$ non-diagonal term is $O(\mu^{-2}h^{-d})=O(h^{1-d})$ as $\mu \ge h^{-\frac{1}{2}}$. Therefore we are done.

As $r=2$ and both assumptions  $\textup{(\ref{19-7-4})}_1$ and $\textup{(\ref{19-7-20})}_3$ are fulfilled, the weak magnetic field approach gives remainder estimate $O(h^{-3})$ as $(\mu h)^3 \le h$ i.e. as $\mu \le h^{-\frac{2}{3}}$ and  an error when we remove $O(\mu^{-1})$ non-diagonal term is $O(\mu^{-2}h^{-4})=O(h^{-3})$  as $\mu \ge h^{-\frac{2}{3}}$. Therefore we are done as well.

Consider case $r=2$ and only  assumption $\textup{(\ref{19-7-4})}_1$ is fulfilled. Then the weak magnetic field approach gives remainder estimate $O(h^{-3})$ as $\mu \le h^{-\frac{1}{4}}$. Then we have a proper remainder estimate with the correction term defined by  
\begin{equation}
\cA_\cT=\cA_0 + \mu^{-1}\cA_1 +\mu^{-2}\cA_2+\mu^{-3}\cA_3;
\label{19-7-48}
\end{equation}
however in virtue of our previous arguments skipping $O(\mu^{-3})$ term would lead to an error $O(\mu^{-4}(\mu h)^{\frac{1}{2}}h^{-4})= O(h^{-3})$ then.  Therefore we arrive to

\begin{theorem}\label{thm-19-7-11}
Let assumptions \textup{(\ref{19-1-4})}--\textup{(\ref{19-1-6})}, \textup{(\ref{19-7-6})} and \textup{(\ref{19-2-38})}  be fulfilled. Let there be no resonances of the $2$-nd  order and let assumption $\textup{(\ref{19-7-4})}_{r-1}$ be fulfilled.  Then 

\medskip\noindent
(i) As $r\ge 3$, $\mu \le \epsilon (h|\log h|)^{-1}$ 
\begin{equation}
\R^\MW \le Ch^{1-2r} + C\mu ^r h^{-r};
\label{19-7-49}
\end{equation} 
(ii) As as  $r=2$, $h^{-\frac{1}{5}}\le \mu \le \epsilon (h|\log h|)^{-1}$
\begin{multline}
\R_1^\MW \Def\\
|\int \Bigl( \tilde{e}(x,x,0) - h^{-d}\cN^\MW (x,0)- 
h^{-d}\cN^\MW_{1\corr }(x,0)\Bigr)\psi (x)\, dx |\le\\[3pt]
Ch^{-3} + C\mu ^{2} h^{-2}
\label{19-7-50}
\end{multline}
where $\cN^\MW_{1\corr }(x,0)$ is defined by \textup{(\ref{19-4-93})} and $\cA_\cT$ is defined by \textup{(\ref{19-7-48})} without the last term;  in particular, as $r=2$, $h^{-\frac{1}{2}}\le \mu \le \epsilon (h|\log h|)^{-1}$
this correction term does not exceed $C\mu^{-\frac{3}{2}}h^{-4}$  and therefore one can skip it as $\mu \ge h^{-\frac{4}{7}}$;

\medskip\noindent
(iii) On the other hand, under assumption $\textup{(\ref{19-7-20})}_{r+1}$ 
\begin{equation}
\R^\MW \le Ch^{1-2r}.
\label{19-7-51}
\end{equation} 
\end{theorem}

\begin{remark}\label{rem-19-7-12}
Further, skipping $O(\mu^{-2})$ term in $\cA_\cT$ leads to an error $O(\mu^{-2}(\mu h)^{\frac{1}{2}}h^{-4})$ which is $O(h^{-3})$ as $\mu \ge h^{-\frac{1}{3}}$. 

Furthermore, skipping $O(\mu^{-1})$ non-diagonal term leads to an error $O(\mu^{-\frac{3}{2}}h^{-4})=O(h^{-3})$ as $\mu \ge h^{-\frac{2}{3}}$.
\end{remark}

\section{Strong magnetic field}
\label{sect-19-7-4}

The case of the strong magnetic field is easy:

\begin{theorem}\label{thm-19-7-13}
Let assumptions \textup{(\ref{19-1-4})}--\textup{(\ref{19-1-6})}, \textup{(\ref{19-7-6})} and \textup{(\ref{19-2-38})}  be fulfilled. Let there be no resonances of the $2$-nd order, and let assumption $\textup{(\ref{19-7-4})}_{r-1}$ be fulfilled. Let $\epsilon (h|\log h|)^{-1}\le \mu \le \epsilon h^{-1}$.

\medskip\noindent
(i) Assume that there no resonances of the $3$-rd order as well. Then  estimate \textup{(\ref{19-7-42})} holds and under assumption $\textup{(\ref{19-7-20})}_{r+1}$  estimate \textup{(\ref{19-7-44})} holds.

\medskip\noindent
(ii) In the general case   estimate \textup{(\ref{19-7-49})} holds and under assumption $\textup{(\ref{19-7-20})}_{r+1}$  estimate \textup{(\ref{19-7-50})} holds.
\end{theorem}

\begin{proof}
An easy proof using arguments of Section~\ref{sect-19-5} is left to the reader.
{\ }
\end{proof}

\section{Conclusion to Section}
\label{sect-19-7-5}
There is no theory of very strong and superstrong magnetic field separate from theory of Section~\ref{sect-19-6} 

\begin{Problem}\label{Problem-19-7-14}
Prove similar results under assumptions \ref{19-7-4-m} and \ref{19-7-20-n} with $0\le m<r-1$  and $1\le m+n \le 2r$.
\end{Problem}

\section{Final remarks}
\label{sect-19-7-6}
%%WORK
\subsection{Vanishing $V$}
\label{sect-19-7-6-1}

\begin{remark}\label{rem-19-7-15}
(i) In Sections~\ref{sect-19-2}--\ref{sect-19-5} we divided by $V$ in the microhyperbolicity or $\fN$-microhyperbolicity assumptions. However, it is not a problem: we could assume in advance that
\begin{equation}
|V|+|\nabla V|\ge \epsilon_0
\label{19-7-52}
\end{equation}
and then the microhyperbolicity or $\fN$-microhyperbolicity conditions  will be needed in zone $\{x:|V(x)|\ge \epsilon'\}$ only.

\medskip\noindent
(ii) We can assume instead that
\begin{equation}
|V|+|\nabla V|\le \epsilon_0 \implies |\det \Hess V |\ge \epsilon_0.
\label{19-7-53}
\end{equation}
Really then due to a rescaling technique we can ensure microhyperbolicity condition (except a small zone) after rescaling as $|V| \le \epsilon '$
and then the microhyperbolicity or $\fN$-microhyperbolicity conditions  will be needed in zone $\{x:|V(x)|\ge \epsilon'\}$ only.
\end{remark}

We leave all the details to the reader the following:

\begin{Problem}\label{Problem-19-7-16}
(i) Using rescaling technique like in Subsection~~\ref{book_new-sect-18-9-5} of \cite{futurebook} get rid off assumption $|V|\ge \epsilon_0$ in the framework of Sections~\ref{sect-19-2}, ~\ref{sect-19-4}, and~\ref{sect-19-5};

\medskip\noindent
(ii) Using rescaling technique get rid off assumption $|V|\ge \epsilon_0$ in the framework of Sections~\ref{sect-19-2}, ~\ref{sect-19-4}, and~\ref{sect-19-5};
\end{Problem}

\subsection{Pointwise asymptotics}
\label{sect-19-7-6-2}

I strongly believe that the following problem is both extremely challenging and interesting:

\begin{Problem}\label{Problem-19-7-17}
For operators of the type considered in this Chapter (namely, Sections~\ref{sect-19-2}--\ref{sect-19-6}) construct theory similar to theory of Chapter~\ref{book_new-sect-16} of \cite{futurebook}.
\end{Problem}

(i) Let us discuss the pilot-model with constant $g^{jk}$, $F_{jk}$, and linear $V$. Then $A$ is a sum of $2\D$-operators and the $U(x,y,t)$ is a product of $2\D$-propagators which after rescaling $x\mapsto \mu x$, $t\mapsto \mu t$ are given by (\ref{book_new-16-1-9})--(\ref{book_new-16-1-10}) of \cite{futurebook} with $t$ replaced by $f_j t$ and $\alpha_j$ replaced by $f_j^{-1}\alpha_j$. Then we come to the oscillatory integral with respect to $t$ of (\ref{book_new-16-1-22}) of \cite{futurebook} type  (and we will need to multiply it later by $\mu^{-r}h^{-r}$) but with a phase function which is the sum of (\ref{book_new-16-1-23})-type (of \cite{futurebook}) expressions
\begin{equation}
\phi (t) = -t^2 \sum_j \alpha_j^2 \cot (f_j t) - t(\tau-\beta),\qquad  \beta=-\mu^{-2}\sum_j \alpha_j^2 f_j^{-1}
\label{19-7-54}
\end{equation}
and with factor $\mu^2 \csc (t)$ replaced by $\mu^{2r} \csc(f_1t)\cdots \csc (f_rt)$. Now analysis of the stationary points becomes really difficult. The stationary phase equation 
\begin{gather}
-2t \sum_j \alpha_j^2 \cot (f_j t) +t^2 \sum_j \alpha_j^2 f_j \csc^2 (f_j t)  - (\tau-\beta)=0\label{19-7-55}\\
\shortintertext{is equivalent to}
-2t  \cot (f_j t) +t^2 \alpha_j^2 f_j \csc^2 (f_j t)  -\tau_j=0,
\label{19-7-56}\\
\sum \tau_j =\tau-\beta.
\label{19-7-57}
\end{gather}
This system is not easy to handle for large $t$ if $f_j$ are not commensurable.

\medskip\noindent
(ii) However in the case when $f_1=\ldots=f_r=1$ which is supposedly the worst case scenario we get almost (\ref{book_new-16-1-23}) of \cite{futurebook} exactly 
\begin{equation}
\phi (t) = - t^2 \cot ( t) +t \mu^{-2}\alpha^2 - t\tau
\label{19-7-58}
\end{equation}
but factor $\mu^{2r}(\csc(t))^r$ with $r\ge 2$ rather than $r=1$ will be a game-changer in all Tauberian estimates as $\mu \le h^{-1}$. 

\begin{Problem}\label{Problem-19-7-18}
(i) Prove that for the pilot-model with $f_1=\ldots=f_r$
\begin{equation}
\R^\T \le C\left\{\begin{aligned}
&\mu^{-1}h^{1-d}\qquad &&\text{as\ \ } \mu_0\le \mu \le h^{-\frac{1}{2}},\\
&\mu h^{2-d}\qquad &&\text{as\ \ }h^{-\frac{1}{2}}\le \mu \le h^{-1},\\
&\mu^{r-\frac{1}{2}}h^{-r+\frac{1}{2}} &&\text{as\ \ } \mu \ge h^{-1};
\label{19-7-59}
\end{aligned}\right.
\end{equation}
(ii) Prove the same in the case of commensurable $f_1,\ldots,f_r$;

\medskip\noindent
(iii) Investigate the case of non-commensurable $f_1,\ldots,f_r$.
\end{Problem}

\begin{Problem}\label{Problem-19-7-19}
(i) Calculate Tauberian expressions using stationary phase methods as $\mu \le h^{\delta -\frac{1}{2}}$.

\medskip\noindent
(ii) Investigate the matching cases when $V$ is non-linear.
\end{Problem}

\bibliographystyle{alpha}

\providecommand{\bysame}{\leavevmode\hbox to3em{\hrulefill}\thinspace}

\vglue .06truein

\begin{tabular}{rrl}
&{\hskip 200 pt} &Department of Mathematics,\cr
&&University of Toronto,\cr
&&40, St.George Str.,\cr
&&Toronto, Ontario M5S 2E4\cr
&&Canada\cr
&&ivrii@math.toronto.edu\cr
&&Fax: (416)978-4107\cr
\end{tabular}

\end{document}